\documentclass[journal,twocolumn]{IEEEtran}

\usepackage[english]{babel}
\usepackage[latin1]{inputenc}
\usepackage[T1]{fontenc}
\usepackage{url}
\usepackage{graphicx}
\usepackage{subfigure}
\usepackage{color}
\usepackage{amsmath,amsthm,amsfonts,amssymb}
\usepackage{colortbl}                   
\usepackage{url}
\usepackage{booktabs}
\usepackage{epsfig}
\usepackage{pstricks}
\usepackage{pst-node}
\usepackage{ifthen}
\usepackage{algorithm}
\usepackage{algpseudocode}

\newboolean{draft}
\newcommand{\isdraft}[2]{\ifthenelse{\boolean{draft}}{#1}{#2}}

\isdraft{\usepackage{setspace}}{}

\theoremstyle{definition}

\theoremstyle{plain}
\newtheorem{Theorem}{Theorem}

\newtheorem{Assumption}{Assumption}

\newcommand{\mypar}[1]{{\bf #1.}}

\isdraft{\renewcommand{\baselinestretch}{1.406}}{}

\title{Distributed Basis Pursuit}

\author{Jo\~ao F.~C.~Mota, Jo\~ao M.~F.~Xavier, Pedro M.~Q.~Aguiar, and~Markus P\"uschel
\IEEEcompsocitemizethanks{
  \IEEEcompsocthanksitem Copyright (c) 2011 IEEE. Personal use of this material is permitted. However, permission to use this material for any other purposes must be obtained from the IEEE by sending a request to pubs-permissions@ieee.org.
  \IEEEcompsocthanksitem Jo\~ao M.~F.~Xavier, Pedro M.~Q.~Aguiar, and Jo\~ao F.~C.~Mota are with Instituto de Sistemas e Rob\'otica (ISR), Instituto Superior T\'ecnico (IST), Technical University of Lisbon, Portugal.
  \IEEEcompsocthanksitem Jo\~ao F.~C.~Mota is also with the Department of Electrical and Computer Engineering at Carnegie Mellon University, USA.
  \IEEEcompsocthanksitem Markus P{\"u}schel is with the Department of
Computer Science at ETH Zurich, Switzerland.
}
\thanks{
  This work was supported by the FCT grant CMU-PT/SIA/0026/2009, PTDC/EEA-ACR/73749/2006 and SFRH/BD/33520/2008 (through the Carnegie Mellon/Portugal Program managed by ICTI) from Funda\c{c}\~{a}o para a Ci\^{e}ncia e Tecnologia and also by ISR/IST plurianual funding (POSC program, FEDER). This work was also supported by NSF through award 0634967.
  }  
}

\begin{document}

\maketitle

\begin{abstract}

We propose a distributed algorithm for solving the optimization
problem Basis Pursuit (BP). BP finds the
least~\isdraft{}{\boldmath}{$\ell_1$}-norm solution of the
underdetermined linear system~\isdraft{}{\boldmath}{$Ax = b$} and is
used, for example, in compressed sensing for reconstruction.  Our algorithm solves BP on a distributed platform such as a
sensor network, and is designed to minimize the communication between
nodes. The algorithm only requires the network to be connected, has no
notion of a central processing node, and no node has access to the
entire matrix~\isdraft{}{\boldmath}{$A$} at any time. We consider two scenarios in which
either the columns or the rows of~\isdraft{}{\boldmath}{$A$} are
distributed among the compute nodes. Our algorithm, named D-ADMM, is a decentralized implementation of the alternating direction method of multipliers. We show through numerical
simulation that our algorithm requires considerably less communications between
the nodes than the state-of-the-art algorithms.
\end{abstract}


\isdraft{\pagebreak}{}

\section{Introduction}
\label{Sec:intro}

    \emph{Basis Pursuit} (BP) is the convex optimization problem~\cite{AtomicDecompBP}
    \begin{equation}\tag{BP}\label{bp}
        \begin{array}[t]{ll}
            \textrm{minimize} & \|x\|_1 \\
            \textrm{subject to} & Ax = b, \\
        \end{array}
    \end{equation}
    where the optimization variable is~$x \in \mathbb{R}^{n}$, $\|x\|_1 =|x_1| + \cdots + |x_{n}|$ is the~$\ell_1$ norm of the vector~$x$, and~$A \in \mathbb{R}^{m \times {n}}$ is a matrix with more columns
    than rows: $m < {n}$. In words, BP seeks the ``smallest'' (in the $\ell_1$ norm sense) solution of the underdetermined linear system~$Ax = b$. To make sure that~$Ax = b$ has at least one solution, we require the following.
    \begin{Assumption}\label{Ass:Fullrank}
        $A$ is full rank.
    \end{Assumption}
    BP has recently attracted attention due to its ability to find the sparsest solution of a linear system under certain conditions (see \cite{SiamReview,CandesTerrTaoDecodingLP}). In particular, BP is a convex relaxation of the combinatorial and nonconvex problem obtained by replacing the $\ell_1$ norm in~\eqref{bp} by the~$\ell_0$ pseudonorm $\|x\|_0$, which counts the number of nonzero elements of~$x$. Note that the linear system~$Ax = b$ has a unique $k$-sparse solution, i.e., a solution whose~$\ell_0$ norm is~$k$, if every set of~$2k$ columns of~$A$ is linearly independent.

    BP belongs to a set of optimization problems that has applications in many areas of engineering. Examples include signal and image denoising and restoration~\cite{AtomicDecompBP,SiamReview}, compression, fitting and approximation of functions~\cite{Boyd:ConvexOpti}, channel estimation and coding~\cite{CandesTerrTaoDecodingLP} and compressed sensing~\cite{RobustUncertaintyPrinciples,DonohoCompressed} (for more applications see for example \cite{Tropp06,SiamReview} and the references therein). In particular, in the recent field of compressed sensing, BP plays a key role in the reconstruction of a signal.

    Notice that Assumption~\ref{Ass:Fullrank} holds with probability one if the entries of~$A$ are independent and identically distributed (i.i.d.) and drawn from some (non-degenerate) probability
    distribution, as commonly seen in compressed sensing~\cite{RobustUncertaintyPrinciples}. Also in compressed sensing, there are several strategies to deal with noisy observations, i.e., when the observation vector~$b$ is corrupted with noise. These include solving variations of~\eqref{bp}, namely BPDN~\cite{AtomicDecompBP} and LASSO~\cite{IntroCompressiveSamplingCandesWakin}.

\mypar{Problem statement and contribution}
    Consider a network (e.g., a sensor network) with~$P$ compute nodes, and partition the matrix~$A$ into~$P$ blocks. Our goal is to solve BP in a distributed way. By distributed we mean that there is no notion of a central processing node and that the $p$th node has only access to the block $A_p$ of $A$ during the execution.

    We partition~$A$ into blocks in two different ways, which we call \emph{row partition} and \emph{column partition}, visualized in Figure~\ref{Fig:IntroPartitionOfA}. In the row partition, the block~$A_p$ contains~$m_p$ rows of~$A$, which implies $m_1 + \cdots + m_P = m$. In the column partition, $A_p$ contains~$n_p$ columns of~$A$, which implies $n_1 + \cdots + n_P =~n$.

    \isdraft{
		\begin{figure}
    \centering
    \begin{pspicture}(8,2.5)
    \rput(1.8,1){
      $
      \begin{bmatrix}
        &  & \phantom{aaaaaaaaaaa} &  &  \\
        &  &  &  & \\
        &  &  &  & \\
       \isdraft{}{ &  &  &  &  \\}
      \end{bmatrix}
      $
    }
    \def\blockmatrix{
      \psframe*[linecolor=black!15!white,fillstyle=solid](0,0)(3.4,0.4667)
    }
    \rput[bl](0.1,0.2){\blockmatrix}
    \rput[bl](0.1,1.334){\blockmatrix}

    \rput(1.8,1.5667){$A_1$}
    \rput(1.8,1.1){$\vdots$}
    \rput(1.8,0.4334){$A_P$}
    
    \rput(1.8,2.23){Row Partition}
    \end{pspicture}
    \hspace{-7cm}
		\begin{pspicture}(8,2.5)
    \rput(6.2,1){
        $
        \begin{bmatrix}
           &  & \phantom{aaaaaaaaaaa} &  &  \\
           &  &  &  &  \\
           &  &  &  &  \\
           \isdraft{}{&  &  &  &  \\}
        \end{bmatrix}
        $
    }
    
    \def\blockmatrix{
        \psframe*[linecolor=black!15!white,fillstyle=solid](0,0)(0.775,1.6)
    }

    \rput[bl](4.5,0.2){\blockmatrix}
    \rput[bl](5.3750,0.2){\blockmatrix}
    \rput[bl](7.1250,0.2){\blockmatrix}

    \rput(4.8875,1){$A_1$}
    \rput(5.7625,1){$A_2$}
    \rput(6.6375,1){$\cdots$}
    \rput(7.5125,1){$A_P$}

    \rput(6.2,2.23){Column Partition}

    \end{pspicture}
    \isdraft{\vspace{-0.3cm}}{}
    \caption{Row partition and column partition of~$A$ into~$P$ blocks. We assume there are~$P$ nodes and the $p$th node stores~$A_p$. In the row partition a block is a set of rows, while in the column partition a block is a set of columns.}
    \label{Fig:IntroPartitionOfA}
    \isdraft{\vspace{-0.5cm}}{}
    \end{figure}
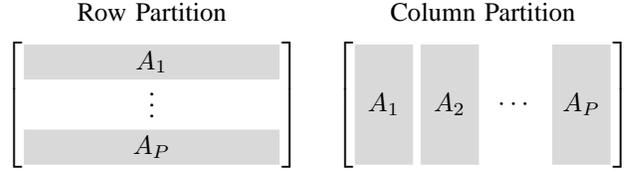
		}{
    \begin{figure}
    \centering
    \begin{pspicture}(8,2.5)
    \rput(1.8,1){
      $
      \begin{bmatrix}
        &  & \phantom{aaaaaaaaaaa} &  &  \\
        &  &  &  & \\
        &  &  &  & \\
       \isdraft{}{ &  &  &  &  \\}
      \end{bmatrix}
      $
    }
    \def\blockmatrix{
      \psframe*[linecolor=black!15!white,fillstyle=solid](0,0)(3.4,0.4667)
    }
    \rput[bl](0.1,0.2){\blockmatrix}
    \rput[bl](0.1,1.334){\blockmatrix}

    \rput(1.8,1.5667){$A_1$}
    \rput(1.8,1.1){$\vdots$}
    \rput(1.8,0.4334){$A_P$}

    \rput(6.2,1){
        $
        \begin{bmatrix}
           &  & \phantom{aaaaaaaaaaa} &  &  \\
           &  &  &  &  \\
           &  &  &  &  \\
           \isdraft{}{&  &  &  &  \\}
        \end{bmatrix}
        $
    }
    \def\blockmatrix{
        \psframe*[linecolor=black!15!white,fillstyle=solid](0,0)(0.775,1.6)
    }

    \rput[bl](4.5,0.2){\blockmatrix}
    \rput[bl](5.3750,0.2){\blockmatrix}
    \rput[bl](7.1250,0.2){\blockmatrix}

    \rput(4.8875,1){$A_1$}
    \rput(5.7625,1){$A_2$}
    \rput(6.6375,1){$\cdots$}
    \rput(7.5125,1){$A_P$}

    \rput(1.8,2.23){Row Partition}
    \rput(6.2,2.23){Column Partition}

    \end{pspicture}
    \isdraft{\vspace{-0.3cm}}{}
    \caption{Row partition and column partition of~$A$ into~$P$ blocks. We assume there are~$P$ nodes and the $p$th node stores~$A_p$. In the row partition a block is a set of rows, while in the column partition a block is a set of columns.}
    \label{Fig:IntroPartitionOfA}
    \isdraft{\vspace{-0.5cm}}{}
    \end{figure}
		}

    In summary: \emph{given a network, we solve BP in a distributed way, either in the row partition or in the column partition.}

    For the network we only require:
    \begin{Assumption}\label{Ass:Network}
        The given network is connected and static.
    \end{Assumption}
    Connected means that for any two nodes there is a path connecting them. Static means that the network topology does not change over time.

We propose an algorithm to solve this problem and show through extensive
simulations that it improves over previous work (discussed
below), by reducing the total number of communications to achieve a
given solution accuracy. The number of communications in distributed algorithms is an important measure of performance. For example, it is well known that communicating with the neighboring nodes is the most energy-consuming task for a wireless sensor~\cite{Akyildiz02-WirelessNetworksASurvey}; as a consequence, many energy-aware algorithms and protocols for wireless sensor networks have been proposed~\cite{Krishnamachari05-NetworkingWirelessSensors}. An energy-aware algorithm minimizes the communications and/or allows the nodes to become idle for some time instants. On distributed supercomputing platforms, on the other hand, computation time is the main concern. In this case, the computational bottleneck is again the communication between the nodes, and thus algorithms requiring less communications have the potential of being faster. 

Before we discuss related work, we provide possible applications of
our algorithm.

\mypar{Application: row partition}
    Given a network of~$P$ interconnected sensors, we try to capture an ultra-wide band but spectrally sparse signal, represented in vector form as~$x \in \mathbb{R}^n$. For simplicity, we assume the $p$th sensor only stores one row~$r_p^\top$ of~$A$, i.e., $m = P$. Each sensor only captures some time samples at a rate far below the Nyquist rate, using for example a random demodulator~\cite{BeyondNyquist,IntroCompressiveSamplingCandesWakin}. One can represent each measurement as the number~$b_p$.  Under certain conditions (\cite{RobustUncertaintyPrinciples,DonohoCompressed,L0L1NearOptimalUniversalEncodingStrategies}), it is possible to recover~$x$ by solving~\eqref{bp} with~$A = [r_1 \cdots r_P]^\top$ and~$b = [b_1 \cdots b_P]^\top$. Further details about the matrix~$A$ and the vector~$b$ can be found in~\cite{IntroCompressiveSamplingCandesWakin}. Since each vector~$r_p$ is associated with a sensor, this corresponds to our row partition case.  This scenario applies, for example, to sparse event detection in wireless networks~\cite{SparseEventDetection}, and to distributed target localization in sensor networks~\cite{DistributedTargetLoc}.

\mypar{Application: column partition}
    The work~\cite{RombergSeismic} introduces a method of speeding up seismic forward modeling in geological applications. The goal is to find the Green's functions of some model of a portion of the earth's surface. Given a set of sources and a set of receivers, from the knowledge of both the emitted and the received signals, the Green's function of the model, represented by~$x$, has to be found. The authors of~\cite{RombergSeismic} propose to solve this problem when all sources emit at the same time and the receivers capture a linear superposition of all signals. The approach is then to solve BP, where a set of columns of~$A$ is associated with a source. Note that a distributed solution makes sense because the sources are physically far apart.

    As another example for the column partition, we interpret BP as finding a sparse representation of a given signal~$b$ with respect to a dictionary of atomic signals (columns of~$A$). It is common to assume that the dictionary (the matrix~$A$) contains several families of functions, e.g., Fourier, DCT, wavelets, to become overcomplete. Suppose that we are given~$P$ processors, each of which is tuned to perform computations for a certain family of functions. In this case, solving BP in a column partition framework would arise naturally.

\mypar{Algorithms for solving BP and related work} Since BP can be
recast as a linear program (LP)~\cite{Boyd:ConvexOpti}, any algorithm
that solves LPs can also solve BP. Among the many algorithms solving
LPs~\cite{BertsimasLinearOptim}, most cannot be
readily adapted to our distributed scenario. For example, the
(distributed) simplex
algorithm~\cite{Hypercube,DistributedLinearProgDataMining}
can solve LPs only in complete networks, i.e., those with a link between
any pair of nodes. In this paper, we aim to solve BP
for every connected network topology.

In recent years, some approaches have been proposed for solving
general optimization problems, including BP, in distributed networks.
For example, \cite{MultiAgent} proposes a method based on
subgradient algorithms, but these are known to converge very slowly.
Other approaches to distributed optimization combine the method of
multipliers (MM) with the nonlinear Gauss-Seidel (NGS) method or with
Jacobi algorithms~\cite{Bertsekas:Parallel}. For example,
\cite{Ruszczynski} uses MM together with a Jacobi-type algorithm named
diagonal quadratic approximation (DQA) to solve, in a distributed way,
convex problems constrained by linear equations. Using a suitable
reformulation of~\eqref{bp}, this method can be applied to our problem
statement. In~\cite{Spars09} we analyzed how well MM together with NGS
solves BP in the row partition scenario; and in~\cite{Mota11ICASSP} we
used a fast gradient algorithm in both loops. The algorithm we propose
here has just one loop and requires considerably fewer iterations to
converge than all the previous approaches.

Fast algorithms solving BP in a non-distributed way include
spgl1~\cite{spgl1}, fpc~\cite{fpc}, LARS~\cite{LARS}, C-SALSA~\cite{CSALSA}, and NESTA~\cite{NESTA}. These are faster than
distributed algorithms but require that~$A$ and~$b$ are available at
the same location. In contrast, a distributed algorithm can solve
problems that can only fit into the combined memory of all the nodes.

The work~\cite{Giannakis} is closest related to ours. It solves the
Basis Pursuit Denoising (BPDN)~\cite{AtomicDecompBP} (a noise-robust
version of BP), which also produces sparse solutions of linear
systems. The algorithm is called D-Lasso and can be adapted to solve our problem. Our simulations show that the algorithm we propose requires systematically less communications than D-Lasso. 

Our algorithm is based on the alternating direction method of multipliers (ADMM). The work~\cite{BoydADMM} also
uses ADMM in a distributed scenario, but is only applicable to networks where all the nodes connect to a central node. Our algorithm, in contrast,
is designed for decentralized scenarios (no central node) and applies to any connected network.

Our type of matrix partitioning has been considered
before in the context of distributed algorithms for
linear
programs~\cite{Hypercube,DistributedLinearProgDataMining}
and in regression of distributed
data~\cite{NedicNewClassDistributedOptimization}.

\section{Row Partition}
\label{Sec:RowPartition}

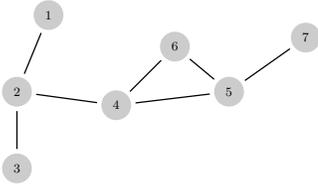
\begin{figure}[h]
        \centering
        \psscalebox{0.60}{
        \begin{pspicture}(7,4)
            \pscircle*[linecolor=black!20!white](1,3.7){0.33}
            \pscircle*[linecolor=black!20!white](0.3,2){0.33}
            \pscircle*[linecolor=black!20!white](0.3,0.3){0.33}
            \pscircle*[linecolor=black!20!white](2.5,1.7){0.33}
            \pscircle*[linecolor=black!20!white](5,2){0.33}
            \pscircle*[linecolor=black!20!white](3.8,3){0.33}
            \pscircle*[linecolor=black!20!white](6.7,3.2){0.33}

            \rput(1,3.7){\footnotesize $1$}
            \rput(0.3,2){\footnotesize $2$}
            \rput(0.3,0.3){\footnotesize $3$}
            \rput(2.5,1.7){\footnotesize $4$}
            \rput(5,2){\footnotesize $5$}
            \rput(3.8,3){\footnotesize $6$}
            \rput(6.7,3.2){\footnotesize $7$}

            \psline(0.4637,2.3976)(0.8363,3.3024)
            \psline(0.3000,1.5700)(0.3000,0.7300)
            \psline(0.7261,1.9419)(2.0739,1.7581)
            \psline(3.4959,2.6959)(2.8041,2.0041)
            \psline(4.5731,1.9488)(2.9269,1.7512)
            \psline(4.6697,2.2753)(4.1303,2.7247)
            \psline(5.3513,2.2480)(6.3487,2.9520)

        \end{pspicture}
        }
        \isdraft{\vspace{-0.1cm}}{}
        \caption{
                 Example of a connected network with~$P=E=7$.
                 The set of edges is~$\mathcal{E} = \{ \{1,2\}, \{2,3\}, \{2,4\}, \{4,5\}, \{4,6\}, \{5,6\}, \{5,7\} \}$. 
                }
        \label{Fig:ArbitraryNetwork}
    \isdraft{\vspace{-0.3cm}}{}
    \end{figure}

In this section we partition the matrix~$A$ by rows:
    $$
    \phantom{A =}
    \begin{pspicture}(4,2)
        \rput(2,1){
          $
          \begin{bmatrix}
            &  & \phantom{aaaaaaaaaaa} &  &  \\
            &  &  &  &  \\
            \isdraft{}{&  &  &  &  \\}
            \isdraft{}{&  &  &  &  \\}
          \end{bmatrix}
          $
        }
        \def\blockmatrix{
          \psframe*[linecolor=black!15!white,fillstyle=solid](0,0)(3.4,0.4667)
        }
        \rput[bl](0.3,0.2){\blockmatrix}
        \rput[bl](0.3,1.334){\blockmatrix}

        \rput(2,1.5667){$A_1$}
        \rput(2,1.1){$\vdots$}
        \rput(2,0.4334){$A_P$}

        \isdraft{
        \rput(-0.6,1){$A = $}
        \rput(4.15,0.9){,}
        }{
        \rput(-0.4,1){$A = $}
        \rput(4.05,0.9){,}
        }
    \end{pspicture}
    \isdraft{\vspace{-0.2cm}}{}
    $$
where each block $A_p \in {\mathbb R}^{m_p\times n}$ contains a subset
of rows of~$A$ such that $m_1 + \cdots + m_P = m$. The vector~$b$ is
partitioned similarly: $b = [b_1^\top \cdots b_P^\top]^\top$. We
assume that~$A_p$ and $b_p$ are available only at the $p$th node of a
connected network with~$P$ compute nodes. We model the network as an
undirected graph~$\mathcal{G} = ({\mathcal V}, {\mathcal E} )$, where
${\mathcal V} = \{ 1, 2, \ldots, P \}$ is the set of nodes
and~$\mathcal{E} \subset {\mathcal V} \times {\mathcal V}$ is the set of edges. We represent the edge connecting nodes~$i$ and~$j$ by~$\{i,j\}$ or~$\{j,i\}$; $E$ is the total number of
edges. See Figure~\ref{Fig:ArbitraryNetwork} for an
example graph. If~$\{i,j\}$ is an edge, then
node~$i$ and node~$j$ can exchange
messages with each other. The
set of neighbors of node $p$ is written
as~$\mathcal{N}_p$, and its degree is $D_p = |\mathcal{N}_p|$.

\mypar{Graph coloring} We assume that a proper coloring~$\mathcal{C} =
\{1,\ldots,C\}$ of the graph is available. This means that each node
is labeled with a number $c \in \mathcal{C}$, which we call color,
such that no adjacent nodes (i.e., neighbors) have the same color. The
minimum number of colors required for a proper coloring of a
graph~$\mathcal{G}$ is its chromatic number~$\chi(\mathcal{G})$.
Coloring a graph with~$\chi(\mathcal{G})$
colors or just computing~$\chi(\mathcal{G})$ is NP-hard
for~$\chi(\mathcal{G}) >
2$~\cite{GareyComputersIntractability}. Several distributed algorithms
for coloring a graph exist~\cite{Kuhn06ComplexityDistributedGraphColoring,Leith06DistributedLearningAlgorithms,Duffy09ComplexityAnalysisDecentralizedGraphColoring,Linial92LocalityDistributedGraphAlgorithms}. For
example, \cite{Kuhn06ComplexityDistributedGraphColoring} determines
a coloring with~$O(D_{\max})$ colors, where~$D_{\max}
= \max_p D_p$, using~$O(D_{\max}/\log^2(D_{\max}) +
\log^\star(P))$ iterations. If more colors are allowed, for example
$O(D_{\max}^2)$, then ~$O(\log^\star( P))$
iterations suffice~\cite{Linial92LocalityDistributedGraphAlgorithms}.
In this paper we assume that a proper coloring~$\mathcal{C}$ with~$C$ colors is
given.

\mypar{Problem reformulation}
To solve BP in a distributed way we first rewrite~\eqref{bp}
to make the row partition explicit:
    \begin{equation}\label{Eq:RPReformulation1}
        \begin{array}{ll}
          \textrm{minimize} & \|x\|_1 \\
          \textrm{subject to} & A_p x = b_p\,,\quad p=1,\ldots,P\,.
        \end{array}
    \end{equation}
The variable~$x$ is coupling the problem. To decouple,
we replace~$x$ with~$P$
copies of~$x$. The $p$th copy is denoted with~$x_p$. To ensure the
necessary global consistency condition~$x_1 = x_2 = \cdots = x_P$, we
enforce the equivalent (since the network is connected)
constraint~$x_i = x_j$ for each edge~$\{i,j\}$ of the
network:
    \begin{equation}\label{Eq:ReformEdges}
        \begin{array}{ll}
          \textrm{minimize} & \frac{1}{P}\sum_{p=1}^P\|x_p\|_1 \\
          \textrm{subject to} & A_p x_p = b_p\,,\quad p=1,\ldots,P \\
          & x_i = x_j\,,\quad \{i,j\} \in \mathcal{E}\,.
        \end{array}
    \end{equation}
The optimization variable is $\bar{x}:=(x_1,\ldots,x_P) \in
(\mathbb{R}^n)^P$. Note that~\eqref{Eq:ReformEdges} can be written more
compactly as
    \begin{equation}\label{Eq:ReformEdgesB}
        \begin{array}{ll}
          \textrm{minimize} & \frac{1}{P}\sum_{p=1}^P\|x_p\|_1 \\
          \textrm{subject to} & A_p x_p = b_p\,,\quad p=1,\ldots,P \\
          & (B^\top \otimes I_n) \bar{x} = 0\,,
        \end{array}
    \end{equation}
    where~$I_n$ is the~$n\times n$ identity matrix, and $\otimes$ is the Kronecker product. The matrix~$B$ is the $P \times E$ node-arc incidence matrix of the graph:
    each edge~$\{i,j\} \in \mathcal{E}$ corresponds to a column in~$B$ with the $i$th and $j$th
    entries equal to~$1$ and~$-1$, respectively.

\mypar{Algorithm for bipartite graphs} We first consider a simple
case: $\mathcal{G}$ is bipartite and hence~$\chi (\mathcal{G}) =
2$. The generalization to any connected graph will be straightforward.
Bipartite graphs include trees and grid graphs.

    Without loss of generality, assume nodes~$1$ to~$c$ have color~$1$ and the remaining have color~$2$. Then, \eqref{Eq:ReformEdgesB} can be written as
    \begin{equation}\label{Eq:ReformEdgesBPartit}
        \begin{array}{ll}
          \textrm{minimize} & \frac{1}{P}\sum_{p=1}^c\|x_p\|_1 + \frac{1}{P}\sum_{p=c+1}^P\|x_p\|_1 \\
          \textrm{subject to} & A_p x_p = b_p\,,\quad p=1,\ldots,P \\
          & (B_1^\top \otimes I_n) \bar{x}_1 + (B_2^\top \otimes I_n) \bar{x}_2 = 0\,,
        \end{array}
    \end{equation}
    where~$\bar{x} = (\bar{x}_1, \bar{x}_2) \in (\mathbb{R}^{n})^c \times (\mathbb{R}^{n})^{P-c}$ and~$B$ is partitioned as
    $
        B = \begin{bmatrix}
              B_1^\top &
              B_2^\top
            \end{bmatrix}^\top
    $.
We propose the alternating direction method of multipliers (ADMM, briefly described in appendix~\ref{App:ADMM}) to
solve~\eqref{Eq:ReformEdgesBPartit}. The augmented Lagrangian
of~\eqref{Eq:ReformEdgesBPartit}, dualizing only the last
constraint, is
    \begin{multline}\label{Eq:RPBipAugmentedLagrangian}
        L(\bar{x}_1,\bar{x}_2;\lambda) = \frac{1}{P}\sum_{p \in \mathcal{C}_1} \|x_p\|_1 + \frac{1}{P}\sum_{p \in \mathcal{C}_2} \|x_p\|_1 + \phi_1(\bar{x}_1,\lambda)
        \isdraft{}{\\}+ \phi_2(\bar{x}_2,\lambda) + \rho \bar{x}_1^\top (B_1 B_2^\top \otimes I_n) \bar{x}_2\,,
    \end{multline}
    where~$\mathcal{C}_1 = \{1,\ldots,c\}$, $\mathcal{C}_2 = \{c+1,\ldots,P\}$, and
    \begin{align*}
          \phi_i(\bar{x}_i,\lambda)
        \isdraft{}{&}=
          \lambda^\top (B_i^\top \otimes I_n)\bar{x}_i + \frac{\rho}{2}\|(B_i^\top \otimes I_n)\bar{x}_i\|^2
        \isdraft{}{\notag
        \\
        &}=
          ((B_i \otimes I_n) \lambda)^\top \bar{x}_i + \frac{\rho}{2} \bar{x}_i^\top(B_i B_i^\top \otimes I_n)\bar{x}_i\,,
    \end{align*}
    for~$i=1,2$. Note that, since nodes in each~$\mathcal{C}_i$ are not neighbors between themselves, $B_i B_i^\top$ is diagonal (with~$D_p$ in the $p$th diagonal entry). Hence,
		\begin{equation}\label{Eq:RPphi}
			\phi_i(\bar{x}_i, \lambda) = \sum_{p \in \mathcal{C}_i} \Bigl(\gamma_p^\top x_p + \frac{\rho}{2}D_p\|x_p\|^2 \Bigr)\,,\quad i = 1,2 \,,
		\end{equation}
		where~$\gamma_p := \sum_{j \in \mathcal{N}_p} \textrm{sign} (j-p) \lambda_{\{p,j\}}$ and~$\textrm{sign}(w)$ gives~$1$ if~$w \geq 0$ and~$-1$ otherwise. We decomposed the dual variable~$\lambda$ into $(\ldots,\lambda_{\{i,j\}},\ldots)$, where~$\lambda_{\{i,j\}} = \lambda_{\{j,i\}}$ is associated with the constraint~$x_i = x_j$. 		
    
    Equations~\eqref{Eq:RPBipAugmentedLagrangian} and~\eqref{Eq:RPphi} show that minimizing~$L(\bar{x}_1,\bar{x}_2;\lambda)$ with respect to (w.r.t.) $\bar{x}_1$ yields~$c$ optimization problems that can be executed in parallel; similarly, minimizing it w.r.t.\ $\bar{x}_2$ yields~$P-c$ parallel optimization problems. Algorithm~\ref{Alg:ADMMBipartite} shows the application of ADMM to our problem. We name our algorithm D-ADMM, after Distributed ADMM. 
    \begin{algorithm}[H]
    \caption{D-ADMM for bipartite graphs}
    \algrenewcommand\algorithmicrequire{\textbf{Initialization:}}
    \label{Alg:ADMMBipartite}
    \begin{algorithmic}[1]
    \small
    \Require for all~$p \in \mathcal{V}$, set $\gamma_{p}^{(1)} = x_p^{(1)} = 0$ and $k=1$
    \Repeat
    \ForAll{$p \in \mathcal{C}_1$ [in parallel]}
    \label{SubAlg:ADMMBipLoop}            
    \State Set $v_p^{(k)} =  \gamma_p^{(k)} -\rho \sum_{j \in \mathcal{N}_p} x_j^{(k)}$ and find
    \label{SubAlg:ADMMBipOptimProb}
           $$
            x_p^{(k+1)} = \begin{array}[t]{cl}
                            \underset{x_p}{\textrm{argmin}} & \frac{1}{P}\|x_p\|_1 + {v_p^{(k)}}^\top x_p + \frac{D_p \rho}{2}\|x_p\|^2\\
                            \textrm{s.t.} & A_p x_p = b_p
                          \end{array}
           $$
    \State Send~$x_p^{(k+1)}$ to $\mathcal{N}_p$
    \label{SubAlg:ADMMBipCommunication}
    \EndFor
    \label{SubAlg:ADMMBipEndLoop}

    \State Repeat \ref{SubAlg:ADMMBipLoop}-\ref{SubAlg:ADMMBipEndLoop} for all~$p \in \mathcal{C}_2$, replacing $x_j^{(k)}$ by $x_j^{(k+1)}$

    \ForAll{$p \in \mathcal{C}_1 \cup \mathcal{C}_2$ [in parallel]} \vspace{0.15cm}
    \hfill
    
        $
            \gamma_p^{(k+1)} = \gamma_p^{(k)} + \rho \sum_{j \in \mathcal{N}_p} (x_p^{(k+1)} -  x_j^{(k+1)})
        $\vspace{0.15cm}
        \label{SubAlg:ADMMBipDualVarUp}
    \EndFor
    \State $k \gets k+1$
    \Until{some stopping criterion is met}
    \end{algorithmic}
    \end{algorithm}    
            
    The optimization problem in step~\ref{SubAlg:ADMMBipOptimProb} results from minimizing the augmented Lagrangian~$L(\bar{x}_1, \bar{x}_2;\lambda)$ w.r.t.\ $x_p$. To derive it, note that~\eqref{Eq:RPphi} enables us to rewrite~$L(\bar{x}_1, \bar{x}_2; \lambda)$ as
    \isdraft{
			$$
			L(\bar{x}_1, \bar{x}_2;\lambda) = \sum_{i=1}^2 \sum_{p \in \mathcal{C}_i} 
			\Bigl(
				\frac{1}{P}\|x_p\|_1 + \gamma_p^\top x_p + \frac{\rho}{2}D_p \|x_p\|^2
			\Bigr)
			\isdraft{}{\\}+
			\rho \,\bar{x}_1 (B_1 B_2^\top \otimes I_n) \bar{x}_2\,.
			$$
		}{
			\begin{multline*}
			L(\bar{x}_1, \bar{x}_2;\lambda) = \sum_{i=1}^2 \sum_{p \in \mathcal{C}_i} 
			\Bigl(
				\frac{1}{P}\|x_p\|_1 + \gamma_p^\top x_p + \frac{\rho}{2}D_p \|x_p\|^2
			\Bigr)
			\isdraft{}{\\}+
			\rho \,\bar{x}_1 (B_1 B_2^\top \otimes I_n) \bar{x}_2\,.
    \end{multline*}
    }
    The~$(ij)$th entry of~$B_1B_2^\top$ is~$-1$ if~$\{i,j\} \in \mathcal{E}$ and~$0$ otherwise. Therefore, $\rho \,\bar{x}_1 (B_1 B_2^\top \otimes I_n) \bar{x}_2 = -\rho\sum_{\{i,j\}\in \mathcal{E}}x_i^\top x_j$. Picking~$p \in \mathcal{C}_i$ for any~$i = 1,2$ and minimizing~$L(\bar{x}_1, \bar{x}_2;\lambda)$ w.r.t.~$x_p$ yields the optimization problem in step~\ref{SubAlg:ADMMBipOptimProb}. Appendix~\ref{App:OptimForEachNode} describes an efficient method for solving this problem.

    Algorithm~\ref{Alg:ADMMBipartite} shows that nodes with the same color operate in parallel, whereas nodes with different colors cannot. In other words, the nodes from~$\mathcal{C}_1$ have to wait for the computation of the nodes from~$\mathcal{C}_2$ and vice-versa. However, at the end of each iteration, every node will have communicated once (sending~$x_p^{(k+1)}$ and receiving~$x_j^{(k+1)}$) with all its neighbors.

    Regarding the dual variable~$\lambda$, its components do not appear explicitly in Algorithm~\ref{Alg:ADMMBipartite}. The reason is that node~$p$ only requires~$\gamma_p = \sum_{j \in \mathcal{N}_p} \textrm{sign} (j-p) \lambda_{\{p,j\}}$ for its optimization problem. According to the canonical form of ADMM, we have to update~$\lambda_{\{i,j\}}$, for each edge~$\{i,j\} \in \mathcal{E}$ as
    \begin{equation}\label{Eq:RPBipLambdaUpdate}
    	\lambda_{\{i,j\}}^{(k+1)} = \lambda_{\{i,j\}}^{(k)} + \rho \,\text{sign}(j-p)(x_i^{(k+1)} - x_j^{(k+1)})\,.
    \end{equation}
    Inserting~\eqref{Eq:RPBipLambdaUpdate} into the expression of~$\gamma_p$ we obtain the update of step~\ref{SubAlg:ADMMBipDualVarUp}. 
    
    The following theorem establishes the convergence of Algorithm~\ref{Alg:ADMMBipartite}.
    \begin{Theorem}\label{Thm:Bipartite}
        Assume the given graph is bipartite. Then, for all~$p$, the sequence~$\{x_p^{(k)}\}$ produced by Algorithm~\ref{Alg:ADMMBipartite} converges to a solution of~\eqref{bp}.
    \end{Theorem}
    \begin{proof}
        We have already seen that when the graph is bipartite~\eqref{bp} is equivalent to~\eqref{Eq:ReformEdgesBPartit}. We now show that~\eqref{Eq:ReformEdgesBPartit} satisfies the conditions of Theorem~\ref{Thm:ADMMConvergence} in appendix~\ref{App:ADMM}. Let~$f_i(\bar{x}_i) = (1/P)\sum_{c \in \mathcal{C}_i}\|x_p\|_1$, for~$i=1,2$. Clearly, $f_1$ and~$f_2$ are real-valued convex functions. Assumption~\ref{Ass:Fullrank} on the rank of the matrix~$A$ implies that~\eqref{bp}, and thus~\eqref{Eq:ReformEdgesBPartit}, is always solvable. Also, the non-dualized equations~$A_p x_p = b_p$  in~\eqref{Eq:ReformEdgesBPartit} define polyhedral sets. 

				Now we have to prove that the matrices~$B_1^\top \otimes I_n$ and~$B_2^\top \otimes I_n$ have full column rank, i.e., that~$B_1^\top$ and~$B_2^\top$ have full column rank. We have seen that~$B_1 B_1^\top$ and~$B_2 B_2^\top$ are diagonal matrices because the nodes within one class are not neighbors. Note that the $p$th entry of the diagonal of~$B_1 B_1^\top$ (or~$B_2 B_2^\top$) is the degree of the $p$th node. Due to Assumption~\ref{Ass:Network}, there are no isolated nodes and thus~$B_1 B_1^\top$  and~$B_2B_2^\top$ are full-rank. The result then follows because $\textrm{rank}\,(B B^\top) = \textrm{rank}\,(B^\top)$ for any matrix~$B$.
\end{proof}

Theorem~\ref{Thm:Bipartite} also shows that after
Algorithm~\ref{Alg:ADMMBipartite} terminates, every node will know a
solution~$x^\star$ of BP.

\mypar{Algorithm for general graphs} We now generalize Algorithm~\ref{Alg:ADMMBipartite} to arbitrary graphs with $\chi(\mathcal{G}) > 2$.
The generalization is straightforward, but we cannot guarantee
convergence as in Theorem~\ref{Thm:Bipartite}. However, in our extensive experiments,
shown later, the resulting algorithm never failed to converge.

Let~$\mathcal{G}$ be a graph with a proper coloring~$\mathcal{C}$ and let $C
= |\mathcal{C}|$ be the number of colors. Let~$\mathcal{C}_c$
be the set of nodes that have color~$c$, $c =1,\ldots,C$. Without loss of generality, suppose the nodes are numbered the following way: $\mathcal{C}_1 =
\{1,\ldots,|\mathcal{C}_1|\}$, $\mathcal{C}_2 =
\{|\mathcal{C}_1|+1,\ldots,|\mathcal{C}_1|+|\mathcal{C}_2|\}$, \ldots,
$\mathcal{C}_C = \{\sum_{c=1}^{C-1}|\mathcal{C}_c| +
1,\ldots,P\}$. This enables a partition of the matrix~$B$ as~$ B =
        \begin{bmatrix}
          B_1^\top &
          \cdots &
          B_C^\top
        \end{bmatrix}^\top\!\!,
    $
making~\eqref{Eq:ReformEdgesB} equivalent to
    \begin{equation}\label{Eq:RPGenProb}
        \begin{array}{ll}
          \textrm{minimize} & \frac{1}{P}\sum_{c =1}^C \sum_{p \in \mathcal{C}_c} \|x_p\|_1 \\
          \textrm{subject to} & A_p x_p = b_p \,,\quad p = 1,\ldots,P \\
                              & \sum_{c=1}^C (B_c^\top \otimes I_n) \bar{x}_c = 0\,,
        \end{array}
    \end{equation}
    where~$\bar{x} = (\bar{x}_1,\ldots,\bar{x}_C)$ is the variable, and~$\bar{x}_c \in (\mathbb{R}^n)^{|\mathcal{C}_c|}$ for~$c=1,\ldots, C$. From the proof of Theorem~\ref{Thm:Bipartite} we know that each matrix~$B_c$ has full row rank. Thus, we can apply the generalized ADMM to solve~\eqref{Eq:RPGenProb} (see Appendix~\ref{App:ADMM}). That leads to the following algorithm.
    \begin{algorithm}[H]
    \caption{D-ADMM for general graphs}
    \algrenewcommand\algorithmicrequire{\textbf{Initialization:}}
    \label{Alg:ADMMGeneral}
    \begin{algorithmic}[1]
    \small
    \Require for all~$p \in \mathcal{V}$, set $\gamma_{p}^{(1)} = x_p^{(1)} = 0$ and $k=1$    
    \Repeat
    \For{$c =1,\ldots,C$}
        \ForAll{$p \in \mathcal{C}_c$ [in parallel]}
            $$
                v_p^{(k)} = \gamma_p^{(k)}-
                \rho \sum_{\begin{subarray}{c}
                             j \in \mathcal{N}_p \\
                             j < p
                           \end{subarray}
                }x_j^{(k+1)} - \rho \sum_{\begin{subarray}{c}
                             j \in \mathcal{N}_p \\
                             j > p
                           \end{subarray}
                }x_j^{(k)}
            $$
        \State and find
            $$
            x_p^{(k+1)} = \begin{array}[t]{cl}
                            \underset{x_p}{\textrm{argmin}} & \frac{1}{P}\|x_p\|_1 + {v_p^{(k)}}^\top x_p + \frac{D_p \rho}{2}\|x_p\|^2\\
                            \textrm{s.t.} & A_p x_p = b_p
                          \end{array}
           $$
           \label{SubAlg:ADMMGenProb}
        \State Send~$x_p^{(k+1)}$ to $\mathcal{N}_p$
        \label{SubAlg:ADMMGenComm}
    \EndFor
    \EndFor

    \ForAll{$p =1,\ldots,P$ [in parallel]} \vspace{0.15cm}
    \hfill
    
        $
            \gamma_p^{(k+1)} = \gamma_p^{(k)} + \rho \sum_{j \in \mathcal{N}_p} (x_p^{(k+1)} -  x_j^{(k+1)})
        $\vspace{0.15cm}
        \label{SubAlg:ADMMGenDualVarUp}
    \EndFor
    \State $k \gets k+1$
    \Until{some stopping criterion is met}
    \end{algorithmic}
    \end{algorithm}

Algorithm~\ref{Alg:ADMMGeneral} is a straightforward generalization of
Algorithm~\ref{Alg:ADMMBipartite}. Now there are~$C$ classes of nodes
and all the nodes in one class ``work'' in parallel, but the classes
cannot work at the same time. Consequently, if we consider the time to
solve one instance of the problem in step~\ref{SubAlg:ADMMGenProb} as
one unit, one (outer) iteration in Algorithm~\ref{Alg:ADMMGeneral}
takes~$C$ units.

\begin{figure}[h]
        \centering
        \subfigure[Undirected]{
        \psscalebox{0.6609}{
        \begin{pspicture}(6,4)
            \rput(0.5,3.0){\rnode{N1}{\pscircle*[linecolor=black!60!white](0,0){0.33}}}
            \rput(3.0,3.0){\rnode{N2}{\pscircle*[linecolor=black!15!white](0,0){0.33}}}
            \rput(5.5,3.0){\rnode{N3}{\pscircle[fillstyle=vlines*,linecolor=black!50!white,hatchcolor=black!50!white](0,0){0.33}}}
            \rput(2.0,1.0){\rnode{N4}{\pscircle[fillstyle=vlines*,linecolor=black!50!white,hatchcolor=black!50!white](0,0){0.33}}}
            \rput(4.0,1.0){\rnode{N5}{\pscircle*[linecolor=black!60!white](0,0){0.33}}}

            \rput(0.5,3.0){\footnotesize \textcolor[rgb]{0.85,0.85,0.85}{$1$}}
            \rput(3.0,3.0){\footnotesize $3$}
            \rput(5.5,3.0){\footnotesize $5$}
            \rput(2.0,1.0){\footnotesize $4$}
            \rput(4.0,1.0){\footnotesize \textcolor[rgb]{0.85,0.85,0.85}{$2$}}

            \ncline[nodesep=0.4cm]{-}{N1}{N2}
            \ncarc[nodesep=0.4cm,arcangle=35,arrowsize=6pt]{-}{N1}{N3}
            \ncline[nodesep=0.4cm]{-}{N1}{N4}
            \ncline[nodesep=0.4cm]{-}{N2}{N3}
            \ncline[nodesep=0.4cm]{-}{N2}{N4}
            \ncline[nodesep=0.4cm]{-}{N2}{N5}
            \ncline[nodesep=0.4cm]{-}{N3}{N5}
            \ncline[nodesep=0.4cm]{-}{N4}{N5}

        \end{pspicture}
        }
        \label{SubFig:Undirected}
        }
        \isdraft{\hspace{4cm}}{\hfill}
        \subfigure[Directed]{
        \psscalebox{0.6609}{
        \begin{pspicture}(6,4)
            \rput(0.5,3.0){\rnode{N1}{\pscircle*[linecolor=black!60!white](0,0){0.33}}}
            \rput(3.0,3.0){\rnode{N2}{\pscircle*[linecolor=black!15!white](0,0){0.33}}}
            \rput(5.5,3.0){\rnode{N3}{\pscircle[fillstyle=vlines*,linecolor=black!50!white,hatchcolor=black!50!white](0,0){0.33}}}
            \rput(2.0,1.0){\rnode{N4}{\pscircle[fillstyle=vlines*,linecolor=black!50!white,hatchcolor=black!50!white](0,0){0.33}}}
            \rput(4.0,1.0){\rnode{N5}{\pscircle*[linecolor=black!60!white](0,0){0.33}}}

            \rput(0.5,3.0){\footnotesize \textcolor[rgb]{0.85,0.85,0.85}{$1$}}
            \rput(3.0,3.0){\footnotesize $3$}
            \rput(5.5,3.0){\footnotesize $5$}
            \rput(2.0,1.0){\footnotesize $4$}
            \rput(4.0,1.0){\footnotesize \textcolor[rgb]{0.85,0.85,0.85}{$2$}}

            \ncline[nodesep=0.4cm,arrowsize=6pt]{->}{N1}{N2}
            \ncarc[nodesep=0.4cm,arcangle=35,arrowsize=6pt]{->}{N1}{N3}
            \ncline[nodesep=0.4cm,arrowsize=6pt]{->}{N1}{N4}
            \ncline[nodesep=0.4cm,arrowsize=6pt]{->}{N2}{N3}
            \ncline[nodesep=0.4cm,arrowsize=6pt]{->}{N2}{N4}
            \ncline[nodesep=0.4cm,arrowsize=6pt]{<-}{N2}{N5}
            \ncline[nodesep=0.4cm,arrowsize=6pt]{<-}{N3}{N5}
            \ncline[nodesep=0.4cm,arrowsize=6pt]{<-}{N4}{N5}

        \end{pspicture}
        }
        \label{SubFig:Directed}
        }
        \caption{
            $(\textrm{a})$ undirected network with $\chi(\mathcal{G}) = 3$ and with classes~$\mathcal{C}_1 = \{1,2\}$, $\mathcal{C}_2 = \{3\}$, $\mathcal{C}_3 = \{4,5\}$; $(\textrm{b})$ directed graph constructed from~$(\textrm{a})$ by assigning a direction to each link: from smallest color node to the largest color node.
                }
        \label{Fig:UndirectedAndDirected}
    \isdraft{\vspace{-0.5cm}}{}
    \end{figure}
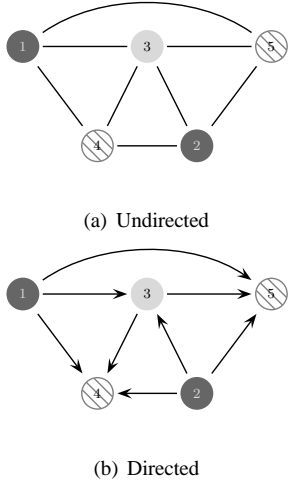

In the bipartite case the coordination between the nodes was straightforward: node~$p$ only works after it has received~$x_j$ from all its neighbors. Here, according to the canonical format of Algorithm~\ref{Alg:ADMMGeneral}, all the nodes in one class should work at the same time. Since these nodes are not neighbors, neither there is a central node to coordinate them, in practice node~$p$ works after having received ~$x_j^{(k+1)}$'s from all its neighbors of lower color. An alternative way to see this is to transform the undirected graph of the network into a directed graph, as shown in
Figure~\ref{Fig:UndirectedAndDirected}. The graph in Figure~\ref{Fig:UndirectedAndDirected}\subref{SubFig:Directed} is constructed from the graph
in Figure~\ref{Fig:UndirectedAndDirected}\subref{SubFig:Undirected} by assigning a direction to each edge~$\{i,j\}$: $i \rightarrow j$ if the color of~$i$ is smaller than the color of~$j$, and~$i \leftarrow j$
otherwise. Then, each node only starts working after having received the~$x_j$'s from all its inward links. In practice, this procedure can reduce the overall execution time since each node does not need to wait for its ``color time.'' As described in step~\ref{SubAlg:ADMMGenComm} (and in contrast to what Figure~\ref{Fig:UndirectedAndDirected}\subref{SubFig:Directed} may suggest), each node sends~$x_p^{k+1}$ to all its neighbors in each iteration. 

As stated earlier, we have no proof of convergence
for Algorithm~\ref{Alg:ADMMGeneral}, only practical evidence.

\section{Column Partition}
\label{Sec:ColumnPartition}

In this section, we adapt the algorithm for the row partition
to the column partition case:
    $$
    \phantom{A =}
    \begin{pspicture}(4,2)
        \rput(2,1){
          $
          \begin{bmatrix}
            &  & \phantom{aaaaaaaaaaa} &  &  \\
            &  &  &  &  \\
            \isdraft{}{&  &  &  &  \\}
            \isdraft{}{&  &  &  &  \\}
          \end{bmatrix}
          $
        }
        \def\blockmatrix{
            \psframe*[linecolor=black!15!white,fillstyle=solid](0,0)(0.775,1.6)
        }
        \rput[bl](0.3,0.2){\blockmatrix}
        \rput[bl](1.1750,0.2){\blockmatrix}
        \rput[bl](2.9250,0.2){\blockmatrix}
        \rput(0.6875,1){$A_1$}
        \rput(1.5625,1){$A_2$}
        \rput(2.4375,1){$\cdots$}
        \rput(3.3125,1){$A_P$}
        \isdraft{
        \rput(-0.6,1){$A = $}
        \rput(4.15,0.9){.}
        }{
        \rput(-0.4,1){$A = $}
        \rput(4.05,0.9){.}
        }
    \end{pspicture}
    \isdraft{\vspace{-0.2cm}}{}
    $$
    Each block $A_p \in \mathbb{R}^{m\times n_p}$ contains a subset of columns of~$A \in \mathbb{R}^{m\times n}$ such that~$n_1 + \cdots + n_P = n$. The block~$A_p$ is only available at the $p$th node of an arbitrary connected network, and the vector~$b \in \mathbb{R}^m$ is known by all the nodes.

\mypar{Duality: pros and cons}
    In section~\ref{Sec:RowPartition} we saw an algorithm that solves BP with a row partition. Here, we want to reutilize that algorithm for BP with a column partition. The first approach to that is to consider the dual problem of~\eqref{bp}:
    \begin{equation}\label{Eq:CPDualBP}
            \begin{array}[t]{ll}
                \textrm{minimize} &  b^\top \lambda \vspace{0.5ex}\\
                \textrm{subject to} & -1_n \leq A^\top \lambda \leq 1_n\,,

            \end{array}
    \end{equation}
where the dual variable is~$\lambda \in \mathbb{R}^m$, and~$1_{n}\in\mathbb{R}^{n}$ is
the vector of all ones. For a derivation
of~\eqref{Eq:CPDualBP}, see for
example~\cite[\S1.3.3]{JoaoMotaThesis}. The matrix~$A$ now
appears in the constraints of~\eqref{Eq:CPDualBP} as~$A^\top$, i.e.,
we can partition the constraint matrix in~\eqref{Eq:CPDualBP} by
rows. The problem is that there is no
straightforward way to recover a solution of~\eqref{bp} from a
solution of~\eqref{Eq:CPDualBP}. Hence we need an alternative
approach.

\mypar{Regularizing BP}
    Consider the following regularized version of~\eqref{bp}:
    \begin{equation}\label{Eq:ColBPRegularized}
        \begin{array}{ll}
          \textrm{minimize} & \|x\|_1 + \frac{\delta}{2}\|x\|^2\\
          \textrm{subject to} & Ax = b\,,
        \end{array}
    \end{equation}
    where~$\delta$ is a small positive number. While~\eqref{bp} may have multiple solutions, \eqref{Eq:ColBPRegularized} just has one, due to the strict convexity of its objective. When~$\delta$ is small enough, \eqref{Eq:ColBPRegularized} selects the least $\ell_2$-norm solution of~\eqref{bp}:
    \begin{Theorem}\label{Thm:Regularization}
        There exists~$\bar{\delta} > 0$ such that the solution of~\eqref{Eq:ColBPRegularized} is a solution of~\eqref{bp} for all~$0 < \delta < \bar{\delta}$.
    \end{Theorem}
    The proof of this theorem is based on exact regularization results for linear programming~\cite{Friedlander,MangasarianLPs}. To prove it, recast~\eqref{bp} as a linear program~\cite{AtomicDecompBP}, regularize it, and then rewrite the resulting problem as~\eqref{Eq:ColBPRegularized}.
    Consequently, we recover a solution of~\eqref{bp} if~\eqref{Eq:ColBPRegularized} is solved for a sufficiently small~$\delta$. The benefit of solving~\eqref{Eq:ColBPRegularized} is that it is immediate to recover the solution of~\eqref{Eq:ColBPRegularized} from its dual solution. We are unaware of any strategy for choosing~$\delta$ without first solving~\eqref{bp}. We will thus  adopt a trial-and-error strategy.

\mypar{Dual problem} We use duality because the dual problem
of~\eqref{Eq:ColBPRegularized} will have terms
involving~$A^\top$. Since~$A$ is partitioned by columns, $A^\top$ will
be partitioned by rows. Therefore, the algorithm for the row partition
will be applicable with some minor modifications.

    The dual problem of~\eqref{Eq:ColBPRegularized} is
    \begin{equation}\label{Eq:CPDualRegularized}
			\isdraft{
				\underset{y}{\text{maximize}} \,\,\, L(y)\,,
			}{
        \begin{array}{cl}
          \textrm{maximize} & L(y) \isdraft{\vspace{-0.25cm}}{}\\
          y &
        \end{array}\,,
			}
    \end{equation}
    where the dual function is
    $
        L(y) = -b^\top y + \inf_x (\|x\|_1 + (A^\top y)^\top x + \frac{\delta}{2}\|x\|^2)
    $, and~$y \in \mathbb{R}^m$ the dual variable.
    To keep the notation consistent with the previous section, we recast~\eqref{Eq:CPDualRegularized} as a minimization problem:
    \begin{equation}\label{Eq:CPDualRegularized2}
			\isdraft{
				\underset{y}{\text{minimize}}\,\,\, b^\top y + \Psi(y)
			}{
        \begin{array}{cl}
          \textrm{minimize} & b^\top y + \Psi(y) \isdraft{\vspace{-0.25cm}}{} \\
          y &
        \end{array}\,,
			}
    \end{equation}
    where
    \begin{equation}\label{Eq:CPPsi}
        \Psi(y) = -\inf_x \left(\|x\|_1 + (A^\top y)^\top x + (\delta/2)\|x\|^2\right).
    \end{equation}
    The objective of the inner optimization problem of~\eqref{Eq:CPPsi} has a unique minimizer for each~$y$, since it is strictly convex. Let~$x(y)$ denote the solution of this problem, for a fixed~$y$. Strong duality holds for~\eqref{Eq:ColBPRegularized} because its objective is convex and its constraints linear~\cite[\S5.2.3]{Boyd:ConvexOpti}, \cite[prop.5.2.1]{Bertsekas:Nonlinear}. Therefore, after we find a solution~$y^\star$ to the dual problem~\eqref{Eq:CPDualRegularized2}, a (primal) solution of~\eqref{Eq:ColBPRegularized} is available as~$x(y^\star)$. This follows directly from the KKT conditions~\cite[\S5.5]{Boyd:ConvexOpti}, \cite[prop.5.1.5]{Bertsekas:Nonlinear}, and we express it in the following theorem.

    \begin{Theorem}\label{Thm:RecoverPrimalFromDual}
        Let~$y^\star$ solve~\eqref{Eq:CPDualRegularized}. Then, $x(y^\star)$ solves~\eqref{Eq:ColBPRegularized}.
    \end{Theorem}

\mypar{Adapting the algorithm}
Now we focus on solving~\eqref{Eq:CPDualRegularized2}.
Let~$x$ be partitioned analogous to~$A$, i.e., $x = (x_1,\ldots,x_P)$, where~$x_p \in \mathbb{R}^{n_p}$. Note that~$\Psi(y)$ can be decomposed as the sum of~$P$ functions:
    $
        \Psi(y) = \Psi_1(y) + \cdots + \Psi_P(y)\,,
    $
    where
    \begin{equation}\label{Eq:CPPsiP}
        \Psi_p(y) = -\inf_{x_p} \|x_p\|_1 + (A_p^\top y)^\top x_p + \frac{\delta}{2}\|x_p\|^2
    \end{equation}
    can only be computed at node~$p$ because~$A_p$ is only known there. We can then rewrite~\eqref{Eq:CPDualRegularized2} as
    $$
    \isdraft{
			\underset{y}{\text{minimize}} \,\,\, \sum_{p=1}^P \bigl(\frac{1}{P}b^\top y + \Psi_p(y)\bigr)
		}{
        \begin{array}{cl}
          \textrm{minimize} & \sum_{p=1}^P \bigl(\frac{1}{P}b^\top y + \Psi_p(y)\bigr) \\
            y &
        \end{array}\,.
		}
    $$
    Notice that~$\Psi_p(y)$ can be easily computed at node~$p$, since the optimization problem defining it has a closed form solution. We now apply the same procedure as in section~\ref{Sec:RowPartition}: we clone the variable~$y$ into several~$y_p$'s, and constrain the problem with~$y_i = y_j$, for all~$\{i,j\} \in \mathcal{E}$. This yields
    \begin{equation}\label{Eq:CPClonedProb}
        \begin{array}{ll}
          \textrm{minimize} & \sum_{p=1}^P \bigl(\frac{1}{P}b^\top y_p + \Psi_p(y_p) \bigr) \vspace{0.5ex} \\
          \textrm{subject to} & (B^\top \otimes I_n)\bar{y} = 0\,,
        \end{array}
    \end{equation}
where the variable is~$\bar{y} = (y_1,\ldots,y_P) \in
(\mathbb{R}^n)^P$. Note the similarity between~\eqref{Eq:CPClonedProb}
and~\eqref{Eq:ReformEdgesB}. Having a proper coloring of the graph, the generalized ADMM is applicable:
    \begin{algorithm}[H]
    \caption{D-ADMM for general graphs (column partition)}
    \algrenewcommand\algorithmicrequire{\textbf{Initialization:}}
    \label{Alg:ADMMGeneralCP}
    \begin{algorithmic}[1]
    \small
    \Require for all~$p \in \mathcal{V}$, set $\gamma_{p}^{(1)} = x_p^{(1)} = 0$ and $k=1$
    \Repeat
    \For{$c =1,\ldots,C$}
        \ForAll{$p \in \mathcal{C}_c$ [in parallel]}
            $$
                v_p^{(k)} = \gamma_p^{(k)}-
                \rho \sum_{\begin{subarray}{c}
                             j \in \mathcal{N}_p \\
                             j < p
                           \end{subarray}
                }x_j^{(k+1)} - \rho \sum_{\begin{subarray}{c}
                             j \in \mathcal{N}_p \\
                             j > p
                           \end{subarray}
                }x_j^{(k)}
            $$
        \State and find
            $$
                y_p^{(k+1)} = \arg\min_{y_p} \Psi_p(y_p) + (v_p^{(k)}+\frac{1}{P}b)^\top y_p + \frac{D_p \rho}{2}\|y_p\|^2
            $$
           \label{SubAlg:ADMMGenProbCP}
        \State Send~$y_p^{(k+1)}$ to $\mathcal{N}_p$
        \label{SubAlg:ADMMGenCommCP}
        \EndFor
    \EndFor

    \ForAll{$p =1,\ldots,P$ [in parallel]}\vspace{0.15cm}
    \hfill
    
        $
            \gamma_p^{(k+1)} = \gamma_p^{(k)} + \rho \sum_{j \in \mathcal{N}_p}(y_p^{(k+1)} -  y_j^{(k+1)})
        $\vspace{0.15cm}
        \label{SubAlg:ADMMGenDualVarUpCP}
    \EndFor
    \State $k \gets k+1$
    \Until{some stopping criterion is met}
    \end{algorithmic}
    \end{algorithm}

    Algorithm~\ref{Alg:ADMMGeneralCP} is similar to Algorithm~\ref{Alg:ADMMGeneral} except for some minor modifications: the size of the variable to be transmitted is smaller (instead of transmitting~$x_p \in \mathbb{R}^n$, now the nodes transmit~$y_p \in \mathbb{R}^m$), and the optimization problem to be solved at each node (see step~\ref{SubAlg:ADMMGenProbCP}) is slightly different. Since that problem is unconstrained and its objective is differentiable, we can solve it directly with the Barzilai-Borwein algorithm~\cite{Barzilai-Borwein} (see appendix~\ref{App:OptimForEachNode} for more details).

Another difference to Algorithm~\ref{Alg:ADMMGeneral} is that after
the algorithm finished (finding an optimal vector~$y^\star$), node~$p$
will not know the entire solution~$x(y^\star)$
to~\eqref{Eq:ColBPRegularized}, but only a portion of it,
$x_p(y^\star)$, as
the solution to the optimization problem defining~$\Psi_p$
in~\eqref{Eq:CPPsiP}. In case we want the entire solution~$x(y^\star)$ to be available in all nodes, just a few additional communications are required because~$x(y^\star)$ is expected to be sparse; furthermore, a spanning tree can be used to spread the~$x_p$'s over the network.

We remark that if the graph is bipartite, then
Algorithm~\ref{Alg:ADMMGeneralCP} is proven to converge to an optimal
solution of~\eqref{Eq:ColBPRegularized} and, if~$\delta$ is small
enough, to a solution of~\eqref{bp}. An important issue is the possible ill-conditioning provoked by a small value of~$\delta$. In fact, a very small value for~$\delta$ may lead to difficulties in finding~$y^{(k+1)}_p$ in step~\ref{SubAlg:ADMMGenProbCP}. Note that this is the only step where~$\delta$ appears. In our simulations, explained in section~\ref{Sec:ExperimentalResults}, we used~$\delta = 10^{-3}$ and this value allowed us to compute solutions to BP with a very large precision, without incurring into numerical problems.

\section{Other Algorithms}
\label{Sec:OtherAlgs}

In this section we overview other methods that solve BP in a truly
distributed way. We only cover the row partition case because
corresponding algorithms for the column partition can always be derived
as shown in the previous section.

We divide the algorithms into two categories according to the number of (nested) loops they have: single-looped and
double-looped. D-ADMM is single-looped and, in each
iteration, every node transmits a vector of size~$n$ to its neighbors.

\mypar{Performance measure: communication steps} We say that a communication step has occurred after all the nodes finish communicating their current estimates to their neighbors. All single-looped algorithms have one communication step per iteration. The double-looped algorithms
have one communication step per iteration
of the inner loop. In all algorithms, the size of the
transmitted vector is~$n$. Another feature common to all algorithms is
that in every iteration (or in every inner iteration, for the
double-looped algorithms) each node has to solve the optimization
problem in step~\ref{SubAlg:ADMMGenProb} of
Algorithm~\ref{Alg:ADMMGeneral} (or Algorithm~\ref{Alg:ADMMGeneralCP},
for the column partition). This means that the algorithms have a
common ground for comparison: if each iteration (or inner iteration,
for the double-looped algorithms) involves one communication step and
all the nodes have to solve a similar optimization problem (same
format, same dimensions, but possibly different parameters), then the
number of iterations (or the sum of inner iterations) becomes a
natural metric to compare the algorithms. We will then compare the algorithms by their number of communication steps, which is equal to the number of iterations in the single-looped algorithms and to the sum of inner iterations in the double-looped algorithms. Note that less communication steps 
can be expected to produce significant energy savings
in scenarios such as sensor networks~\cite{Akyildiz02-WirelessNetworksASurvey}.

Although data is transmitted in every communication step, the quantity of the transmitted data might actually decrease with the iterations. The reason is because the solution to BP is sparse and, at some point, the nodes' estimates start being sparse, allowing a possible compression of the transmitted data (e.g., just transmit the nonzero entries).

We start with describing the single-looped algorithms.

\mypar{Subgradient} Nedi$\acute{\textrm{c}}$ and Ozdaglar were the
first to propose a subgradient-based algorithm to solve general
convex optimization problems in a completely distributed
way~\cite{Nedic09DistributedSubgradientMultiAgent}. However, they only
addressed unconstrained optimization problems, which is not our
case. Instead, we will use the method proposed in~\cite{MultiAgent},
which generalizes~\cite{Nedic09DistributedSubgradientMultiAgent} to
problems with private constraints in each node. That is,
\cite{MultiAgent} solves
    $$
        \begin{array}{ll}
          \textrm{minimize} & \sum_{p=1}^P f_p(x) \\
          \textrm{subject to} & x \in \cap_{p=1}^P X_p\,,
        \end{array}
    $$
    where each~$f_p$ is convex and each~$X_p$ is a closed convex set. This method combines consensus algorithms~\cite{DeGroot74ReachingConsensus} with subgradient algorithms~\cite[Ch.6]{Bertsekas:Nonlinear}, and for each node~$p$, it takes the form
    \begin{equation}\label{Eq:OASubgradAlg}
        x_p^{(k+1)} = \Bigl[c_{pp}^{(k)} x_p^{(k)} + \sum_{j \in \mathcal{N}_p} c_{pj}^{(k)} x_j^{(k)} - \alpha^{(k)} g_p^{(k)}\Bigr]_{X_p}^+\,,
    \end{equation}
    where~$c_{ij}$ are positive weights such that $\sum_{i} c_{ij}^{(k)} = \sum_{j} c_{ij}^{(k)} = 1$, the sequence~$\{\alpha^{(k)} > 0\,:\, k=1,2,\ldots\}$ is square summable but not summable, and~$[p]_{X}^+$ is the projection of the point~$p$ onto the set~$X$:
    $
      [p]_{X}^+ = \arg\min_x \, \{ \frac{1}{2} \|x - p\|^2 \,:\, x \in X\}\,.
    $
    The vector~$g_p^{(k)}$ is a subgradient of~$f_p$ at the point $c_{pp}^{(k)} x_p^{(k)} + \sum_{j \in \mathcal{N}_p} c_{pj}^{(k)} x_j^{(k)}$.

    We apply~\eqref{Eq:OASubgradAlg} directly to problem~\eqref{Eq:RPReformulation1}, where we see~$\|x\|_1$ as~$\|x\|_1 = \frac{1}{P}\|x\|_1 + \cdots + \frac{1}{P}\|x\|_1$; in other words, we set $f_p(x) = \frac{1}{P}\|x\|_1$. We choose~$\alpha^{(k)} = 1/(k+1)$ for the step-size sequence. In our case, since the network is static (Assumption~\ref{Ass:Network}), the weights~$c_{ij}$ are constant: for every~$p$, $c_{pi} = 1/(D_p+1)$ for~$i \in \mathcal{N}_p \cup \{p\}$, and~$0$ otherwise. The implementation of~\eqref{Eq:OASubgradAlg} in a network is now straightforward: first, node~$p$ transmits~$x_p^{(k)}$ to its neighbors and receives~$x_j^{(k)}$ from them; then, it updates its variable with~\eqref{Eq:OASubgradAlg}. These two steps are repeated until convergence.

    While~\eqref{Eq:OASubgradAlg} is proven to be robust to link failures, its convergence speed is too slow in practice.

\mypar{D-Lasso}
    As mentioned in section~\ref{Sec:intro}, Bazerque and Giannakis~\cite{Giannakis} proposed a distributed algorithm that solves a problem similar to ours. Here, we adapt it to solve BP. The starting point is problem~\eqref{Eq:ReformEdges}, which by introducing a new variable~$z_{ij}$ for each edge~$\{i,j\} \in \mathcal{E}$, is reformulated as
    \begin{equation}\label{Eq:OADLasso1}
        \begin{array}{ll}
          \textrm{minimize} & \frac{1}{P} \sum_{p=1}^P \|x_p\|_1 \\
          \textrm{subject to} & A_p x_p = b_p\,,\quad p=1,\ldots,P \\
          \isdraft{&x_i = z_{ij}, \,\,\, x_j = z_{ij}\,,\quad  \{i,j\} \in \mathcal{E}.}{
          & x_i = z_{ij}\,,\quad \{i,j\} \in \mathcal{E}, \\
          & x_j = z_{ij}\,,\quad \{i,j\} \in \mathcal{E}.
					}
        \end{array}
    \end{equation}
    This problem is solved with ADMM by dualizing its last two constraints. We consider the problem partitioned in terms of the variable~$\bar{z} = (\ldots,z_{ij},\ldots)$ and~$\bar{x}=(\ldots,x_p,\ldots)$. In short, ADMM minimizes the augmented Lagrangian of~\eqref{Eq:OADLasso1} w.r.t. $\bar{z}$ and then minimizes it w.r.t. $\bar{x}$, using the new value of~$\bar{z}$. The minimization w.r.t. $\bar{z}$ has a closed form solution. After some manipulations, the algorithm for an arbitrary node~$p$ is:
    \begin{algorithm}[H]
    \caption{D-Lasso (node~$p$)}
    \algrenewcommand\algorithmicrequire{\textbf{Initialization:}}
    \label{Alg:DLasso}
    \begin{algorithmic}[1]
    \small
    \Require for all~$p \in \mathcal{V}$, set $\gamma_{p}^{(1)} = x_p^{(1)} = 0$ and $k=1$
    \Repeat
    \ForAll{$p=1,\ldots,P$ [in parallel]}
    \hfill
    
		\State set
            $
                v_p^{(k)} = \gamma_p^{(k)}-
                \rho \sum_{j \in \mathcal{N}_p \cup \{p\}} x_j^{(k)}
            $
        and find
            $$
            x_p^{(k+1)} = \begin{array}[t]{cl}
                            \underset{x_p}{\textrm{argmin}} & \frac{1}{P}\|x_p\|_1 + {v_p^{(k)}}^\top x_p + \rho D_p\|x_p\|^2\\
                            \textrm{s.t.} & A_p x_p = b_p
                          \end{array}
           $$
    \label{Step:DLasso-Optm}
    
    \State Send~$x_p^{(k+1)}$ to~$\mathcal{N}_p$, and receive~$x_j^{(k+1)}$, $j \in \mathcal{N}_p$
    \EndFor
    \ForAll{$p =1,\ldots,P$ [in parallel]} \vspace{0.15cm}
    \hfill
    
        $
            \gamma_p^{(k+1)} = \gamma_p^{(k)} + \rho \sum_{j \in \mathcal{N}_p} (x_p^{(k+1)} -  x_j^{(k+1)})
        $\vspace{0.15cm}
        \label{SubAlg:ADMMGenDualVarUp}
    \EndFor
    \State $k \gets k+1$
    \Until{some stopping criterion is met}
    \end{algorithmic}
    \end{algorithm}
    Although D-Lasso and D-ADMM (Algorithm~\ref{Alg:ADMMGeneral}) have a similar format, they are different. For example, D-Lasso is synchronous and D-ADMM asynchronous, and the parameters of the optimization problem each node solves are different in both algorithms. Also, D-ADMM is proven to converge for bipartite graphs only, while D-Lasso is proven to converge for any connected graph. In the next section, we will see that, in practice, D-ADMM converges in less iterations than D-Lasso, despite their common underlying algorithm.

    We now move to the double-looped algorithms.

\mypar{Double-looped algorithms}
    All double-looped algorithms we will see have the same theoretical foundation, but use different subalgorithms. Namely, all solve the following dual problem of~\eqref{Eq:ReformEdgesB}:
    \begin{equation}\label{Eq:OADualProb}
        \begin{array}{cl}
          \textrm{maximize} & L(\lambda) \isdraft{\vspace{-0.25cm}}{} \\
          \lambda &
        \end{array}\,,
    \end{equation}
    where~$L(\lambda)$ is the augmented dual function
    \begin{equation}\label{Eq:OADualFunct}
        L(\lambda) =
        \begin{array}[t]{cl}
          \inf & \sum_{p=1}^P \frac{1}{P}\|x_p\|_1 + \sum_{\{i,j\} \in \mathcal{E}} \phi_{\lambda_{\{i,j\}}} (x_i - x_j) \\
          \textrm{s.t.} & A_p x_p = b_p\,,\quad p = 1,\ldots,P\,,
        \end{array}
    \end{equation}
    where $\phi_{\lambda}(z) = \lambda^\top z + \frac{\rho}{2}\|z\|^2$, and~$\rho$ is a positive parameter.
    The algorithms have an outer loop that solves~\eqref{Eq:OADualProb}, and an inner loop that solves the optimization problem in~\eqref{Eq:OADualFunct}.

    We consider three distributed, double-looped algorithms~\cite{Spars09,Ruszczynski,Mota11ICASSP} to solve~\eqref{Eq:OADualProb}, and thus~\eqref{Eq:ReformEdgesB} because strong duality holds. While~\cite{Spars09,Mota11ICASSP} were designed to solve BP, \cite{Ruszczynski} was designed to solve more general problems. We thus have to adapt the latter to our problem. The algorithms described in~\cite{Spars09,Ruszczynski,Mota11ICASSP} will be denoted respectively by MM/NGS (method of multipliers and nonlinear Gauss-Seidel), MM/DQA (method of multipliers and diagonal quadratic approximation), and DN (double Nesterov).

All algorithms solve~\eqref{Eq:OADualProb} with an iterative scheme
in the outer loop. As in D-ADMM, the dual
variable~$\lambda$ consists of several variables~$\lambda_{\{i,j\}}$
associated with the edges~$\{i,j\} \in \mathcal{E}$. It can be
shown that the dual function~$L(\lambda)$ in~\eqref{Eq:OADualFunct}
is differentiable and that its gradient~$\nabla
L(\lambda) = (\ldots,x_i(\lambda) - x_j(\lambda),\ldots)$ is Lipschitz
continuous with constant~$1/\rho$~\cite{Konnov07}. The
vector~$\bar{x}(\lambda) :=
(x_1(\lambda),x_2(\lambda),\ldots,x_P(\lambda))$ solves the
optimization problem in~\eqref{Eq:OADualFunct} for a
fixed~$\lambda$. The algorithm for solving this inner problem will
be the inner loop and is considered later. These nice
properties of~$L(\lambda)$ enable the edge-wise application of the gradient
method~\cite[\S1.2]{Bertsekas:Nonlinear}
    \begin{equation}\label{Eq:OADLGradMethod}
        \lambda_{\{i,j\}}^{(k+1)} = \lambda_{\{i,j\}}^{(k)} + \rho \nabla_{\lambda_{\{i,j\}}} L(\lambda^{(k)})\,,
    \end{equation}
    or the edge-wise application of Nesterov's method~\cite{Nesterov}
    \begin{equation}\label{Eq:OADLNesterov}
        \begin{array}{ll}
          \lambda_{\{i,j\}}^{(k+1)} &= \eta_{\{i,j\}}^{(k)} + \rho \nabla_{\eta_{\{i,j\}}} L(\eta^{(k)}) \vspace{1ex}\\
          \eta_{\{i,j\}}^{(k+1)} &= \lambda_{\{i,j\}}^{(k+1)} + \frac{k-1}{k+2}(\lambda_{\{i,j\}}^{(k+1)} - \lambda_{\{i,j\}}^{(k)})\,,
        \end{array}
    \end{equation}
    to solve~\eqref{Eq:OADualProb}. Nesterov's method is proven to be faster than the gradient method. When we use the gradient method~\eqref{Eq:OADLGradMethod} to solve a dual problem, where duality here is seen in the augmented Lagrangian sense, the resulting algorithm is called method of multipliers (MM)~\cite[p.408]{Bertsekas:Nonlinear}. While MM/NGS and MM/DQA use MM for their outer loop, DN uses~\eqref{Eq:OADLNesterov}.

    So far, we assumed that a solution of the optimization problem in~\eqref{Eq:OADualFunct}, for a given~$\lambda$, was available. Nevertheless, solving this problem in a distributed way is more challenging than solving~\eqref{Eq:OADualProb} (when~$\nabla L(\lambda)$ is readily available). The reason is that we cannot decouple the term~$\sum_{\{i,j\} \in \mathcal{E}} \phi_{\lambda_{\{i,j\}}} (x_i - x_j)$ into a sum of~$P$ functions, each one depending only on~$x_p$. Both MM/NGS and MM/DQA use an iterative method that optimizes the objective of~\eqref{Eq:OADualFunct} w.r.t. one block variable~$x_p$, while keeping the other blocks fixed. More concretely, let~$g_{\lambda}(x_1,\ldots,x_P)$ denote the objective of~\eqref{Eq:OADualFunct} when~$\lambda$ is fixed. MM/NGS uses the nonlinear Gauss-Seidel (NGS) method~\cite[\S3.3.5]{Bertsekas:Parallel}\cite{Tseng}:
    \begin{align}
          x_1^{(t+1)} &= \arg\min_{x_1 \in X_1} g_{\lambda}(x_1,x_2^{(t)},x_3^{(t)},\ldots,x_P^{(t)})
          \notag\\
          x_2^{(t+1)} &= \arg\min_{x_2 \in X_2} g_{\lambda}(x_1^{(t+1)},x_2,x_3^{(t)},\ldots,x_P^{(t)})
          \notag\\
          \vspace{-0.2cm}
          &\phantom{=} \vdots
          \label{Eq:OADLNGS}
          \vspace{-0.2cm}
          \\
          x_P^{(t+1)} &= \arg\min_{x_P \in X_P} g_{\lambda}(x_1^{(t+1)},x_2^{(t+1)},x_3^{(t+1)},\ldots,x_P)\,,
          \notag
    \end{align}
    where~$X_p := \{x_p\,:\, A_p x_p = b_p\}$, $p=1,\ldots,P$. It can be proven that any limit point of the sequence generated by~\eqref{Eq:OADLNGS} solves problem~\eqref{Eq:OADualFunct}; see~\cite{Tseng,JoaoMotaThesis}. Each optimization problem in~\eqref{Eq:OADLNGS} is solved at one node. It turns out that these are equivalent to the problem in step~\ref{SubAlg:ADMMGenProb} of Algorithm~\ref{Alg:ADMMGeneral}. Note that the nodes in~\eqref{Eq:OADLNGS} cannot operate in parallel, akin to the algorithm we propose here. MM/DQA, on the other hand, solves the problem in~\eqref{Eq:OADualFunct} with a parallel scheme called diagonal quadratic approximation (DQA):
    \begin{align}
          &u_1 = \arg\min_{x_1 \in X_1} g_{\lambda}(x_1,x_2^{(t)},x_3^{(t)},\ldots,x_P^{(t)})
          \notag\\
          &u_2 = \arg\min_{x_2 \in X_2} g_{\lambda}(x_1^{(t)},x_2,x_3^{(t)},\ldots,x_P^{(t)})
          \notag\\
          \vspace{-0.2cm}
          &\phantom{u_2 =} \vdots
          \label{Eq:OADLDQA}
          \\
          &u_P = \arg\min_{x_P \in X_P} g_{\lambda}(x_1^{(t)},x_2^{(t)},x_3^{(t)},\ldots,x_P)
          \notag
          \vspace{-0.2cm}
          \\
          &x_p^{(t+1)} = \tau u_p + (1-\tau) x_p^{(t)}\,,\quad p=1,\ldots,P\,,
          \notag
    \end{align}
    where~$\tau = 1/P$. For a proof that~\eqref{Eq:OADLDQA} solves~\eqref{Eq:OADualFunct} see~\cite{Ruszczynski,JoaoMotaThesis}. The difference between~\eqref{Eq:OADLNGS} and~\eqref{Eq:OADLDQA} is that the latter allows all the nodes to operate in parallel, and after the minimization step, each node combines the solution of the optimization problem it has just solved with the previous estimate of the solution: $x_p^{(t)}$. Note that a communication step has to occur after each iteration.

    Regarding DN, we made some modifications to the inner loop of the method proposed in~\cite{Mota11ICASSP}, so that we could get an algorithm comparable with what we propose here.

\mypar{Double Nesterov (DN)}
    In~\cite{Mota11ICASSP}, BP is recast as a linear program by increasing the size of the variable to~$2n$. The result is that the problem defining the dual function has a differentiable objective with a Lipschitz continuous gradient, and thus Nesterov's method is directly applicable. However, the size of the variable transmitted in each communication step is~$2n$. Here, we do not recast BP as an LP. As seen before, the dual problem~\eqref{Eq:OADualProb} is solved with Nesterov's method~\eqref{Eq:OADLNesterov} in the outer loop. Now, to solve the optimization problem in~\eqref{Eq:OADualFunct}, Nesterov's method is not applicable because the objective is not differentiable. However, that objective can be written as the sum of a nondifferentiable function $h(\bar{x}) = \sum_{p=1}^P \frac{1}{P}\|x_p\|_1$ with a differentiable one $g(\bar{x}) = \sum_{\{i,j\} \in \mathcal{E}} \phi_{\lambda_{\{i,j\}}} (x_i - x_j)$. The gradient of~$g(\bar{x})$ w.r.t. $x_p$ is $\nabla_{x_p} g(\bar{x}) = \gamma_p + \rho D_p x_p - \rho \sum_{j \in \mathcal{N}_p}x_j$. Therefore, to compute~$\nabla_{x_p} g(\bar{x})$, each node needs only to communicate with its neighbors. The gradient~$\nabla g(\bar{x})$ is Lipschitz continuous with constant~$\rho \lambda_{\max}(\mathcal{L})$, where~$\lambda_{\max}(\mathcal{L})$ denotes the maximum eigenvalue of the graph Laplacian. FISTA~\cite{BeckTeboulleFISTA} is an algorithm that adapts Nesterov's method to this scenario. It operates the following way:
    \begin{algorithm}[H]
    \caption{FISTA (for node~$p$)}
    \algrenewcommand\algorithmicrequire{\textbf{Initialization:}}
    \label{Alg:FISTA}
    \begin{algorithmic}[1]
    \small
    \Require choose $\alpha = 1/(\rho \lambda_{\max}(\mathcal{L}))$, $x_p^{(0)} = y_p^{(0)} = 0$, $t = 0$
    \Repeat
    \State $u_p = y_p^{(t)} - \alpha\nabla g(y_p^{(t)})$
    \State
        $
            x_p^{(t+1)} = \arg\min_{x_p} \frac{1}{2\alpha}\|x_p-u_p\|^2 + h(x_p)
        $
        \label{SubAlg:FISTAOptim}
    \State $y_p^{(k+1)} = x_p^{(k+1)} + \frac{k-1}{k+2}\left(x_p^{(k+1)} - x_p^{(k)}\right)$
    \State $k \gets k+1$
    \Until{some stopping criterion is met}
    \end{algorithmic}
    \end{algorithm}
This modification to~\cite{Mota11ICASSP} allows us to compare the
resulting algorithm with ours, because the size of the
variable is now~$n$. Furthermore, the problem
in step~\ref{SubAlg:FISTAOptim} is equivalent to the one in
step~\ref{SubAlg:ADMMGenProb} of Algorithm~\ref{Alg:ADMMGeneral}.

    \mypar{Tuning parameter \boldmath{$\rho$}} Note that all algorithms (except the subgradient) share the same tuning parameter~$\rho$, because all are based on an augmented Lagrangian reformulation. It is known that~$\rho$ influences the convergence rate of augmented Lagrangian methods. Albeit there are self-adaptive schemes to update~$\rho$ during the algorithm~\cite{BoydADMM,He00ADMWithSlefAdaptativePenaltyForVIs,Liao04DecompositionMethodVariableParameterVI}, making the algorithms less sensitive to~$\rho$, we were not able to implement these schemes in a distributed scenario. We will hence assume~$\rho$ is constant during the execution of the algorithms.
    
    \mypar{Execution times in wireless networks} In contrast with all the algorithms described here (except MM/NGS), D-ADMM assumes a coloring scheme based on which the nodes operate asynchronously. Suppose all the algorithms are implemented on an ideal network, where packet collisions do not occur, i.e., two neigboring nodes can transmit messages at the same time without causing interference at the reception. If a communication step by D-ADMM takes~$T$ time units, then a communication step by the other algorithms takes~$T/C$ units, where~$C$ is the number of colors we used for the network (we are ignoring the optimizations that can be made from the procedure described in Figure~\ref{Fig:UndirectedAndDirected}). Therefore, although D-ADMM requires less communication steps, as shown next, it might actually take longer than competing algorithms. However, in a real wireless network, packet collisions occur and medium-access (MAC) protocols have to be implemented to avoid them. Hence, synchronous algorithms cannot operate synchronously in wireless networks. 
    The execution time of an algorithm, among other factors, is highly dependent on the MAC protocol. Comparing execution times is thus beyond the scope of this paper.


\section{Experimental Results}
\label{Sec:ExperimentalResults}

In this section we compare our algorithm against the prior work
discussed in the previous section and listed in~Table~\ref{Tab:Algorithms}.
We focus on the row-partitioned case since the algorithm for
the column partition is derived from it. We start describing how the data and the networks were generated, and how the experiments were carried out. In the first type of experiments we compare all the algorithms on moderate-sized networks (around~$50$ nodes) and conclude that D-ADMM and D-Lasso are the ``fastest'' algorithms. In the second type of experiments we compare only these two algorithms in a more thorough way for the same networks, and we also see how their performance varies as the network size increases (from~$2$ nodes to~$1024$ nodes). Finally, we address the column partition case.
    
    \begin{table}
    \centering
    \caption{
             Algorithms for comparison in the simulations.
            }
    \label{Tab:Algorithms}
    \isdraft{\vspace{-0.3cm}}{\smallskip}
        \renewcommand{\arraystretch}{1.3}
        \isdraft{\scriptsize}{}
        \begin{tabular}{@{}lrlcc@{}}
          \toprule[1pt]
          Acronym & & Algorithm(s) & & Source \\
          \midrule
          D-ADMM      & & Alternating direction MM               & & This paper \\
          Subgradient & & Subgradient method                     & & \cite{MultiAgent} \\
          D-Lasso     & & Alternating direction MM               & & \cite{Giannakis} \\
          MM/NGS      & & MM + nonlinear Gauss-Seidel            & & \cite{Spars09} \\
          MM/DQA      & & MM + diagonal quadratic approximation  & & \cite{Ruszczynski} \\
          DN          & & Nesterov + Nesterov                    & & \cite{Mota11ICASSP} \\
          \bottomrule[1pt]
        \end{tabular}
        \isdraft{\vspace{-0.6cm}}{}
    \end{table}
    
    \begin{table}
    \isdraft{\begin{minipage}[t]{0.46\linewidth}}
    \centering
    \caption{
             Scenarios for row partition experiments.
            }
    \label{Tab:Scenarios}
    \isdraft{\vspace{-0.3cm}}{\smallskip}
        \renewcommand{\arraystretch}{1.3}
        \isdraft{\scriptsize}{}
        \begin{tabular}{@{}crcrrrcrc@{}}
          \toprule[1pt]
          Scenario & & Sparco Id & & $m$\phantom{a} & & $n$ & & $P$\\
          \midrule
           1 & & -----    & & 500   & & 2000 & & 50 \\
           2 & & \phantom{11}7   & & 600   & & 2560 & & 50 \\
           3 & & \phantom{11}3   & & 1024  & & 2048 & & 64 \\
           4 & & 902 & & 200   & & 1000 & & 50 \\
           5 & & \phantom{1}11  & & 256   & & 1024 & & 64 \\
          \bottomrule[1pt]
        \end{tabular}
        \isdraft{\vspace{-0.6cm}}{}
    \isdraft{\end{minipage}}{}
    \isdraft{}{\end{table}
    \begin{table}
    }
    \isdraft{
    \hfill
    \begin{minipage}[t]{0.5\linewidth}}{}
    \centering
    \caption{
             Network models for the experiments.
            }
    \label{Tab:Networks}
    \isdraft{\vspace{-0.3cm}}{\smallskip}
        \renewcommand{\arraystretch}{1.3}
        \isdraft{\scriptsize}{}
        \begin{tabular}{@{}crlrl@{}}
          \toprule[1pt]
          Network number & & Model & & Parameters \\
          \midrule
           1 & & Erd\H os-R\'enyi  & & $p = 0.25$ \\
           2 & & Erd\H os-R\'enyi  & & $p = 0.75$ \\
           3 & & Watts-Strogatz  & & $(n,p) = (4,0.6)$ \\
           4 & & Watts-Strogatz  & & $(n,p) = (2,0.8)$ \\
           5 & & Barabasi-Albert & & -----------------\\
           6 & & Geometric       & & $d = 0.75$ \\
           7 & & Lattice         & & -----------------\\
          \bottomrule[1pt]
        \end{tabular}
        \isdraft{\vspace{-0.6cm}}{}
    \isdraft{\end{minipage}}{}
    \end{table}

\mypar{Experimental setup} We considered five distinct scenarios with
different dimensions and different types of data, shown in
Table~\ref{Tab:Scenarios}. The data (matrix~$A \in \mathbb{R}^{m
  \times n}$ and vector~$b \in \mathbb{R}^m$) was taken from the
Sparco toolbox~\cite{Sparco}, except in scenario~$1$, where we used
a $500 \times 2000$ matrix with i.i.d.\ Gaussian entries with zero
mean and variance~$1/\sqrt{m}$. In each scenario, each node
stores~$m_p = m/P$ rows of~$A$. We ensured that $m_p = m/P$ is an
integer by considering two values for~$P$: $50$ and~$64$, chosen
depending on the scenario.

In the following, $x^\star$ denotes the solution of BP obtained by the Sparco toolbox, or in scenario~$1$, the one obtained by CVX~\cite{cvx}, solving BP as a linear program. Note that due to the dimensions of the matrices and their randomness/structure, $x^\star$ is guaranteed to be unique with overwhelming probability.

For each scenario we ran all algorithms for the seven different networks
shown in Table~\ref{Tab:Networks}. For each
network in Table \ref{Tab:Networks} we generated two networks: one
with~$50$ nodes (used in scenarios with $P=50$), the other with~$64$ nodes
(used in scenarios with $P=64$).
The parameters of the networks were chosen so that the
generated network would be connected with high probability.
Only for network~$4$, $P = 50$ we did not get a connected network at
first, so we changed the parameters to~$(3,0.8)$.
If the generated network had self-connections or multiple
edges between the same pair of nodes, we would remove them. We also generated~$10$ networks with~$2^i$ nodes ($i=1,\ldots,10$), all following the model of network~$3$. These are used in the type~II experiments (explained below).

The Erd\H os-R\'enyi model~\cite{Erdos59RandomGraphs} has one
parameter~$p$, which specifies the probability of any two nodes in the
network being connected. The Watts-Strogatz
model~\cite{Watts98CollectiveDynamicsSmallWorld} has two parameters:
the number of neighbors~$n$ and the rewiring probability~$p$.  First
it creates a lattice where every node is connected with~$n$ other
nodes; then, every link is rewired, or not, with probability~$p$. If a
rewiring occurs in link~$\{i,j\}$, then we pick
node~$i$ or~$j$ (with equal probability) and connect it with other
node in the network, chosen uniformly. The Barabasi-Albert
model~\cite{Barabasi99EmergenceOfScaling} starts with one node; at
each step, one node is added to the network and is connected to one of
the nodes already in the network. However, the probability of the new
node ``choosing'' to connect to the other nodes is not uniform: it is
proportional to the nodes' degrees such that the new node has a
greater probability of connecting to the nodes with larger degrees. The
geometric model~\cite{Penrose} deploys~$P$ nodes randomly (uniformly)
in the unit square; then, two nodes are connected if their distance is
less than~$d$. Finally, the Lattice model has no randomness. For~$P$ nodes, 
it generates a rectangular grid graph in the plane such that the shape is as square as possible
($5\times 10$ for $P=50$ and $8\times 8$ for $P=64$).  Each node has
four neighbors except for the borders.  This lattice network is the
only one guaranteed to be bipartite, and thus
Algorithm~\ref{Alg:ADMMGeneral} is only guaranteed to converge for
this network.

We used an heuristic from the Matgraph toolbox~\cite{Matgraph} to find a coloring for these networks. It is then possible that the number of colors is larger
than~$\chi(\mathcal{G})$ . We checked that the optimal solution of two colors was found for the Lattice model.

	 \begin{table}
	 	\centering
	 	\caption{
	 		Types of experiments.
	 	}
	 	\label{Tab:TypeExperiments}
	 	\isdraft{\vspace{-0.3cm}}{\smallskip}
	 	\renewcommand{\arraystretch}{1.3}
	 	\isdraft{\scriptsize}{}
	 	\begin{tabular}{@{}crl@{}}
	 	\toprule[1pt]
	 	Type of experiment & & Value of $\rho$ \\
	 	\midrule
	 	I  & & $\rho = 1$ for D-ADMM and D-Lasso \\
	 	   & & $\rho = 10$ for MM/NGS, MM/DQA, and DN \\[2mm]
	 	II & & $\rho \in \{10^{-3}, 10^{-2}, 10^{-1}, 10^{0}, 10^{1}\}$ \\
	 	   & & the value that leads to the best results is picked \\
	 	\bottomrule[1pt]
	 	\end{tabular}
	 	\isdraft{\vspace{-0.6cm}}{}
	 \end{table}

\mypar{Results} As mentioned before, we keep the parameter~$\rho$
fixed during the execution of the algorithms (except for the
subgradient method, which has no~$\rho$). We picked~$\rho$ in two
different ways, yielding two types of experiments, shown in
Table~\ref{Tab:TypeExperiments}. In type~I, $\rho$ was always the same
for all scenarios and all networks: $\rho = 1$ for D-ADMM
(Algorithm~\ref{Alg:ADMMGeneral}) and for D-Lasso
(Algorithm~\ref{Alg:DLasso}), and $\rho = 10$ or the double-looped
algorithms MM/NGS, MM/DQA, and DN. These values were
chosen based on some pre-testing. In the type II experiments, given a
fixed scenario and network, we execute each algorithm for
several $\rho$'s and pick the one that yields the best result.  For the
type II experiments, we only considered the best two algorithms: D-ADMM
and D-Lasso.

The two types of experiments reflect two different philosophies in the
assessment of algorithms that depend on parameters: type I represents
real-world applications (the parameters are tuned for known data and
are then used unchanged); type II is more suited to assess the true
capabilities of the algorithm.

    \begin{figure*}
      
        \centering
        \isdraft{
        \begin{pspicture}(\linewidth,0.4cm)	    
            \rput(0.22\linewidth,0.9){\large \textbf{ \textsf{Accuracy:} \boldmath{$\mathsf{1\,\%}$}}}
            \rput(0.74\linewidth,0.9){\large \textbf{ \textsf{Accuracy:} \boldmath{$\mathsf{10^{-3}\,\%}$}}}
        \end{pspicture}
        }{
        \begin{pspicture}(\linewidth,0.7cm)	    
            \rput(0.22\linewidth,0.7){\large \textbf{ \textsf{Accuracy:} \boldmath{$\mathsf{1\,\%}$}}}
            \rput(0.74\linewidth,0.7){\large \textbf{ \textsf{Accuracy:} \boldmath{$\mathsf{10^{-3}\,\%}$}}}
        \end{pspicture}
        }
        
        \vspace{-1cm}
        
        \subfigure[Scenario 1]{\label{SubFig:RPExperimentsScen1}
            
            \isdraft{
	      \begin{pspicture}(0.48\linewidth,4.6cm)
                \rput(0.225\linewidth,2.55){\includegraphics[scale=0.37]{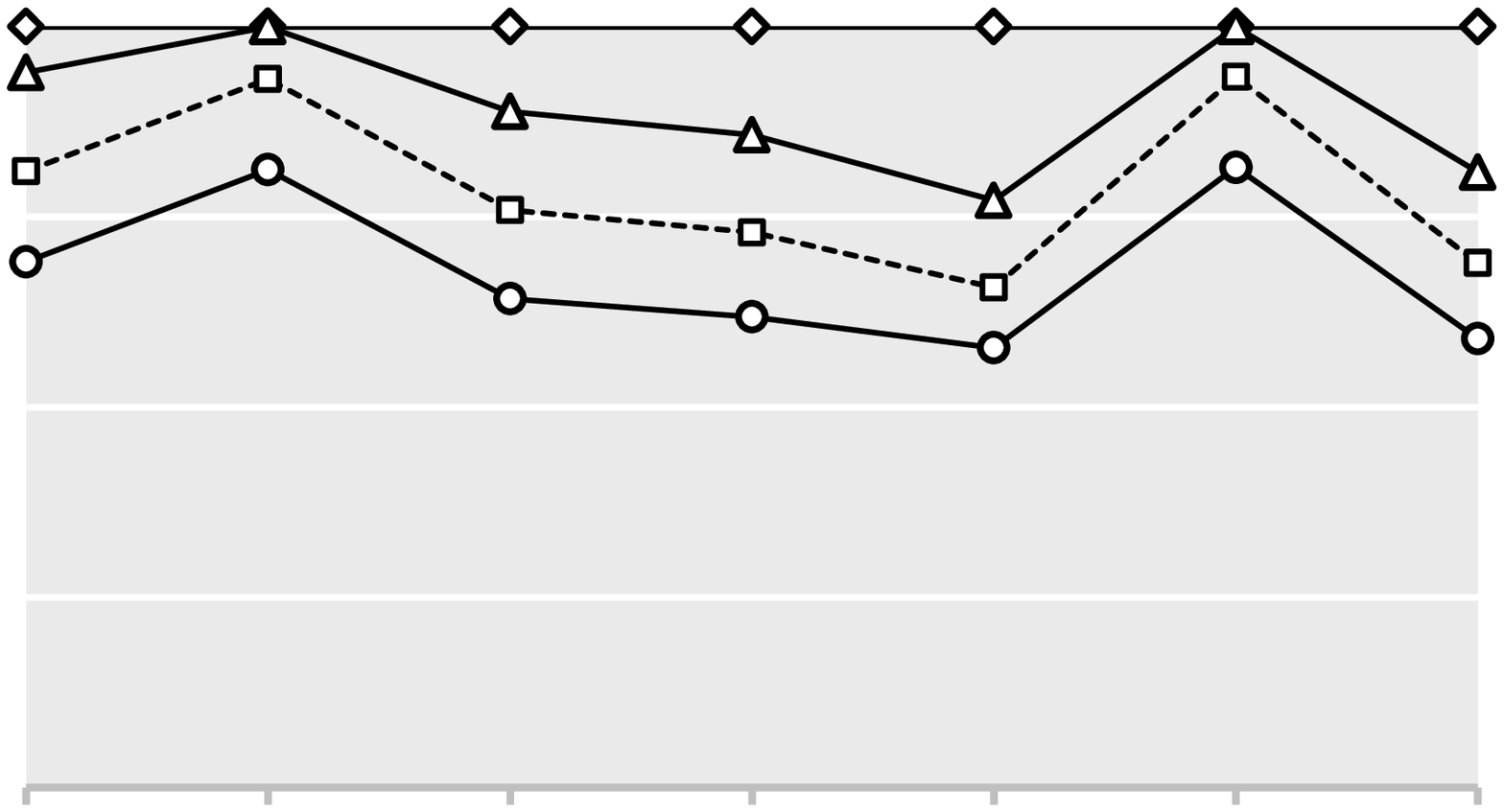}}
                \rput(0.225\linewidth,0.11){\footnotesize \textbf{\sf Network number}}
                \rput[bl](0.07,4.56){\mbox{\footnotesize \textbf{{\sf Communication steps}}}}
                \rput[l](0.05,0.91){\scriptsize $\mathsf{10^{0}}$}
                \rput[l](0.05,1.75){\scriptsize $\mathsf{10^{1}}$}
                \rput[l](0.05,2.58){\scriptsize $\mathsf{10^{2}}$}
                \rput[l](0.05,3.41){\scriptsize $\mathsf{10^{3}}$}
                \rput[l](0.05,4.25){\scriptsize $\mathsf{10^{4}}$}

                \rput[t](0.604,0.65){\scriptsize $\mathsf{1}$}
                \rput[t](1.630,0.65){\scriptsize $\mathsf{2}$}
                \rput[t](2.673,0.65){\scriptsize $\mathsf{3}$}
                \rput[t](3.721,0.65){\scriptsize $\mathsf{4}$}
                \rput[t](4.764,0.65){\scriptsize $\mathsf{5}$}
                \rput[t](5.800,0.65){\scriptsize $\mathsf{6}$}
                \rput[t](6.84,0.65){\scriptsize $\mathsf{7}$}

                \rput[l](7,2.81){\tiny \textbf{\sf D-ADMM}}
                \rput[l](7,3.19){\tiny \textbf{\sf D-Lasso}}
                \rput[l](7,3.59){\tiny \textbf{\sf MM/NGS}}
                \rput[l](7,4.32){\tiny \textbf{\sf MM/DQA, DN,}}
                \rput[l](7,4.06){\tiny \textbf{\sf Subgradient}}
              \end{pspicture}  
            }{
	      \begin{pspicture}(0.48\linewidth,4.9cm)
	        \rput(0.21\linewidth,2.55){\includegraphics[scale=0.38]{figures/RPCen1Acc2.eps}}
		\rput(0.21\linewidth,0.11){\footnotesize \textbf{\sf Network number}}
		\rput[bl](0.1,4.56){\mbox{\footnotesize \textbf{{\sf Communication steps}}}}
		\rput[l](0.08,0.91){\scriptsize $\mathsf{10^{0}}$}
		\rput[l](0.08,1.75){\scriptsize $\mathsf{10^{1}}$}
		\rput[l](0.08,2.58){\scriptsize $\mathsf{10^{2}}$}
		\rput[l](0.08,3.41){\scriptsize $\mathsf{10^{3}}$}
		\rput[l](0.08,4.25){\scriptsize $\mathsf{10^{4}}$}

		\rput[t](0.610,0.61){\scriptsize $\mathsf{1}$}
		\rput[t](1.665,0.61){\scriptsize $\mathsf{2}$}
		\rput[t](2.737,0.61){\scriptsize $\mathsf{3}$}
		\rput[t](3.818,0.61){\scriptsize $\mathsf{4}$}
		\rput[t](4.892,0.61){\scriptsize $\mathsf{5}$}
		\rput[t](5.953,0.61){\scriptsize $\mathsf{6}$}
		\rput[t](7.030,0.61){\scriptsize $\mathsf{7}$}

		\rput[l](7.18,2.81){\scriptsize \textbf{\sf D-ADMM}}
		\rput[l](7.18,3.19){\scriptsize \textbf{\sf D-Lasso}}
		\rput[l](7.18,3.59){\scriptsize \textbf{\sf MM/NGS}}
		\rput[l](7.18,4.32){\scriptsize \textbf{\sf MM/DQA, DN,}}
		\rput[l](7.18,4.06){\scriptsize \textbf{\sf Subgradient}}
	      \end{pspicture}
            }
            
            \hspace{0.4cm}
            
            \isdraft{
	      \begin{pspicture}(0.48\linewidth,4.6cm) 
		\rput(0.225\linewidth,2.55){\includegraphics[scale=0.37]{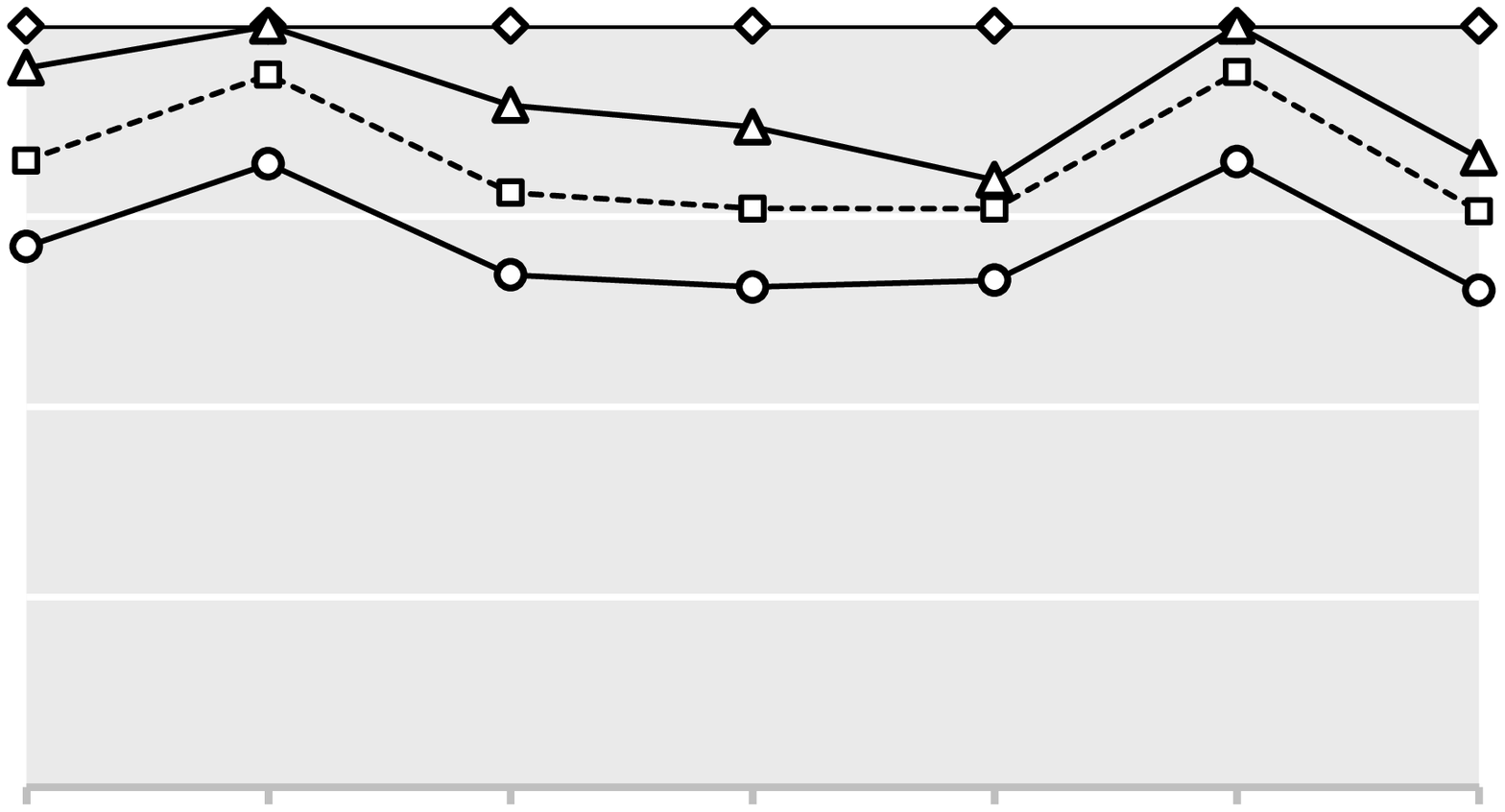}}
                \rput(0.225\linewidth,0.11){\footnotesize \textbf{\sf Network number}}
                \rput[bl](0.07,4.56){\footnotesize \textbf{{\sf Communication steps}}}
                \rput[l](0.05,0.91){\scriptsize $\mathsf{10^{0}}$}
                \rput[l](0.05,1.75){\scriptsize $\mathsf{10^{1}}$}
                \rput[l](0.05,2.58){\scriptsize $\mathsf{10^{2}}$}
                \rput[l](0.05,3.41){\scriptsize $\mathsf{10^{3}}$}
                \rput[l](0.05,4.25){\scriptsize $\mathsf{10^{4}}$}
                
                \rput[t](0.604,0.65){\scriptsize $\mathsf{1}$}
                \rput[t](1.630,0.65){\scriptsize $\mathsf{2}$}
                \rput[t](2.673,0.65){\scriptsize $\mathsf{3}$}
                \rput[t](3.721,0.65){\scriptsize $\mathsf{4}$}
                \rput[t](4.764,0.65){\scriptsize $\mathsf{5}$}
                \rput[t](5.800,0.65){\scriptsize $\mathsf{6}$}
                \rput[t](6.84,0.65){\scriptsize $\mathsf{7}$}
                
                \rput[l](7,3.03){\tiny \textbf{\sf D-ADMM}}
                \rput[l](7,3.38){\tiny \textbf{\sf D-Lasso}}
                \rput[l](7,3.66){\tiny \textbf{\sf MM/NGS}}
                \rput[l](7,4.32){\tiny \textbf{\sf MM/DQA, DN,}}
                \rput[l](7,4.06){\tiny \textbf{\sf Subgradient}}
              \end{pspicture}
            }{
	      \begin{pspicture}(0.48\linewidth,4.9cm) 
                \rput(0.21\linewidth,2.55){\includegraphics[scale=0.38]{figures/RPCen1Acc5.eps}}
                \rput(0.21\linewidth,0.11){\footnotesize \textbf{\sf Network number}}
                \rput[bl](0.1,4.56){\footnotesize \textbf{{\sf Communication steps}}}
                \rput[l](0.08,0.91){\scriptsize $\mathsf{10^{0}}$}
                \rput[l](0.08,1.75){\scriptsize $\mathsf{10^{1}}$}
                \rput[l](0.08,2.58){\scriptsize $\mathsf{10^{2}}$}
                \rput[l](0.08,3.41){\scriptsize $\mathsf{10^{3}}$}
                \rput[l](0.08,4.25){\scriptsize $\mathsf{10^{4}}$}

                \rput[t](0.610,0.61){\scriptsize $\mathsf{1}$}
                \rput[t](1.665,0.61){\scriptsize $\mathsf{2}$}
                \rput[t](2.737,0.61){\scriptsize $\mathsf{3}$}
                \rput[t](3.818,0.61){\scriptsize $\mathsf{4}$}
                \rput[t](4.892,0.61){\scriptsize $\mathsf{5}$}
                \rput[t](5.953,0.61){\scriptsize $\mathsf{6}$}
                \rput[t](7.030,0.61){\scriptsize $\mathsf{7}$}

                \rput[l](7.18,3.03){\scriptsize \textbf{\sf D-ADMM}}
                \rput[l](7.18,3.38){\scriptsize \textbf{\sf D-Lasso}}
                \rput[l](7.18,3.66){\scriptsize \textbf{\sf MM/NGS}}
                \rput[l](7.18,4.32){\scriptsize \textbf{\sf MM/DQA, DN,}}
                \rput[l](7.18,4.06){\scriptsize \textbf{\sf Subgradient}}
              \end{pspicture}
             }
        }
        \\
        \smallskip
        \subfigure[Scenario 2]{
            
            \isdraft{
	      \begin{pspicture}(0.48\linewidth,4.6cm)
		\rput(0.225\linewidth,2.55){\includegraphics[scale=0.37]{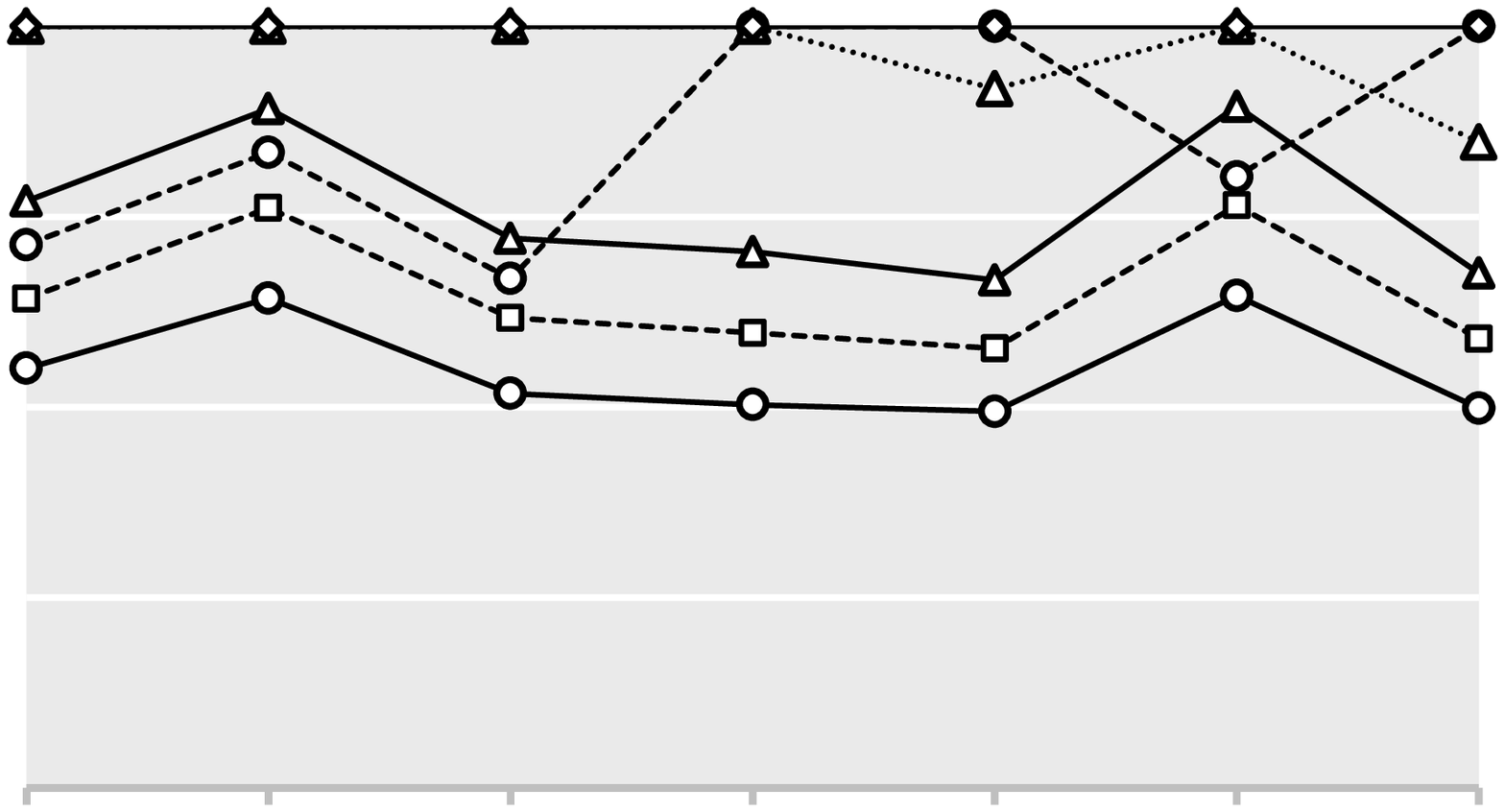}}
                \rput(0.225\linewidth,0.11){\footnotesize \textbf{\sf Network number}}
                \rput[bl](0.07,4.56){\footnotesize \textbf{{\sf Communication steps}}}
                \rput[l](0.05,0.91){\scriptsize $\mathsf{10^{0}}$}
                \rput[l](0.05,1.75){\scriptsize $\mathsf{10^{1}}$}
                \rput[l](0.05,2.58){\scriptsize $\mathsf{10^{2}}$}
                \rput[l](0.05,3.41){\scriptsize $\mathsf{10^{3}}$}
                \rput[l](0.05,4.25){\scriptsize $\mathsf{10^{4}}$}

                \rput[t](0.604,0.65){\scriptsize $\mathsf{1}$}
                \rput[t](1.630,0.65){\scriptsize $\mathsf{2}$}
                \rput[t](2.673,0.65){\scriptsize $\mathsf{3}$}
                \rput[t](3.721,0.65){\scriptsize $\mathsf{4}$}
                \rput[t](4.764,0.65){\scriptsize $\mathsf{5}$}
                \rput[t](5.800,0.65){\scriptsize $\mathsf{6}$}
                \rput[t](6.84,0.65){\scriptsize $\mathsf{7}$}

                \rput[l](7,2.49){\tiny \textbf{\sf D-ADMM}}
                \rput[l](7,2.81){\tiny \textbf{\sf D-Lasso}}
                \rput[l](7,3.11){\tiny \textbf{\sf MM/NGS}}
                \rput[l](7,3.93){\tiny \textbf{\sf DN}}
                \psline[linewidth=0.5pt]{-}(6.9,3.94)(6.67,4.02)
                \rput[l](7,3.62){\tiny \textbf{\sf MM/DQA}}
                \rput[l](7,4.23){\tiny \textbf{\sf Subgradient}}
              \end{pspicture} 
            }{
	      \begin{pspicture}(0.48\linewidth,4.9cm)
                \rput(0.21\linewidth,2.55){\includegraphics[scale=0.38]{figures/RPCen2Acc2.eps}}
                \rput(0.21\linewidth,0.11){\footnotesize \textbf{\sf Network number}}
                \rput[bl](0.1,4.56){\footnotesize \textbf{{\sf Communication steps}}}
                \rput[l](0.08,0.91){\scriptsize $\mathsf{10^{0}}$}
                \rput[l](0.08,1.75){\scriptsize $\mathsf{10^{1}}$}
                \rput[l](0.08,2.58){\scriptsize $\mathsf{10^{2}}$}
                \rput[l](0.08,3.41){\scriptsize $\mathsf{10^{3}}$}
                \rput[l](0.08,4.25){\scriptsize $\mathsf{10^{4}}$}

                \rput[t](0.610,0.61){\scriptsize $\mathsf{1}$}
                \rput[t](1.665,0.61){\scriptsize $\mathsf{2}$}
                \rput[t](2.737,0.61){\scriptsize $\mathsf{3}$}
                \rput[t](3.818,0.61){\scriptsize $\mathsf{4}$}
                \rput[t](4.892,0.61){\scriptsize $\mathsf{5}$}
                \rput[t](5.953,0.61){\scriptsize $\mathsf{6}$}
                \rput[t](7.030,0.61){\scriptsize $\mathsf{7}$}

                \rput[l](7.18,2.49){\scriptsize \textbf{\sf D-ADMM}}
                \rput[l](7.18,2.81){\scriptsize \textbf{\sf D-Lasso}}
                \rput[l](7.18,3.11){\scriptsize \textbf{\sf MM/NGS}}
                \rput[l](7.18,3.95){\scriptsize \textbf{\sf DN}}
                \psline[linewidth=0.5pt]{-}(7.1,3.96)(6.87,4.07)
                \rput[l](7.18,3.67){\scriptsize \textbf{\sf MM/DQA}}
                \rput[l](7.18,4.23){\scriptsize \textbf{\sf Subgradient}}
              \end{pspicture}
            }
            
        \hspace{0.4cm}
            
	    \isdraft{ 
	      \begin{pspicture}(0.48\linewidth,4.6cm)
	        \rput(0.225\linewidth,2.55){\includegraphics[scale=0.37]{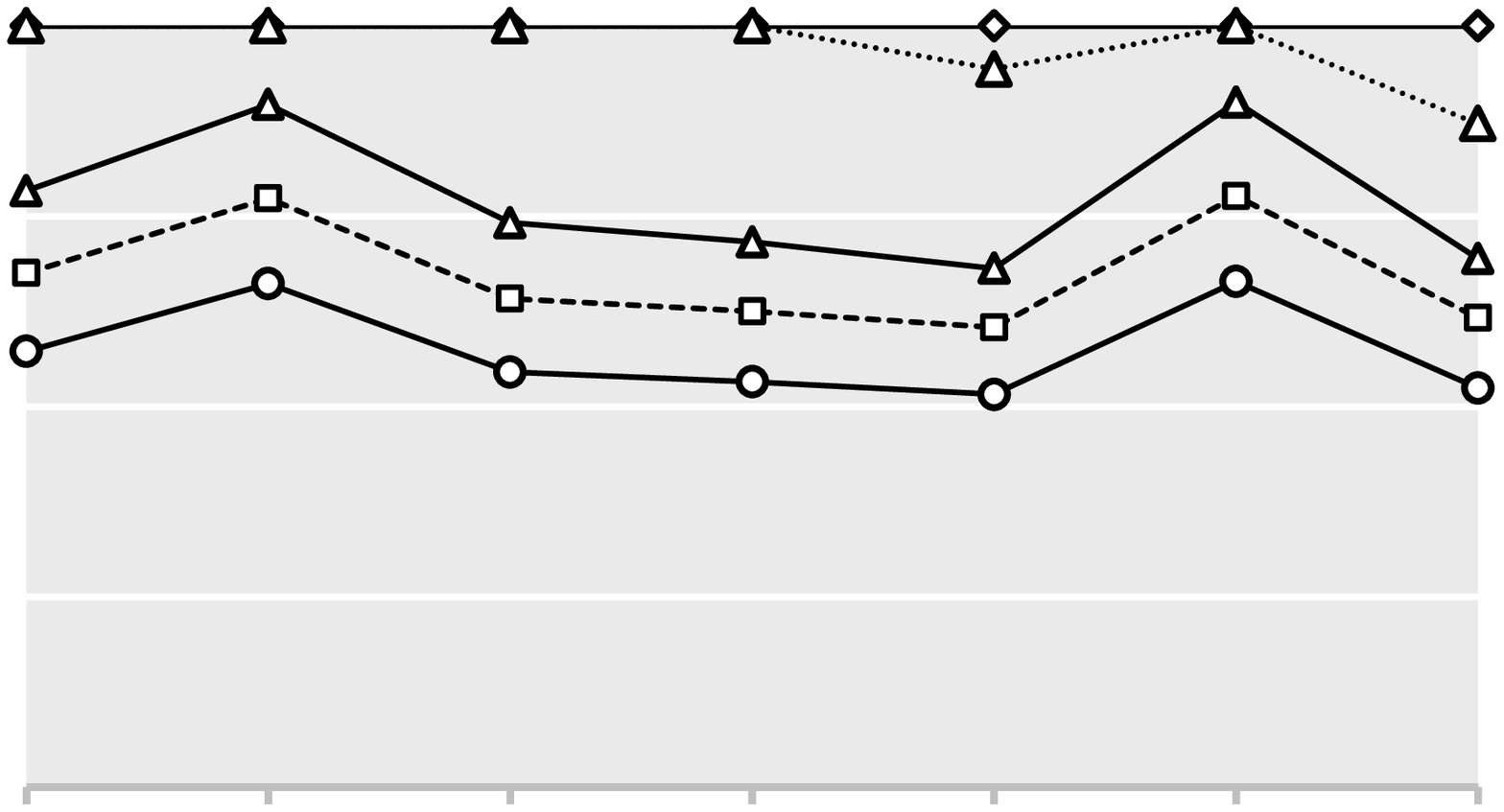}}
                \rput(0.225\linewidth,0.11){\footnotesize \textbf{\sf Network number}}
                \rput[bl](0.07,4.56){\footnotesize \textbf{{\sf Communication steps}}}
                \rput[l](0.05,0.91){\scriptsize $\mathsf{10^{0}}$}
                \rput[l](0.05,1.75){\scriptsize $\mathsf{10^{1}}$}
                \rput[l](0.05,2.58){\scriptsize $\mathsf{10^{2}}$}
                \rput[l](0.05,3.41){\scriptsize $\mathsf{10^{3}}$}
                \rput[l](0.05,4.25){\scriptsize $\mathsf{10^{4}}$}

                \rput[t](0.604,0.65){\scriptsize $\mathsf{1}$}
                \rput[t](1.630,0.65){\scriptsize $\mathsf{2}$}
                \rput[t](2.673,0.65){\scriptsize $\mathsf{3}$}
                \rput[t](3.721,0.65){\scriptsize $\mathsf{4}$}
                \rput[t](4.764,0.65){\scriptsize $\mathsf{5}$}
                \rput[t](5.800,0.65){\scriptsize $\mathsf{6}$}
                \rput[t](6.84,0.65){\scriptsize $\mathsf{7}$}

                \rput[l](7,2.57){\tiny \textbf{\sf D-ADMM}}
                \rput[l](7,2.93){\tiny \textbf{\sf D-Lasso}}
                \rput[l](7,3.22){\tiny \textbf{\sf MM/NGS}}
                \rput[l](7,3.73){\tiny \textbf{\sf MM/DQA}}
                \rput[l](7,4.32){\tiny \textbf{\sf DN,}}
                \rput[l](7,4.06){\tiny \textbf{\sf Subgradient}}
              \end{pspicture}
	    }{
	      \begin{pspicture}(0.48\linewidth,4.9cm)
                \rput(0.21\linewidth,2.55){\includegraphics[scale=0.38]{figures/RPCen2Acc5.eps}}
                \rput(0.21\linewidth,0.11){\footnotesize \textbf{\sf Network number}}
                \rput[bl](0.1,4.56){\footnotesize \textbf{{\sf Communication steps}}}
                \rput[l](0.08,0.91){\scriptsize $\mathsf{10^{0}}$}
                \rput[l](0.08,1.75){\scriptsize $\mathsf{10^{1}}$}
                \rput[l](0.08,2.58){\scriptsize $\mathsf{10^{2}}$}
                \rput[l](0.08,3.41){\scriptsize $\mathsf{10^{3}}$}
                \rput[l](0.08,4.25){\scriptsize $\mathsf{10^{4}}$}

                \rput[t](0.610,0.61){\scriptsize $\mathsf{1}$}
                \rput[t](1.665,0.61){\scriptsize $\mathsf{2}$}
                \rput[t](2.737,0.61){\scriptsize $\mathsf{3}$}
                \rput[t](3.818,0.61){\scriptsize $\mathsf{4}$}
                \rput[t](4.892,0.61){\scriptsize $\mathsf{5}$}
                \rput[t](5.953,0.61){\scriptsize $\mathsf{6}$}
                \rput[t](7.030,0.61){\scriptsize $\mathsf{7}$}

                \rput[l](7.18,2.57){\scriptsize \textbf{\sf D-ADMM}}
                \rput[l](7.18,2.93){\scriptsize \textbf{\sf D-Lasso}}
                \rput[l](7.18,3.22){\scriptsize \textbf{\sf MM/NGS}}
                \rput[l](7.18,3.76){\scriptsize \textbf{\sf MM/DQA}}
                \rput[l](7.18,4.32){\scriptsize \textbf{\sf DN,}}
                \rput[l](7.18,4.06){\scriptsize \textbf{\sf Subgradient}}
              \end{pspicture}
            }
        }
        \\
        \smallskip
        \subfigure[Scenario 3]{\label{SubFig:RPExperimentsScen3}
        
            \isdraft{
	      \begin{pspicture}(0.48\linewidth,4.6cm)
		\rput(0.225\linewidth,2.55){\includegraphics[scale=0.37]{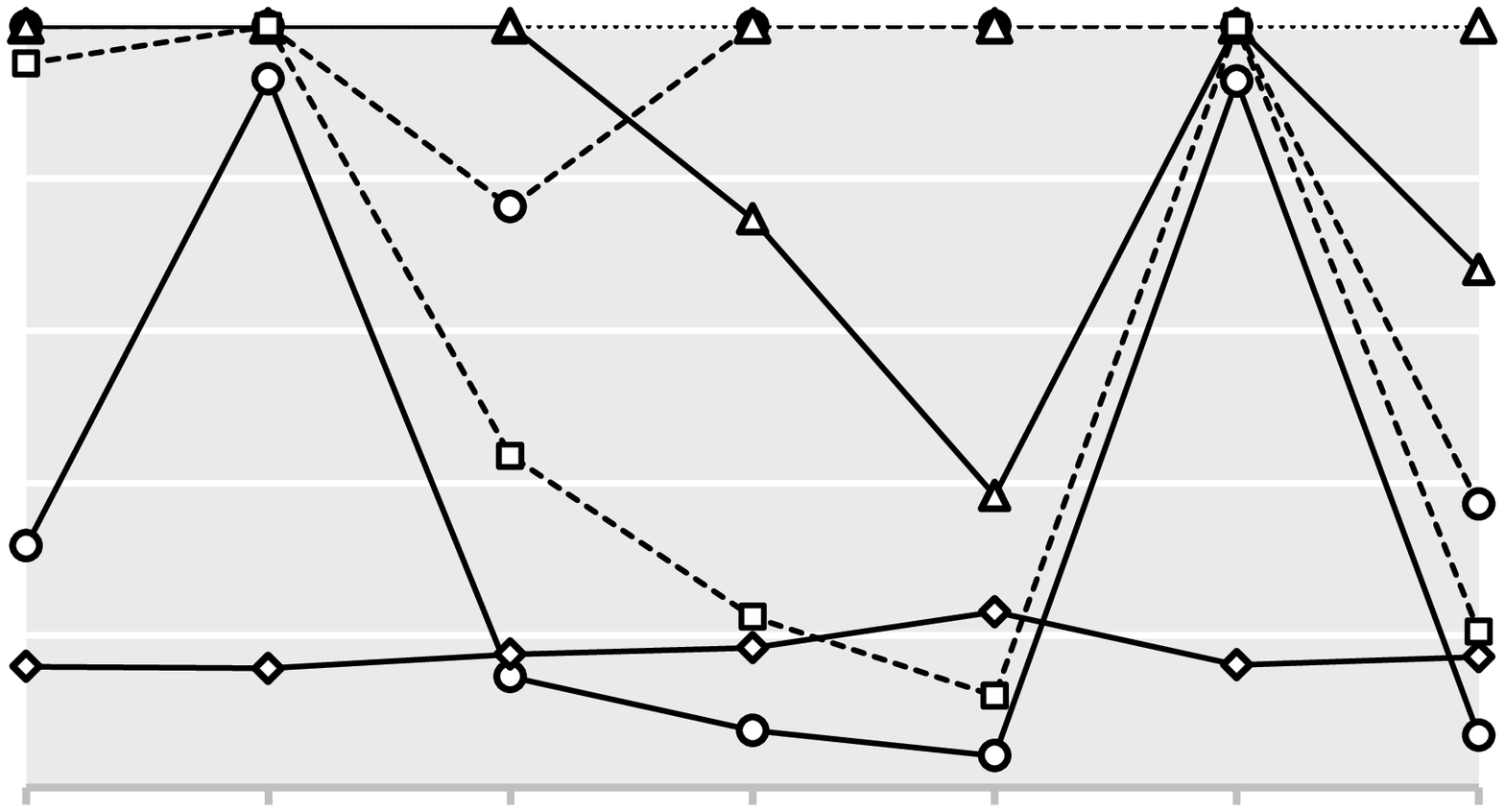}}
		\rput(0.225\linewidth,0.11){\footnotesize \textbf{\sf Network number}}
		\rput[bl](0.17,4.65){\footnotesize \textbf{{\sf Communication steps}}}
                \rput[r](0.42,0.88){\scriptsize $\mathsf{0}$}
                \rput[r](0.42,1.547){\scriptsize $\mathsf{2}$}
                \rput[r](0.42,2.22){\scriptsize $\mathsf{4}$}
                \rput[r](0.42,2.90){\scriptsize $\mathsf{6}$}
                \rput[r](0.42,3.57){\scriptsize $\mathsf{8}$}
                \rput[r](0.42,4.20){\scriptsize $\mathsf{10}$}

                \rput[t](0.604,0.65){\scriptsize $\mathsf{1}$}
                \rput[t](1.630,0.65){\scriptsize $\mathsf{2}$}
                \rput[t](2.673,0.65){\scriptsize $\mathsf{3}$}
                \rput[t](3.721,0.65){\scriptsize $\mathsf{4}$}
                \rput[t](4.764,0.65){\scriptsize $\mathsf{5}$}
                \rput[t](5.800,0.65){\scriptsize $\mathsf{6}$}
                \rput[t](6.84,0.65){\scriptsize $\mathsf{7}$}

                \rput[lb](0.58,4.35){\tiny $\mathsf{\times 10^{3}}$}

                \rput[l](7,1.08){\tiny \textbf{\sf D-ADMM}}
                \rput[l](7,1.33){\tiny \textbf{\sf Subgradient}}
                \rput[l](7,1.60){\tiny \textbf{\sf D-Lasso}}
                \rput[l](7,2.12){\tiny \textbf{\sf DN}}
                \rput[l](7,3.15){\tiny \textbf{\sf MM/NGS}}
                \rput[l](7,4.22){\tiny \textbf{\sf MM/DQA}}
              \end{pspicture}
            }{
              \begin{pspicture}(0.48\linewidth,5.0cm)
                \rput(0.21\linewidth,2.55){\includegraphics[scale=0.38]{figures/RPCen3Acc2.eps}}
		\rput(0.21\linewidth,0.11){\footnotesize \textbf{\sf Network number}}
		\rput[bl](0.17,4.74){\footnotesize \textbf{{\sf Communication steps}}}
                \rput[r](0.42,0.88){\scriptsize $\mathsf{0}$}
                \rput[r](0.42,1.547){\scriptsize $\mathsf{2}$}
                \rput[r](0.42,2.22){\scriptsize $\mathsf{4}$}
                \rput[r](0.42,2.90){\scriptsize $\mathsf{6}$}
                \rput[r](0.42,3.57){\scriptsize $\mathsf{8}$}
                \rput[r](0.42,4.20){\scriptsize $\mathsf{10}$}

                \rput[t](0.610,0.61){\scriptsize $\mathsf{1}$}
                \rput[t](1.665,0.61){\scriptsize $\mathsf{2}$}
                \rput[t](2.737,0.61){\scriptsize $\mathsf{3}$}
                \rput[t](3.818,0.61){\scriptsize $\mathsf{4}$}
                \rput[t](4.892,0.61){\scriptsize $\mathsf{5}$}
                \rput[t](5.953,0.61){\scriptsize $\mathsf{6}$}
                \rput[t](7.030,0.61){\scriptsize $\mathsf{7}$}

                \rput[lb](0.58,4.43){\tiny $\mathsf{\times 10^{3}}$}

                \rput[l](7.18,1.08){\scriptsize \textbf{\sf D-ADMM}}
                \rput[l](7.18,1.33){\scriptsize \textbf{\sf Subgradient}}
                \rput[l](7.18,1.60){\scriptsize \textbf{\sf D-Lasso}}
                \rput[l](7.18,2.12){\scriptsize \textbf{\sf DN}}
                \rput[l](7.18,3.15){\scriptsize \textbf{\sf MM/NGS}}
                \rput[l](7.18,4.22){\scriptsize \textbf{\sf MM/DQA}}
              \end{pspicture}
            }    
            
            \hspace{0.4cm}
            
	    \isdraft{
	      \begin{pspicture}(0.48\linewidth,4.6cm)
		\rput(0.225\linewidth,2.55){\includegraphics[scale=0.37]{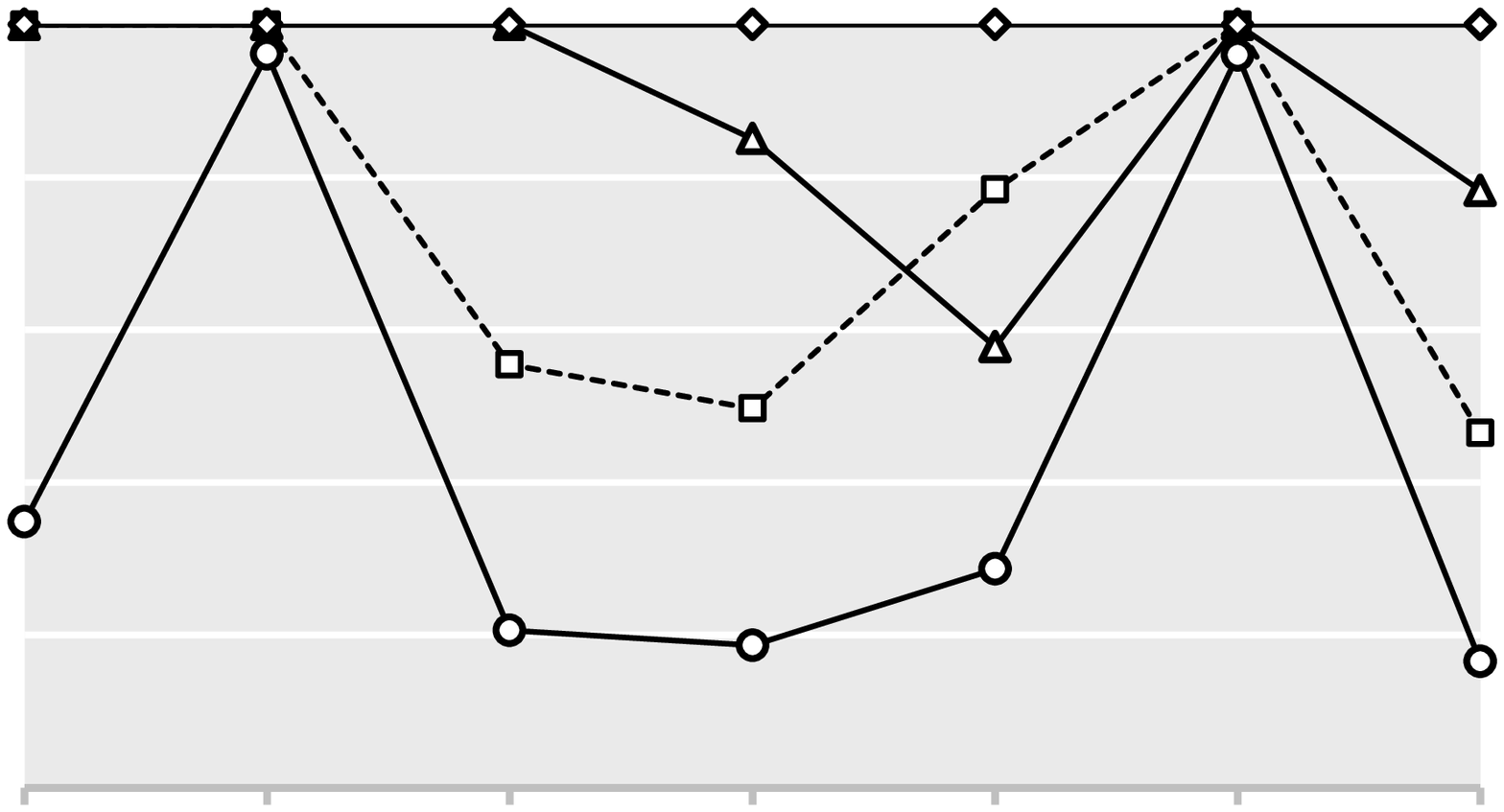}}
		\rput(0.225\linewidth,0.11){\footnotesize \textbf{\sf Network number}}
		\rput[bl](0.17,4.65){\footnotesize \textbf{{\sf Communication steps}}}
		\rput[r](0.42,0.88){\scriptsize $\mathsf{0}$}
                \rput[r](0.42,1.547){\scriptsize $\mathsf{2}$}
                \rput[r](0.42,2.22){\scriptsize $\mathsf{4}$}
                \rput[r](0.42,2.90){\scriptsize $\mathsf{6}$}
                \rput[r](0.42,3.57){\scriptsize $\mathsf{8}$}
                \rput[r](0.42,4.20){\scriptsize $\mathsf{10}$}

                \rput[t](0.604,0.65){\scriptsize $\mathsf{1}$}
                \rput[t](1.630,0.65){\scriptsize $\mathsf{2}$}
                \rput[t](2.673,0.65){\scriptsize $\mathsf{3}$}
                \rput[t](3.721,0.65){\scriptsize $\mathsf{4}$}
                \rput[t](4.764,0.65){\scriptsize $\mathsf{5}$}
                \rput[t](5.800,0.65){\scriptsize $\mathsf{6}$}
                \rput[t](6.84,0.65){\scriptsize $\mathsf{7}$}

                \rput[lb](0.58,4.35){\tiny $\mathsf{\times 10^{3}}$}

                \rput[l](7,1.43){\tiny \textbf{\sf D-ADMM}}
                \rput[l](7,2.44){\tiny \textbf{\sf D-Lasso}}
                \rput[l](7,3.5){\tiny \textbf{\sf MM/NGS}}
		\rput[l](7,4.32){\tiny \textbf{\sf MM/DQA, DN,}}
                \rput[l](7,4.06){\tiny \textbf{\sf Subgradient}}
              \end{pspicture}
	    }{
	      \begin{pspicture}(0.48\linewidth,5.0cm)
                \rput(0.21\linewidth,2.55){\includegraphics[scale=0.38]{figures/RPCen3Acc5.eps}}
		\rput(0.21\linewidth,0.11){\footnotesize \textbf{\sf Network number}}
		\rput[bl](0.17,4.74){\footnotesize \textbf{{\sf Communication steps}}}
                \rput[r](0.42,0.88){\scriptsize $\mathsf{0}$}
                \rput[r](0.42,1.547){\scriptsize $\mathsf{2}$}
                \rput[r](0.42,2.22){\scriptsize $\mathsf{4}$}
                \rput[r](0.42,2.90){\scriptsize $\mathsf{6}$}
                \rput[r](0.42,3.57){\scriptsize $\mathsf{8}$}
                \rput[r](0.42,4.20){\scriptsize $\mathsf{10}$}

                \rput[t](0.610,0.61){\scriptsize $\mathsf{1}$}
                \rput[t](1.665,0.61){\scriptsize $\mathsf{2}$}
                \rput[t](2.737,0.61){\scriptsize $\mathsf{3}$}
                \rput[t](3.818,0.61){\scriptsize $\mathsf{4}$}
                \rput[t](4.892,0.61){\scriptsize $\mathsf{5}$}
                \rput[t](5.953,0.61){\scriptsize $\mathsf{6}$}
                \rput[t](7.030,0.61){\scriptsize $\mathsf{7}$}

                \rput[lb](0.58,4.43){\tiny $\mathsf{\times 10^{3}}$}

                \rput[l](7.18,1.43){\scriptsize \textbf{\sf D-ADMM}}
                \rput[l](7.18,2.44){\scriptsize \textbf{\sf D-Lasso}}
                \rput[l](7.18,3.5){\scriptsize \textbf{\sf MM/NGS}}
		\rput[l](7.18,4.32){\scriptsize \textbf{\sf MM/DQA, DN,}}
                \rput[l](7.18,4.06){\scriptsize \textbf{\sf Subgradient}}
              \end{pspicture}
            }
            
        }
        \\
        \smallskip
        \subfigure[Scenario 4]{
            
            \isdraft{
	      \begin{pspicture}(0.48\linewidth,4.6cm)
		\rput(0.225\linewidth,2.55){\includegraphics[scale=0.37]{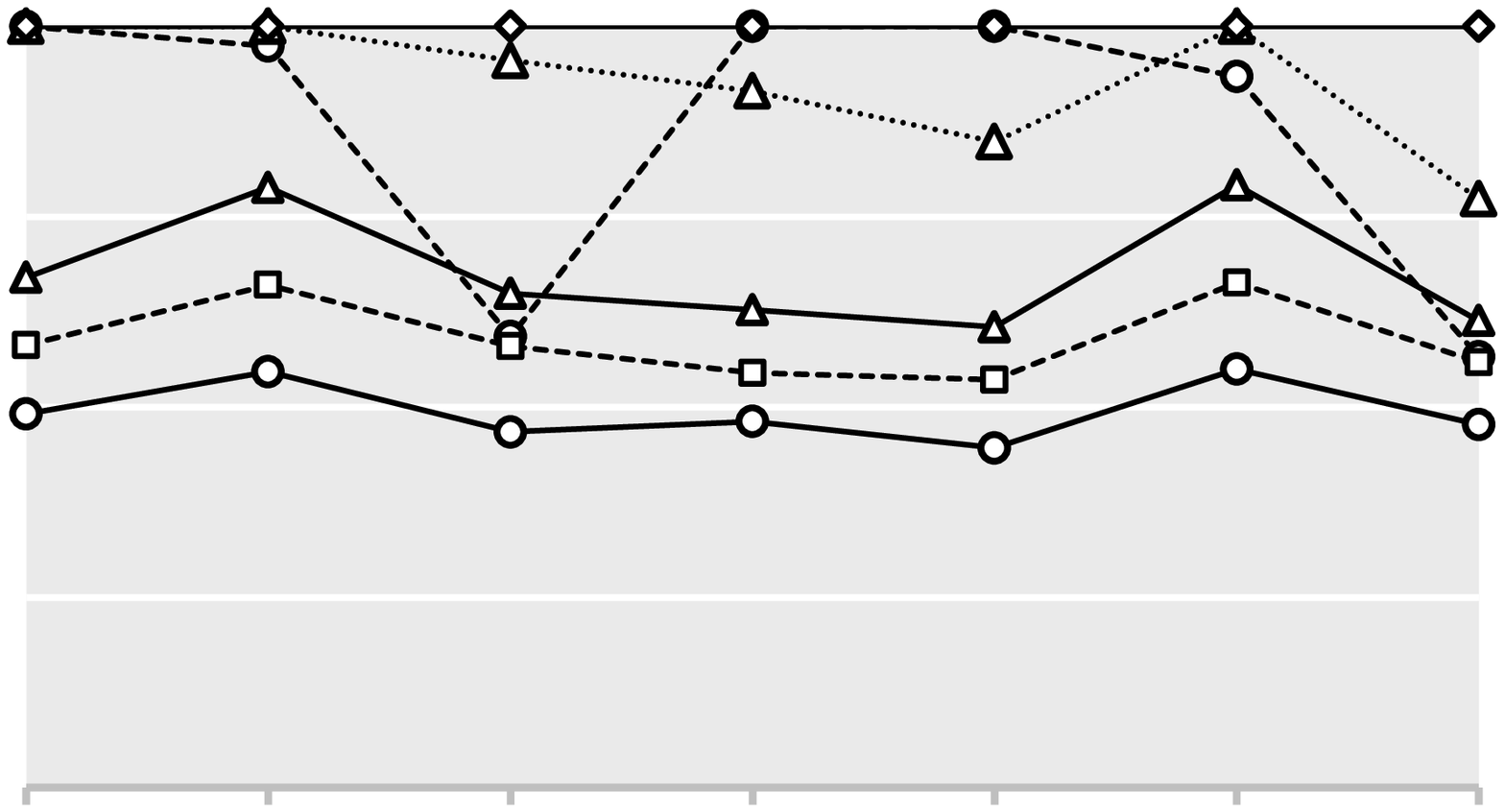}}
                \rput(0.225\linewidth,0.11){\footnotesize \textbf{\sf Network number}}
                \rput[bl](0.07,4.56){\footnotesize \textbf{{\sf Communication steps}}}
                \rput[l](0.05,0.91){\scriptsize $\mathsf{10^{0}}$}
                \rput[l](0.05,1.75){\scriptsize $\mathsf{10^{1}}$}
                \rput[l](0.05,2.58){\scriptsize $\mathsf{10^{2}}$}
                \rput[l](0.05,3.41){\scriptsize $\mathsf{10^{3}}$}
                \rput[l](0.05,4.25){\scriptsize $\mathsf{10^{4}}$}

                \rput[t](0.604,0.65){\scriptsize $\mathsf{1}$}
                \rput[t](1.630,0.65){\scriptsize $\mathsf{2}$}
                \rput[t](2.673,0.65){\scriptsize $\mathsf{3}$}
                \rput[t](3.721,0.65){\scriptsize $\mathsf{4}$}
                \rput[t](4.764,0.65){\scriptsize $\mathsf{5}$}
                \rput[t](5.800,0.65){\scriptsize $\mathsf{6}$}
                \rput[t](6.84,0.65){\scriptsize $\mathsf{7}$}

                \rput[l](7,2.00){\tiny \textbf{\sf D-ADMM}}
                \psline[linewidth=0.5pt]{-}(6.95,2.05)(6.5,2.5)
                \rput[l](7,2.55){\tiny \textbf{\sf D-Lasso}}
                \psline[linewidth=0.5pt]{-}(6.95,2.58)(6.5,2.8)
                \rput[l](7,2.95){\tiny \textbf{\sf MM/NGS}}
                \rput[l](7,3.44){\tiny \textbf{\sf MM/DQA}}
                \rput[l](7,2.75){\tiny \textbf{\sf DN}}
                \rput[l](7,4.22){\tiny \textbf{\sf Subgradient}}
              \end{pspicture}
            }{
	      \begin{pspicture}(0.48\linewidth,4.9cm)
                \rput(0.21\linewidth,2.55){\includegraphics[scale=0.38]{figures/RPCen4Acc2.eps}}
                \rput(0.21\linewidth,0.11){\footnotesize \textbf{\sf Network number}}
                \rput[bl](0.1,4.56){\footnotesize \textbf{{\sf Communication steps}}}
                \rput[l](0.08,0.91){\scriptsize $\mathsf{10^{0}}$}
                \rput[l](0.08,1.75){\scriptsize $\mathsf{10^{1}}$}
                \rput[l](0.08,2.58){\scriptsize $\mathsf{10^{2}}$}
                \rput[l](0.08,3.41){\scriptsize $\mathsf{10^{3}}$}
                \rput[l](0.08,4.25){\scriptsize $\mathsf{10^{4}}$}

                \rput[t](0.610,0.61){\scriptsize $\mathsf{1}$}
                \rput[t](1.665,0.61){\scriptsize $\mathsf{2}$}
                \rput[t](2.737,0.61){\scriptsize $\mathsf{3}$}
                \rput[t](3.818,0.61){\scriptsize $\mathsf{4}$}
                \rput[t](4.892,0.61){\scriptsize $\mathsf{5}$}
                \rput[t](5.953,0.61){\scriptsize $\mathsf{6}$}
                \rput[t](7.030,0.61){\scriptsize $\mathsf{7}$}

                \rput[l](7.18,2.1){\scriptsize \textbf{\sf D-ADMM}}
                \psline[linewidth=0.5pt]{-}(7.17,2.15)(6.7,2.5)
                \rput[l](7.18,2.53){\scriptsize \textbf{\sf D-Lasso}}
                \psline[linewidth=0.5pt]{-}(7.17,2.57)(6.7,2.8)
                \rput[l](7.18,2.99){\scriptsize \textbf{\sf MM/NGS}}
                \rput[l](7.18,3.46){\scriptsize \textbf{\sf MM/DQA}}
                \rput[l](7.18,2.76){\scriptsize \textbf{\sf DN}}
                \rput[l](7.18,4.22){\scriptsize \textbf{\sf Subgradient}}
              \end{pspicture}
            }
            
            \hspace{0.4cm}
            
	    \isdraft{
	      \begin{pspicture}(0.48\linewidth,4.6cm)
		\rput(0.225\linewidth,2.55){\includegraphics[scale=0.37]{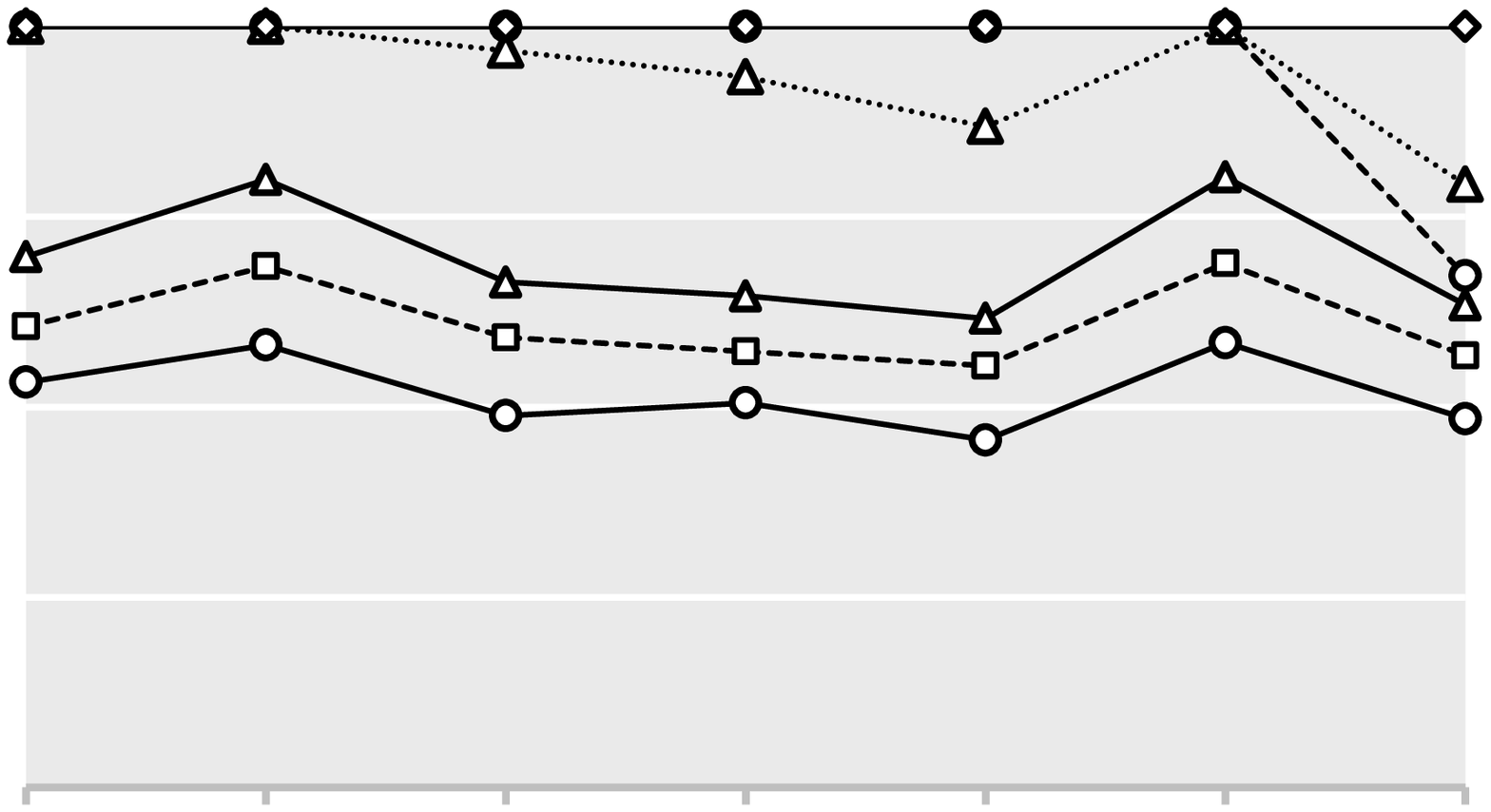}}
                \rput(0.225\linewidth,0.11){\footnotesize \textbf{\sf Network number}}
                \rput[bl](0.07,4.56){\footnotesize \textbf{{\sf Communication steps}}}
                \rput[l](0.05,0.91){\scriptsize $\mathsf{10^{0}}$}
                \rput[l](0.05,1.75){\scriptsize $\mathsf{10^{1}}$}
                \rput[l](0.05,2.58){\scriptsize $\mathsf{10^{2}}$}
                \rput[l](0.05,3.41){\scriptsize $\mathsf{10^{3}}$}
                \rput[l](0.05,4.25){\scriptsize $\mathsf{10^{4}}$}

                \rput[t](0.604,0.65){\scriptsize $\mathsf{1}$}
                \rput[t](1.630,0.65){\scriptsize $\mathsf{2}$}
                \rput[t](2.673,0.65){\scriptsize $\mathsf{3}$}
                \rput[t](3.721,0.65){\scriptsize $\mathsf{4}$}
                \rput[t](4.764,0.65){\scriptsize $\mathsf{5}$}
                \rput[t](5.800,0.65){\scriptsize $\mathsf{6}$}
                \rput[t](6.84,0.65){\scriptsize $\mathsf{7}$}

                \rput[l](7,2.46){\tiny \textbf{\sf D-ADMM}}
                \rput[l](7,2.74){\tiny \textbf{\sf D-Lasso}}
                \rput[l](7,2.96){\tiny \textbf{\sf MM/NGS}}
                \rput[l](7,3.5){\tiny \textbf{\sf MM/DQA}}
                \rput[l](7,3.14){\tiny \textbf{\sf DN}}
                \rput[l](7,4.22){\tiny \textbf{\sf Subgradient}}
              \end{pspicture}
	    }{
	      \begin{pspicture}(0.48\linewidth,4.9cm)
                \rput(0.21\linewidth,2.55){\includegraphics[scale=0.38]{figures/RPCen4Acc5.eps}}
                \rput(0.21\linewidth,0.11){\footnotesize \textbf{\sf Network number}}
                \rput[bl](0.1,4.56){\footnotesize \textbf{{\sf Communication steps}}}
                \rput[l](0.08,0.91){\scriptsize $\mathsf{10^{0}}$}
                \rput[l](0.08,1.75){\scriptsize $\mathsf{10^{1}}$}
                \rput[l](0.08,2.58){\scriptsize $\mathsf{10^{2}}$}
                \rput[l](0.08,3.41){\scriptsize $\mathsf{10^{3}}$}
                \rput[l](0.08,4.25){\scriptsize $\mathsf{10^{4}}$}

                \rput[t](0.610,0.61){\scriptsize $\mathsf{1}$}
                \rput[t](1.665,0.61){\scriptsize $\mathsf{2}$}
                \rput[t](2.737,0.61){\scriptsize $\mathsf{3}$}
                \rput[t](3.818,0.61){\scriptsize $\mathsf{4}$}
                \rput[t](4.892,0.61){\scriptsize $\mathsf{5}$}
                \rput[t](5.953,0.61){\scriptsize $\mathsf{6}$}
                \rput[t](7.030,0.61){\scriptsize $\mathsf{7}$}

                \rput[l](7.18,2.45){\scriptsize \textbf{\sf D-ADMM}}
                \rput[l](7.18,2.72){\scriptsize \textbf{\sf D-Lasso}}
                \rput[l](7.18,2.95){\scriptsize \textbf{\sf MM/NGS}}
                \rput[l](7.18,3.5){\scriptsize \textbf{\sf MM/DQA}}
                \rput[l](7.18,3.15){\scriptsize \textbf{\sf DN}}
                \rput[l](7.18,4.22){\scriptsize \textbf{\sf Subgradient}}
               \end{pspicture}
            }
                            
        }
        \isdraft{\vspace{-0.3cm}}{}
        \caption[]{
        	Type I experiments: number of communication steps to reach accuracies of~$1\%$ and~$10^{-3}\%$ as a function of the network (see Table~\ref{Tab:Networks}).
        }
    \end{figure*}
    \begin{figure*}
        \centering
        \subfigure[Scenario 5]{
	  \isdraft{
	    \begin{pspicture}(0.48\linewidth,4.6cm)
                \rput(0.225\linewidth,2.55){\includegraphics[scale=0.37]{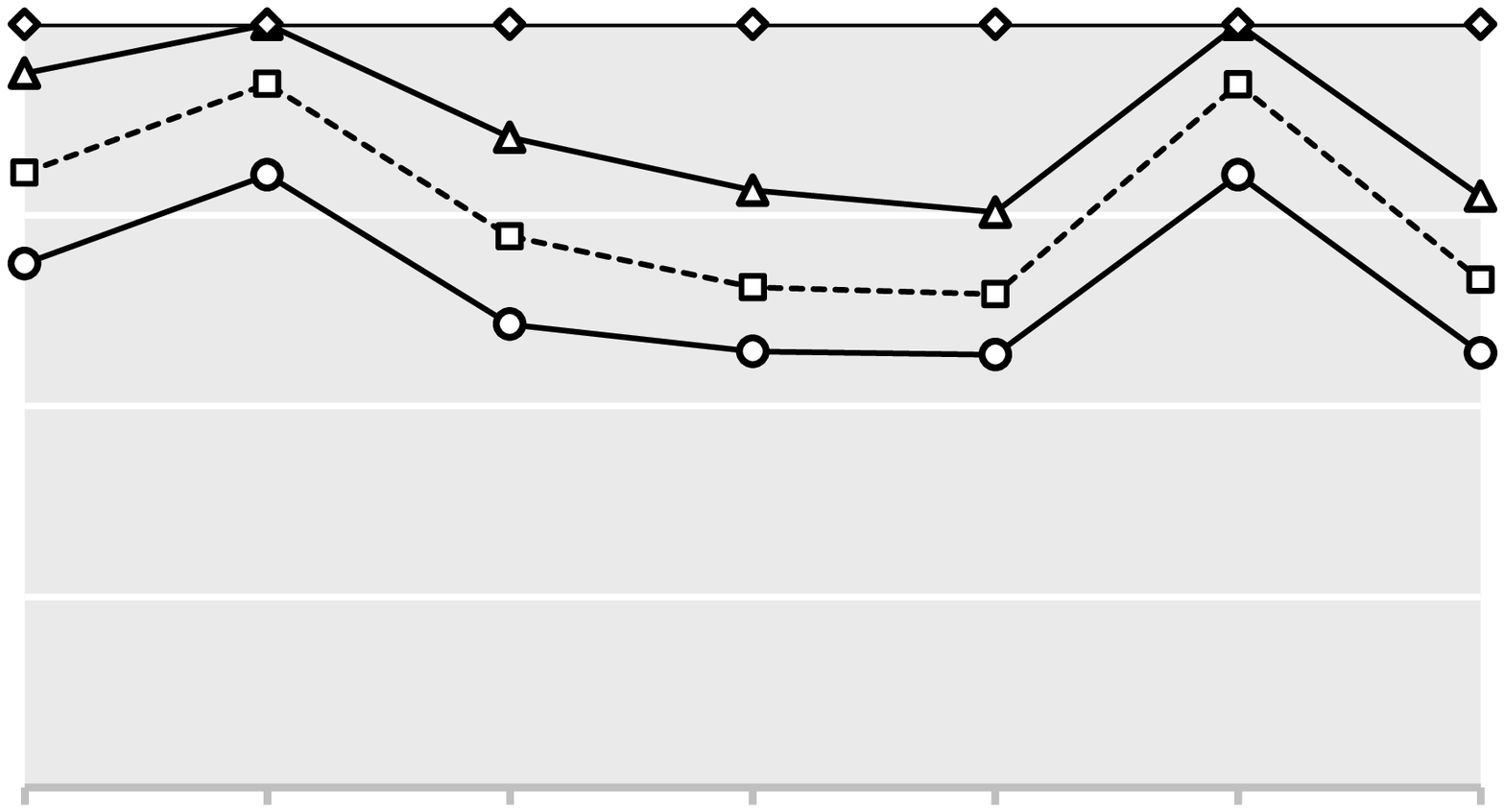}}
                \rput(0.225\linewidth,0.11){\footnotesize \textbf{\sf Network number}}
                \rput[bl](0.07,4.56){\footnotesize \textbf{{\sf Communication steps}}}
                \rput[l](0.05,0.91){\scriptsize $\mathsf{10^{0}}$}
                \rput[l](0.05,1.75){\scriptsize $\mathsf{10^{1}}$}
                \rput[l](0.05,2.58){\scriptsize $\mathsf{10^{2}}$}
                \rput[l](0.05,3.41){\scriptsize $\mathsf{10^{3}}$}
                \rput[l](0.05,4.25){\scriptsize $\mathsf{10^{4}}$}

                \rput[t](0.604,0.65){\scriptsize $\mathsf{1}$}
                \rput[t](1.630,0.65){\scriptsize $\mathsf{2}$}
                \rput[t](2.673,0.65){\scriptsize $\mathsf{3}$}
                \rput[t](3.721,0.65){\scriptsize $\mathsf{4}$}
                \rput[t](4.764,0.65){\scriptsize $\mathsf{5}$}
                \rput[t](5.800,0.65){\scriptsize $\mathsf{6}$}
                \rput[t](6.84,0.65){\scriptsize $\mathsf{7}$}

                \rput[l](7,2.76){\tiny \textbf{\sf D-ADMM}}
                \rput[l](7,3.12){\tiny \textbf{\sf D-Lasso}}
                \rput[l](7,3.47){\tiny \textbf{\sf MM/NGS}}
		\rput[l](7,4.32){\tiny \textbf{\sf MM/DQA, DN,}}
                \rput[l](7,4.06){\tiny \textbf{\sf Subgradient}}
            \end{pspicture}
	  }{
            \begin{pspicture}(0.48\linewidth,4.9cm)
                \rput(0.21\linewidth,2.55){\includegraphics[scale=0.38]{figures/RPCen5Acc2.eps}}
                \rput(0.21\linewidth,0.11){\footnotesize \textbf{\sf Network number}}
                \rput[bl](0.1,4.56){\footnotesize \textbf{{\sf Communication steps}}}
                \rput[l](0.08,0.91){\scriptsize $\mathsf{10^{0}}$}
                \rput[l](0.08,1.75){\scriptsize $\mathsf{10^{1}}$}
                \rput[l](0.08,2.58){\scriptsize $\mathsf{10^{2}}$}
                \rput[l](0.08,3.41){\scriptsize $\mathsf{10^{3}}$}
                \rput[l](0.08,4.25){\scriptsize $\mathsf{10^{4}}$}

                \rput[t](0.610,0.61){\scriptsize $\mathsf{1}$}
                \rput[t](1.665,0.61){\scriptsize $\mathsf{2}$}
                \rput[t](2.737,0.61){\scriptsize $\mathsf{3}$}
                \rput[t](3.818,0.61){\scriptsize $\mathsf{4}$}
                \rput[t](4.892,0.61){\scriptsize $\mathsf{5}$}
                \rput[t](5.953,0.61){\scriptsize $\mathsf{6}$}
                \rput[t](7.030,0.61){\scriptsize $\mathsf{7}$}

                \rput[l](7.18,2.76){\scriptsize \textbf{\sf D-ADMM}}
                \rput[l](7.18,3.12){\scriptsize \textbf{\sf D-Lasso}}
                \rput[l](7.18,3.47){\scriptsize \textbf{\sf MM/NGS}}
		\rput[l](7.18,4.32){\scriptsize \textbf{\sf MM/DQA, DN,}}
                \rput[l](7.18,4.06){\scriptsize \textbf{\sf Subgradient}}
            \end{pspicture}
          }
        \hspace{0.4cm}
	  \isdraft{
	    \begin{pspicture}(0.48\linewidth,4.6cm)
                \rput(0.225\linewidth,2.55){\includegraphics[scale=0.37]{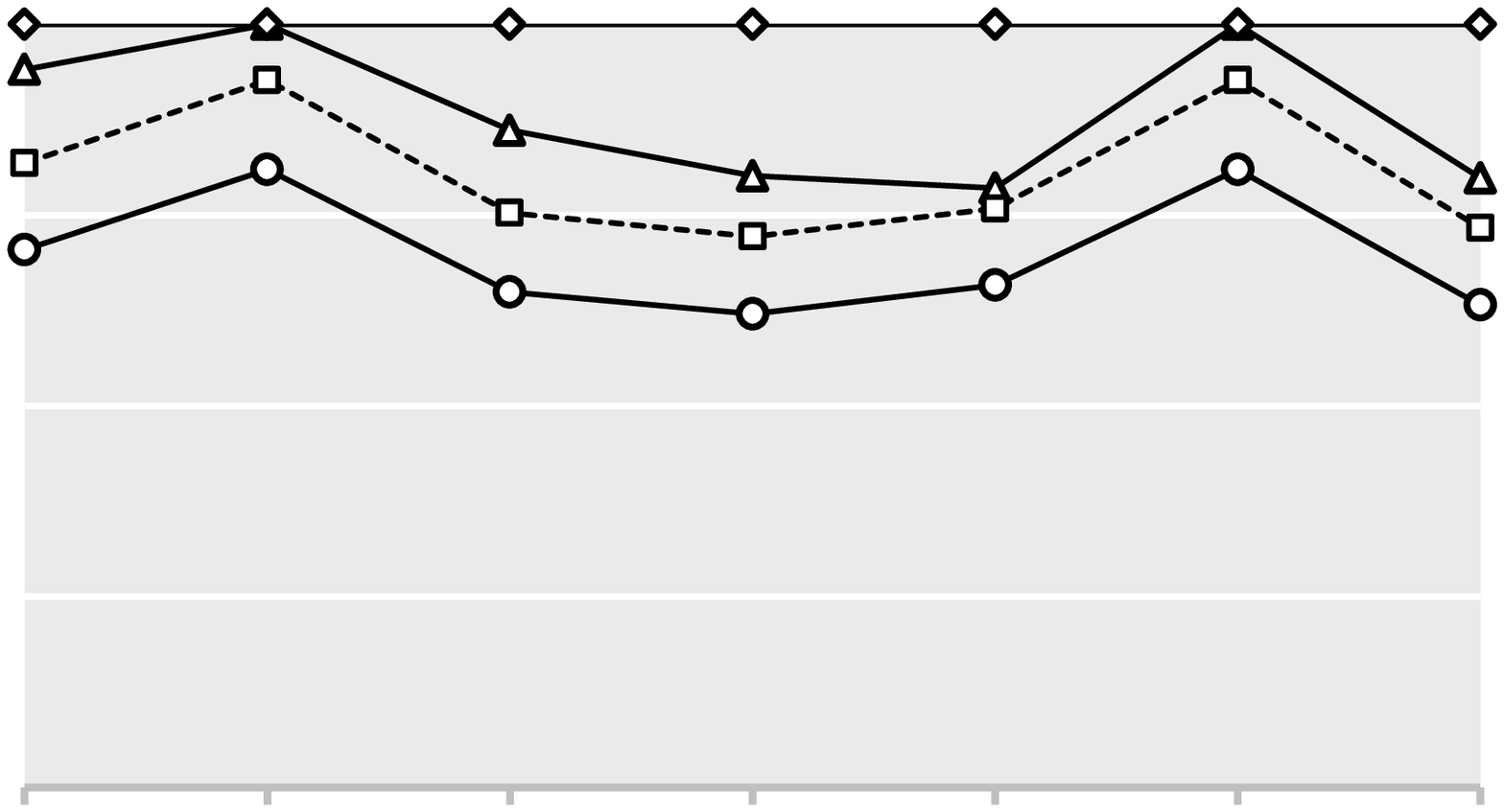}}
                \rput(0.225\linewidth,0.11){\footnotesize \textbf{\sf Network number}}
                \rput[bl](0.07,4.56){\footnotesize \textbf{{\sf Communication steps}}}
                \rput[l](0.05,0.91){\scriptsize $\mathsf{10^{0}}$}
                \rput[l](0.05,1.75){\scriptsize $\mathsf{10^{1}}$}
                \rput[l](0.05,2.58){\scriptsize $\mathsf{10^{2}}$}
                \rput[l](0.05,3.41){\scriptsize $\mathsf{10^{3}}$}
                \rput[l](0.05,4.25){\scriptsize $\mathsf{10^{4}}$}

                \rput[t](0.604,0.65){\scriptsize $\mathsf{1}$}
                \rput[t](1.630,0.65){\scriptsize $\mathsf{2}$}
                \rput[t](2.673,0.65){\scriptsize $\mathsf{3}$}
                \rput[t](3.721,0.65){\scriptsize $\mathsf{4}$}
                \rput[t](4.764,0.65){\scriptsize $\mathsf{5}$}
                \rput[t](5.800,0.65){\scriptsize $\mathsf{6}$}
                \rput[t](6.84,0.65){\scriptsize $\mathsf{7}$}

                \rput[l](7,2.97){\tiny \textbf{\sf D-ADMM}}
                \rput[l](7,3.32){\tiny \textbf{\sf D-Lasso}}
                \rput[l](7,3.59){\tiny \textbf{\sf MM/NGS}}
                \rput[l](7,4.32){\tiny \textbf{\sf MM/DQA, DN,}}
                \rput[l](7,4.06){\tiny \textbf{\sf Subgradient}}

            \end{pspicture}
	  }{
            \begin{pspicture}(0.48\linewidth,4.9cm)
                \rput(0.21\linewidth,2.55){\includegraphics[scale=0.38]{figures/RPCen5Acc5.eps}}
                \rput(0.21\linewidth,0.11){\footnotesize \textbf{\sf Network number}}
                \rput[bl](0.1,4.56){\footnotesize \textbf{{\sf Communication steps}}}
                \rput[l](0.08,0.91){\scriptsize $\mathsf{10^{0}}$}
                \rput[l](0.08,1.75){\scriptsize $\mathsf{10^{1}}$}
                \rput[l](0.08,2.58){\scriptsize $\mathsf{10^{2}}$}
                \rput[l](0.08,3.41){\scriptsize $\mathsf{10^{3}}$}
                \rput[l](0.08,4.25){\scriptsize $\mathsf{10^{4}}$}

                \rput[t](0.610,0.61){\scriptsize $\mathsf{1}$}
                \rput[t](1.665,0.61){\scriptsize $\mathsf{2}$}
                \rput[t](2.737,0.61){\scriptsize $\mathsf{3}$}
                \rput[t](3.818,0.61){\scriptsize $\mathsf{4}$}
                \rput[t](4.892,0.61){\scriptsize $\mathsf{5}$}
                \rput[t](5.953,0.61){\scriptsize $\mathsf{6}$}
                \rput[t](7.030,0.61){\scriptsize $\mathsf{7}$}

                \rput[l](7.18,2.97){\scriptsize \textbf{\sf D-ADMM}}
                \rput[l](7.18,3.32){\scriptsize \textbf{\sf D-Lasso}}
                \rput[l](7.18,3.59){\scriptsize \textbf{\sf MM/NGS}}
                \rput[l](7.18,4.32){\scriptsize \textbf{\sf MM/DQA, DN,}}
                \rput[l](7.18,4.06){\scriptsize \textbf{\sf Subgradient}}

            \end{pspicture}
         }
        }
        \caption{
        	Type I experiments: number of communication steps to reach accuracies of~$1\%$ and~$10^{-3}\%$ as a function of the network (see Table~\ref{Tab:Networks}).
        }
        \label{Fig:RPExperimentsScenarios}
    \end{figure*}

\mypar{Type I experiments} Figure~\ref{Fig:RPExperimentsScenarios}
shows the results of the type~I experiments. The left-hand
(resp.\ right-hand) side plots show, for each network, the number of
communication steps until each algorithm achieves a precision of~$1\%$
(resp.~$10^{-3}\%$) at a randomly selected node~$p$. This means we
count the number of communication steps until $\|x_p^{(k)} -
x^\star\|/\|x^\star\| \leq 10^{-2}$ or~$10^{-5}$. We allowed a
maximum number of~$10^{4}$ communication steps. 

    In Figure~\ref{Fig:RPExperimentsScenarios} we observe that the behavior of the algorithms in all  scenarios, except in scenario~$3$, is identical, so we will focus only on scenarios~$1$ and~$3$. Figure~\ref{Fig:RPExperimentsScenarios}(\subref*{SubFig:RPExperimentsScen1}) shows that, for scenario~$1$, D-ADMM requires the least number of communications to achieve both accuracies regardless the network. We can also see that for this scenario MM/DQA, DN, and Subgradient always reached the maximum number of~$10^4$ iterations before achieving any of the prescribed accuracies. As stated before, the behavior of the algorithms for the remaining scenarios (except scenario~$3$) is very similar. In scenario~$3$, Figure~\ref{Fig:RPExperimentsScenarios}(\subref*{SubFig:RPExperimentsScen3}), we see a different behavior: while D-ADMM required less communications than any of the $\rho$-dependent algorithms, the Subgradient required less communications to achieve the accuracy~$1\%$ for networks~$1$, $2$, and~$6$. However, if we let the algorithms continue executing, the Subgradient reaches the maximum number of communications before achieving the~$10^{-3}\%$ of accuracy, as can be seen in the right-hand plot of Figure~\ref{Fig:RPExperimentsScenarios}(\subref*{SubFig:RPExperimentsScen3}). Note that the relative behavior of the remaining algorithms is roughly the same for both accuracies.

    \isdraft{
      \begin{figure}
	\centering
        \subfigure{\label{SubFig:RPScenario1Network4}
            \begin{pspicture}(0.48\linewidth,5.5cm)
                \rput(0.23\linewidth,3.135){\includegraphics[scale=0.38]{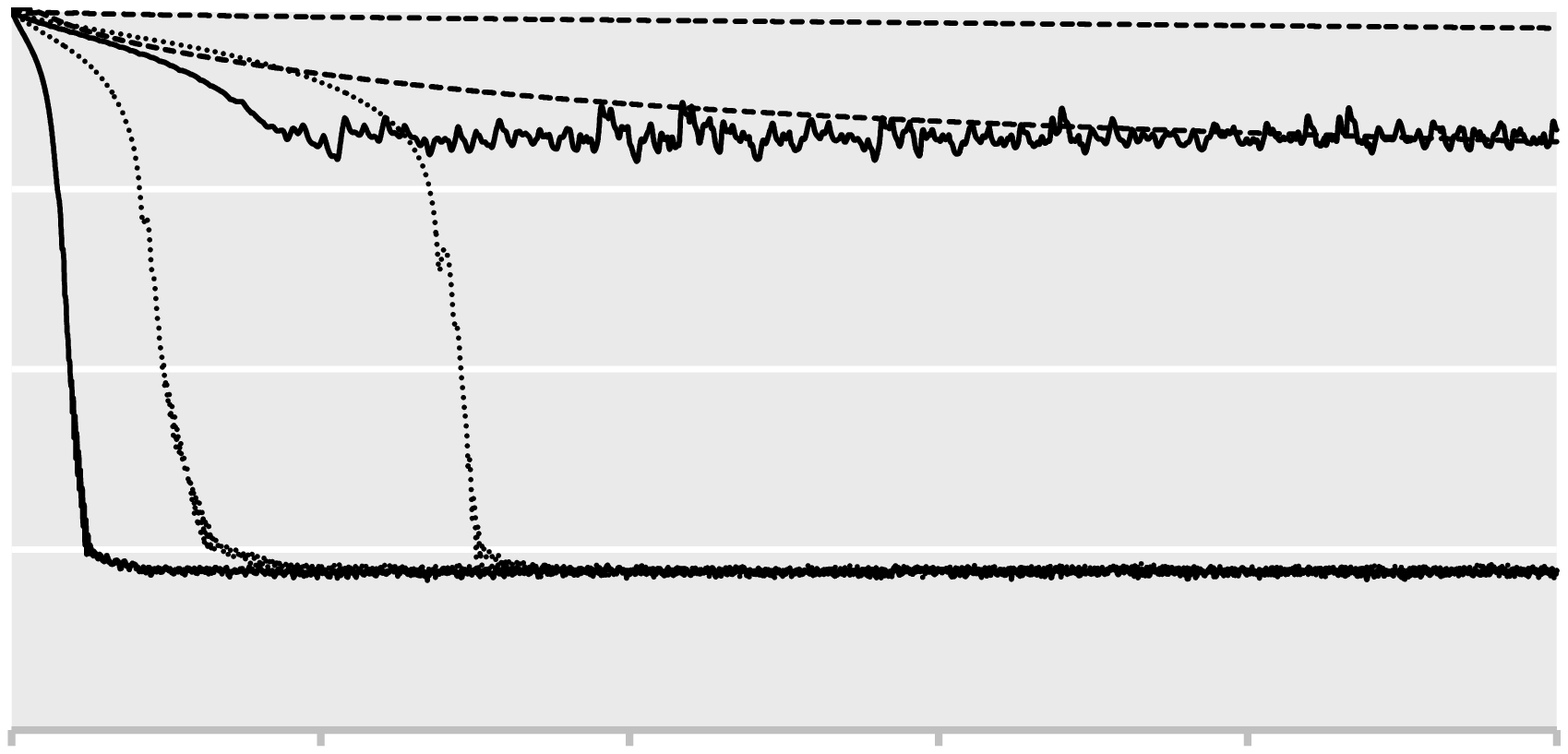}}
                \rput(0.23\linewidth,0.79){\footnotesize \textbf{\sf Communication steps}}
		\rput(0.23\linewidth,0.17){\footnotesize ($\textrm{a}$) Scenario~$1$, network~$4$}
                \rput[bl](-0.1,5.04){\footnotesize \textbf{{\sf Relative error:}} $\|x^{(k)}_p - x^\star\|/\|x^\star\|$}
                \rput[r](0.48,1.66){\scriptsize $\mathsf{10^{-8}}$}
                \rput[r](0.48,2.40){\scriptsize $\mathsf{10^{-6}}$}
                \rput[r](0.48,3.193){\scriptsize $\mathsf{10^{-4}}$}
                \rput[r](0.48,3.98){\scriptsize $\mathsf{10^{-2}}$}
                \rput[r](0.48,4.7){\scriptsize $\mathsf{10^{0\phantom{-}}}$}

		\rput[t](0.573,1.4){\scriptsize $\mathsf{0}$}
		\rput[t](1.87,1.4){\scriptsize $\mathsf{2000}$}
		\rput[t](3.15,1.4){\scriptsize $\mathsf{4000}$}
		\rput[t](4.44,1.4){\scriptsize $\mathsf{6000}$}
		\rput[t](5.74,1.4){\scriptsize $\mathsf{8000}$}
		\rput[t](7,1.4){\scriptsize $\mathsf{10000}$}

                \rput[l](5.9,3.7){\tiny \textbf{\sf MM/DQA}}
                \psline[linewidth=0.5pt](5.9,3.9)(5.5,4.5)
                \rput[l](1.5,1.8){\tiny \textbf{\sf D-ADMM}}
                \psline[linewidth=0.5pt](1.5,1.9)(1.26,2.2)
                \rput[l](1.3,3.4){\tiny \textbf{\sf D-Lasso}}
                \rput[l](2.65,2.8){\tiny \textbf{\sf MM/NGS}}
                \rput[l](3.0,3.55){\tiny \textbf{\sf Subgradient}}
                \psline[linewidth=0.5pt](3.0,3.7)(2.52,4.35)
                \rput[l](4.3,3.94){\tiny \textbf{\sf DN}}
               	
            \end{pspicture}
        }
        \hfill
        \subfigure{\label{SubFig:RPScenario3Network3}
            \begin{pspicture}(0.48\linewidth,5.5cm)
                \rput(0.23\linewidth,3.135){\includegraphics[scale=0.38]{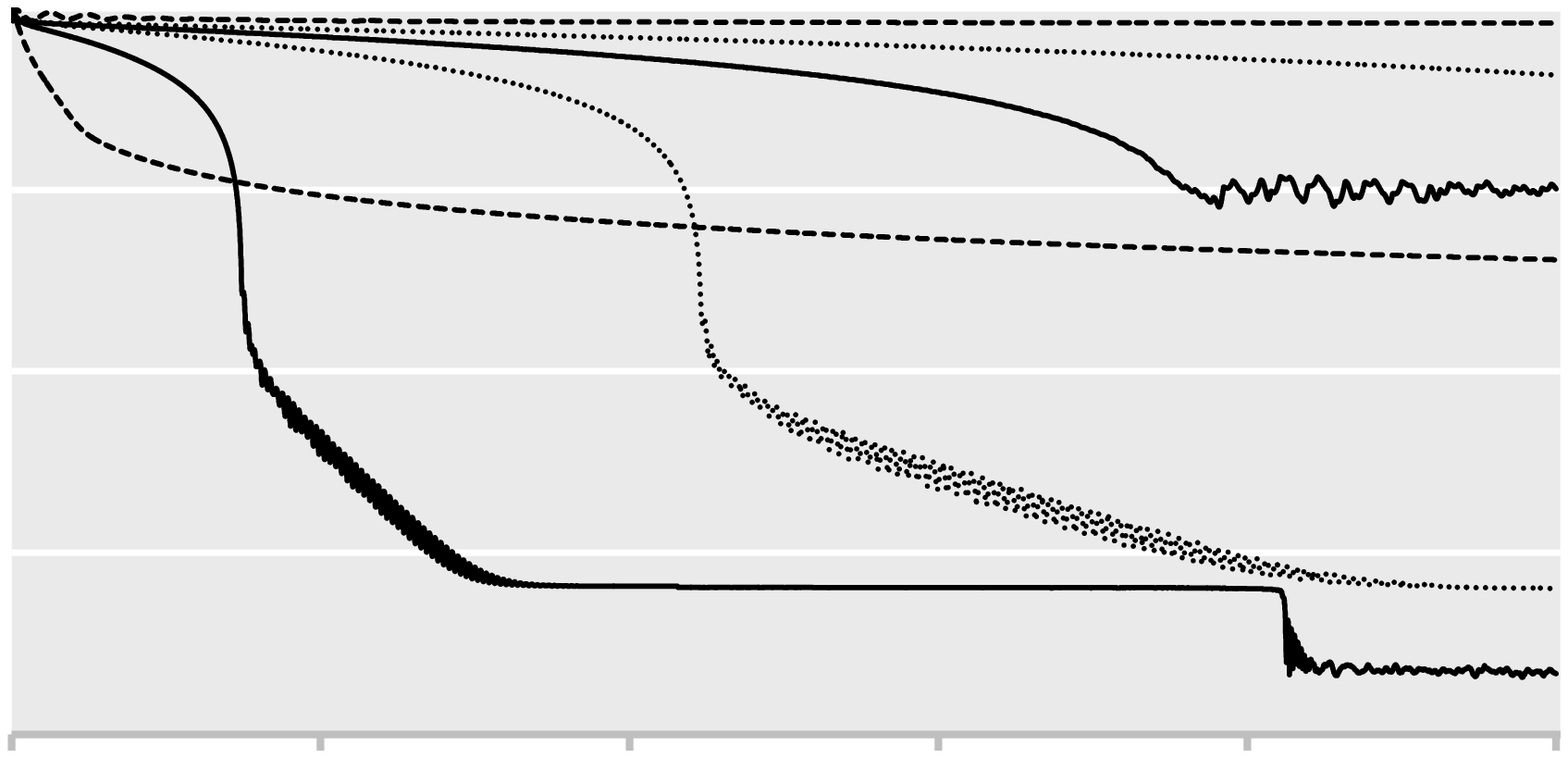}}
                \rput(0.23\linewidth,0.79){\footnotesize \textbf{\sf Communication steps}}
		\rput(0.23\linewidth,0.17){\footnotesize ($\textrm{b}$) Scenario~$3$, network~$3$}
                \rput[bl](-0.1,5.04){\footnotesize \textbf{{\sf Relative error:}} $\|x^{(k)}_p - x^\star\|/\|x^\star\|$}

                \rput[r](0.48,1.66){\scriptsize $\mathsf{10^{-8}}$}
                \rput[r](0.48,2.40){\scriptsize $\mathsf{10^{-6}}$}
                \rput[r](0.48,3.193){\scriptsize $\mathsf{10^{-4}}$}
                \rput[r](0.48,3.98){\scriptsize $\mathsf{10^{-2}}$}
                \rput[r](0.48,4.7){\scriptsize $\mathsf{10^{0\phantom{-}}}$}

		\rput[t](0.573,1.4){\scriptsize $\mathsf{0}$}
		\rput[t](1.87,1.4){\scriptsize $\mathsf{2000}$}
		\rput[t](3.15,1.4){\scriptsize $\mathsf{4000}$}
		\rput[t](4.44,1.4){\scriptsize $\mathsf{6000}$}
		\rput[t](5.74,1.4){\scriptsize $\mathsf{8000}$}
		\rput[t](7,1.4){\scriptsize $\mathsf{10000}$}
		
                \rput[l](7.1,1.84){\tiny \textbf{\sf D-ADMM}}
                \rput[l](7.1,2.24){\tiny \textbf{\sf D-Lasso}}
                \rput[l](7.1,3.63){\tiny \textbf{\sf Subgradient}}
                \rput[l](7.1,3.98){\tiny \textbf{\sf DN}}
                \rput[l](7.1,4.4){\tiny \textbf{\sf MM/NGS}}
                \rput[l](7.1,4.74){\tiny \textbf{\sf MM/DQA}}
            \end{pspicture}
        }
        \caption{Type I experiments: errors along the iterations (communication steps) of the algorithms for fixed scenarios and networks.}
        \label{Fig:RPErrorsAlongIterations}
        \isdraft{\vspace{-0.5cm}}{}
	\end{figure}
    }{
    \begin{figure*}
	\centering
        \subfigure{\label{SubFig:RPScenario1Network4}
            \begin{pspicture}(0.48\linewidth,5.5cm)
                \rput(0.219\linewidth,3.135){\includegraphics[scale=0.40]{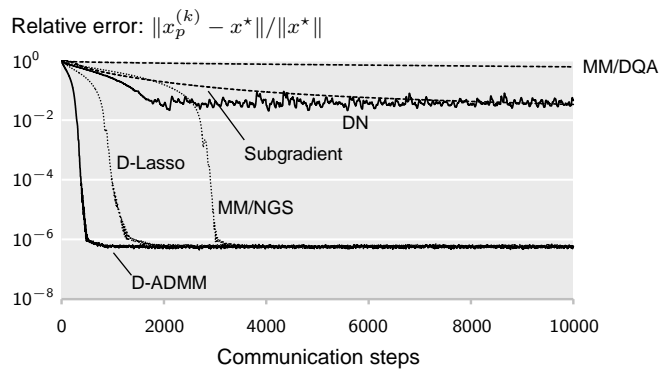}}
                \rput(0.22\linewidth,0.79){\footnotesize \textbf{\sf Communication steps}}
		\rput(0.22\linewidth,0.17){\footnotesize ($\textrm{a}$) Scenario~$1$, network~$4$}
                \rput[bl](-0.1,5.04){\footnotesize \textbf{{\sf Relative error:}} $\|x^{(k)}_p - x^\star\|/\|x^\star\|$}
                \rput[r](0.48,1.62){\scriptsize $\mathsf{10^{-8}}$}
                \rput[r](0.48,2.40){\scriptsize $\mathsf{10^{-6}}$}
                \rput[r](0.48,3.193){\scriptsize $\mathsf{10^{-4}}$}
                \rput[r](0.48,3.98){\scriptsize $\mathsf{10^{-2}}$}
                \rput[r](0.48,4.77){\scriptsize $\mathsf{10^{0\phantom{-}}}$}

		\rput[t](0.573,1.3){\scriptsize $\mathsf{0}$}
		\rput[t](1.939,1.3){\scriptsize $\mathsf{2000}$}
		\rput[t](3.293,1.3){\scriptsize $\mathsf{4000}$}
		\rput[t](4.648,1.3){\scriptsize $\mathsf{6000}$}
		\rput[t](6.028,1.3){\scriptsize $\mathsf{8000}$}
		\rput[t](7.384,1.3){\scriptsize $\mathsf{10000}$}

                \rput[l](7.50,4.68){\scriptsize \textbf{\sf MM/DQA}}
                \rput[l](1.5,1.8){\scriptsize \textbf{\sf D-ADMM}}
                \psline[linewidth=0.5pt](1.5,1.9)(1.26,2.2)
                \rput[l](1.3,3.4){\scriptsize \textbf{\sf D-Lasso}}
                \rput[l](2.65,2.8){\scriptsize \textbf{\sf MM/NGS}}
                \rput[l](3.0,3.5){\scriptsize \textbf{\sf Subgradient}}
                \psline[linewidth=0.5pt](3.0,3.7)(2.52,4.35)
                \rput[l](4.3,3.94){\scriptsize \textbf{\sf DN}}
               	
            \end{pspicture}
        }
        \hspace{0.12cm}
        \subfigure{\label{SubFig:RPScenario3Network3}
            \begin{pspicture}(0.48\linewidth,5.5cm)
                \rput(0.219\linewidth,3.135){\includegraphics[scale=0.40]{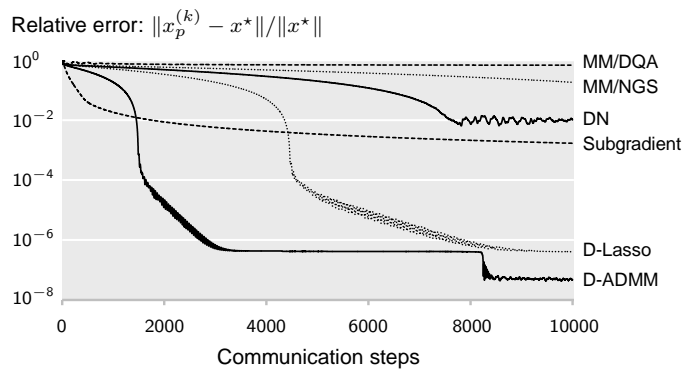}}
                \rput(0.22\linewidth,0.79){\footnotesize \textbf{\sf Communication steps}}
		\rput(0.22\linewidth,0.17){\footnotesize ($\textrm{b}$) Scenario~$3$, network~$3$}
                \rput[bl](-0.1,5.04){\footnotesize \textbf{{\sf Relative error:}} $\|x^{(k)}_p - x^\star\|/\|x^\star\|$}

                \rput[r](0.48,1.62){\scriptsize $\mathsf{10^{-8}}$}
                \rput[r](0.48,2.40){\scriptsize $\mathsf{10^{-6}}$}
                \rput[r](0.48,3.193){\scriptsize $\mathsf{10^{-4}}$}
                \rput[r](0.48,3.98){\scriptsize $\mathsf{10^{-2}}$}
                \rput[r](0.48,4.77){\scriptsize $\mathsf{10^{0\phantom{-}}}$}

		\rput[t](0.573,1.3){\scriptsize $\mathsf{0}$}
		\rput[t](1.939,1.3){\scriptsize $\mathsf{2000}$}
		\rput[t](3.293,1.3){\scriptsize $\mathsf{4000}$}
		\rput[t](4.648,1.3){\scriptsize $\mathsf{6000}$}
		\rput[t](6.028,1.3){\scriptsize $\mathsf{8000}$}
		\rput[t](7.384,1.3){\scriptsize $\mathsf{10000}$}
		
                \rput[l](7.50,1.84){\scriptsize \textbf{\sf D-ADMM}}
                \rput[l](7.50,2.24){\scriptsize \textbf{\sf D-Lasso}}
                \rput[l](7.50,3.63){\scriptsize \textbf{\sf Subgradient}}
                \rput[l](7.50,3.98){\scriptsize \textbf{\sf DN}}
                \rput[l](7.50,4.4){\scriptsize \textbf{\sf MM/NGS}}
                \rput[l](7.50,4.74){\scriptsize \textbf{\sf MM/DQA}}
            \end{pspicture}
        }
        \caption{Type I experiments: errors along the iterations (communication steps) of the algorithms for fixed scenarios and networks.}
        \label{Fig:RPErrorsAlongIterations}
	\end{figure*}
	}

    In Figure~\ref{Fig:RPErrorsAlongIterations} we show how the error of the estimate~$x_p$ at a random node~$p$ varies along the iterations, for each algorithm. Figure~\ref{Fig:RPErrorsAlongIterations}\subref{SubFig:RPScenario1Network4} shows the error for scenario~$1$ when the algorithms are executed in network number~$4$. Notice that the number of communications to achieve accuracies of~$1\%$ and~$10^{-3}\%$ agree with the plots of Figure~\ref{Fig:RPExperimentsScenarios}\subref{SubFig:RPExperimentsScen1}, for example D-ADMM takes less than~$10^{3}$ communication steps to achieve a~$10^{-5}$ precision. Figure~\ref{Fig:RPErrorsAlongIterations}\subref{SubFig:RPScenario3Network3} shows the errors for scenario~$3$ when we use network~$3$ (cf.\ with the plots of Figure~\ref{Fig:RPExperimentsScenarios}\subref{SubFig:RPExperimentsScen3}). Note the similarity of the curves of D-ADMM and D-Lasso: they have the same shape but the D-ADMM error is always smaller. This might happen because both methods use the same internal algorithm, albeit applied to different reformulations. Finally, note in Figure~\ref{Fig:RPErrorsAlongIterations}\subref{SubFig:RPScenario3Network3} how the error of the Subgradient evolves for scenario~$3$, network~$3$: the rate of convergence is very fast at the beginning, but after the first~$1000$ iterations it becomes very slow. This agrees with what was observed in Figure~\ref{Fig:RPExperimentsScenarios}\subref{SubFig:RPExperimentsScen3}.

  \begin{figure*}
    \centering
    \subfigure{\label{SubFig:TypeIIRPScen12}
      \isdraft{
	\begin{pspicture}(0.46\linewidth,5.3cm)
                \rput(0.225\linewidth,3.135){\includegraphics[scale=0.37]{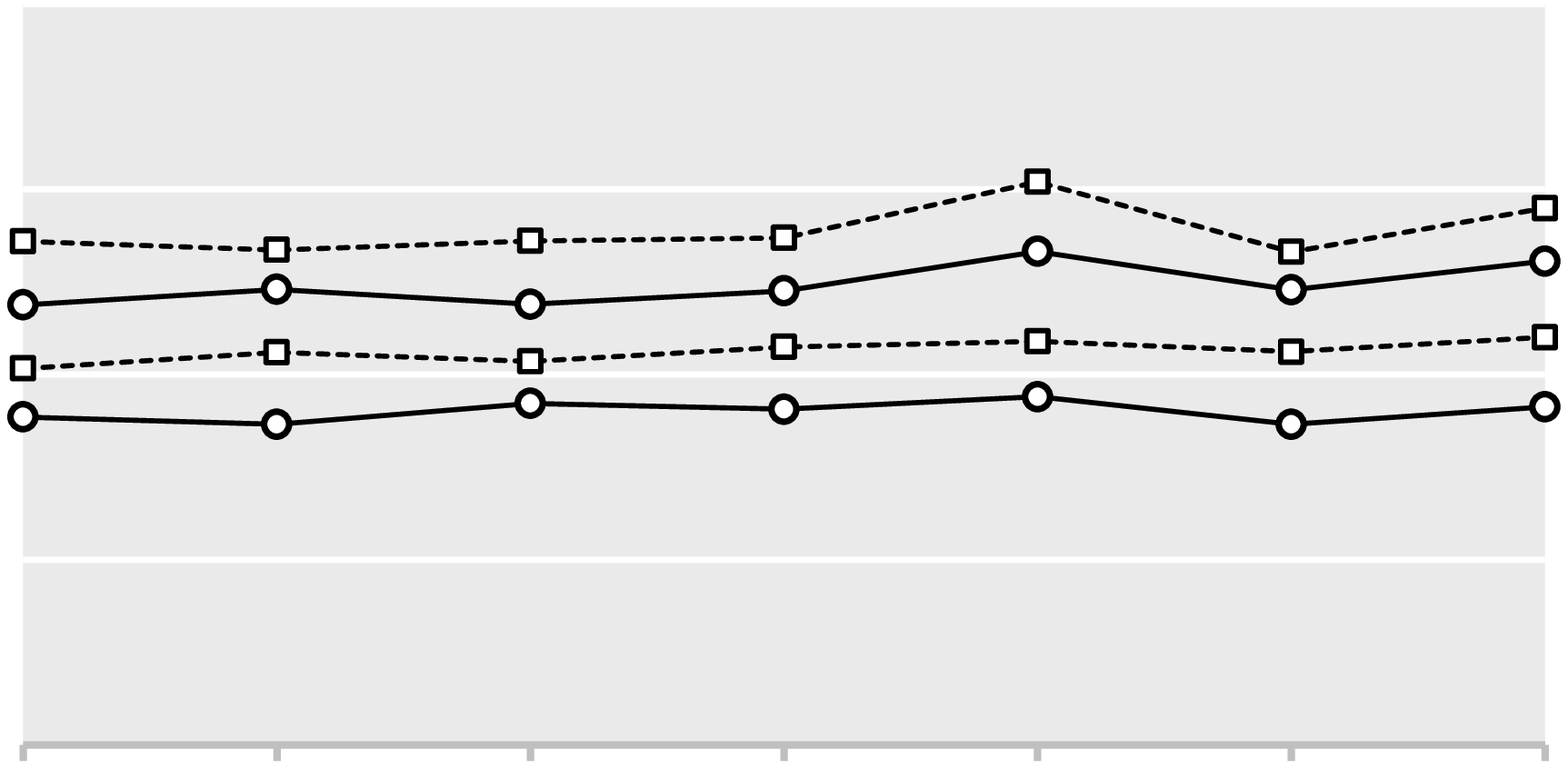}}
                \rput(0.225\linewidth,0.79){\footnotesize \textbf{\sf Network number}}
                \rput(0.225\linewidth,0.17){\footnotesize ($\textrm{a}$) Scenarios~$1$ and~$2$}

                \rput[bl](0.07,5.06){\footnotesize \textbf{{\sf Communication steps}}}
                \rput[l](0.05,1.62){\scriptsize $\mathsf{10^0}$}
                \rput[l](0.05,2.41){\scriptsize $\mathsf{10^1}$}
                \rput[l](0.05,3.21){\scriptsize $\mathsf{10^2}$}
                \rput[l](0.05,4.00){\scriptsize $\mathsf{10^3}$}
                \rput[l](0.05,4.75){\scriptsize $\mathsf{10^4}$}

                \rput[t](0.604,1.42){\scriptsize $\mathsf{1}$}
                \rput[t](1.630,1.42){\scriptsize $\mathsf{2}$}
                \rput[t](2.673,1.42){\scriptsize $\mathsf{3}$}
                \rput[t](3.721,1.42){\scriptsize $\mathsf{4}$}
                \rput[t](4.764,1.42){\scriptsize $\mathsf{5}$}
                \rput[t](5.800,1.42){\scriptsize $\mathsf{6}$}
                \rput[t](6.84,1.42){\scriptsize $\mathsf{7}$}

                \rput[l](7,3.65){\tiny \textbf{\sf D-ADMM: 1}}
                \rput[l](7,3.91){\tiny \textbf{\sf D-Lasso: 1}}
                \rput[l](7,3.05){\tiny \textbf{\sf D-ADMM: 2}}
                \rput[l](7,3.35){\tiny \textbf{\sf D-Lasso: 2}}
            \end{pspicture}
      }{
	\begin{pspicture}(0.46\linewidth,5.5cm)
                \rput(0.214\linewidth,3.135){\includegraphics[scale=0.38]{figures/RPBestRhoCen12.eps}}
                \rput(0.21\linewidth,0.79){\footnotesize \textbf{\sf Network number}}
                \rput(0.21\linewidth,0.17){\footnotesize ($\textrm{a}$) Scenarios~$1$ and~$2$}

                \rput[bl](0.08,5.06){\footnotesize \textbf{{\sf Communication steps}}}
                \rput[l](0.08,1.62){\scriptsize $\mathsf{10^0}$}
                \rput[l](0.08,2.41){\scriptsize $\mathsf{10^1}$}
                \rput[l](0.08,3.21){\scriptsize $\mathsf{10^2}$}
                \rput[l](0.08,4.00){\scriptsize $\mathsf{10^3}$}
                \rput[l](0.08,4.75){\scriptsize $\mathsf{10^4}$}

                \rput[t](0.652,1.34){\scriptsize $\mathsf{1}$}
                \rput[t](1.732,1.34){\scriptsize $\mathsf{2}$}
                \rput[t](2.803,1.34){\scriptsize $\mathsf{3}$}
                \rput[t](3.876,1.34){\scriptsize $\mathsf{4}$}
                \rput[t](4.960,1.34){\scriptsize $\mathsf{5}$}
                \rput[t](6.044,1.34){\scriptsize $\mathsf{6}$}
                \rput[t](7.111,1.34){\scriptsize $\mathsf{7}$}

                \rput[l](7.29,3.65){\scriptsize \textbf{\sf D-ADMM: 1}}
                \rput[l](7.29,3.91){\scriptsize \textbf{\sf D-Lasso: 1}}
                \rput[l](7.29,3.05){\scriptsize \textbf{\sf D-ADMM: 2}}
                \rput[l](7.29,3.35){\scriptsize \textbf{\sf D-Lasso: 2}}
            \end{pspicture}
        }
        }
        \isdraft{\hfill}{\hspace{0.12cm}}
        \subfigure{\label{SubFig:TypeIIRPScen345}
	  \isdraft{
	    \begin{pspicture}(0.46\linewidth,5.3cm)
                \rput(0.225\linewidth,3.135){\includegraphics[scale=0.37]{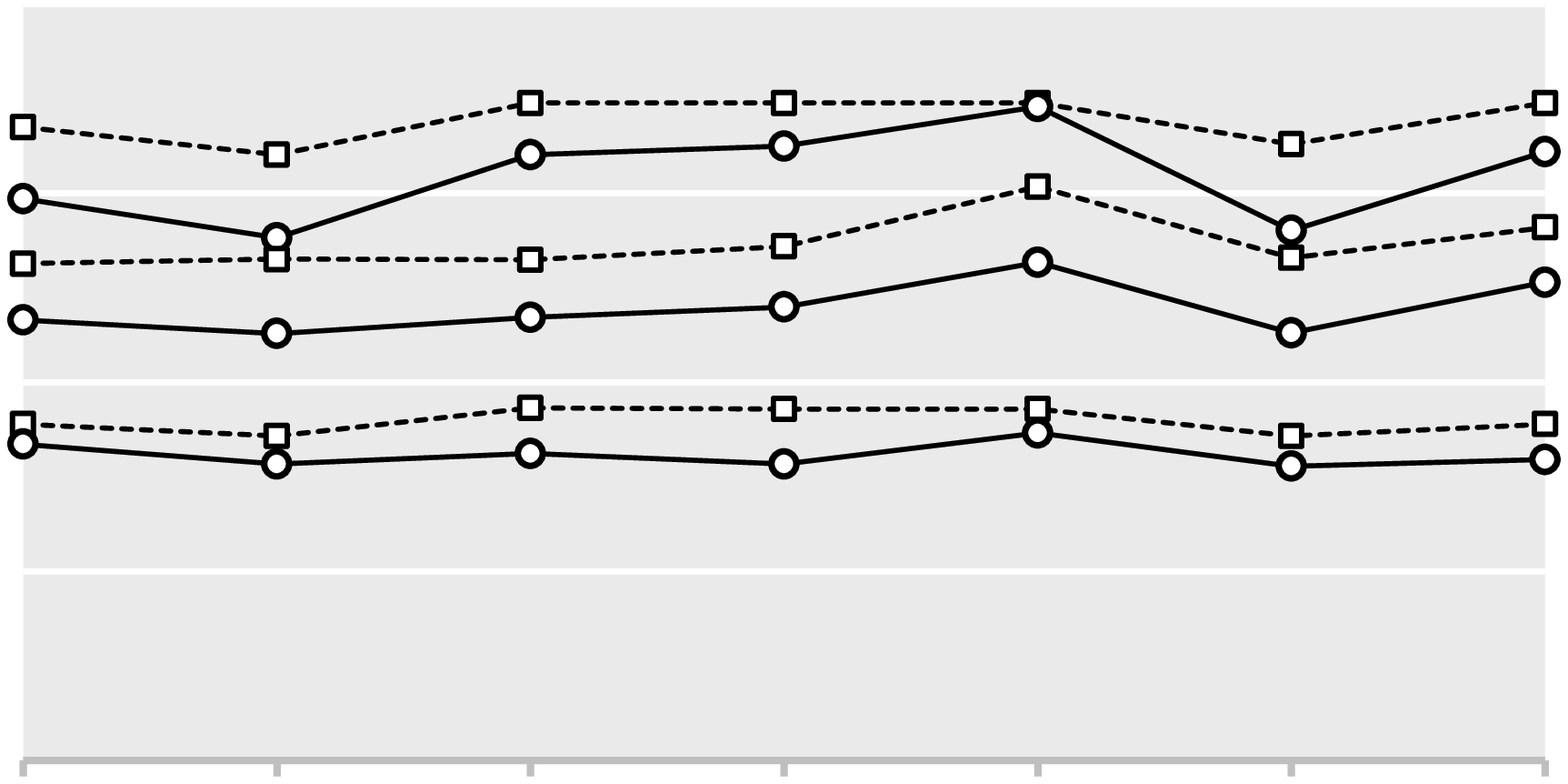}}
                \rput(0.225\linewidth,0.79){\footnotesize \textbf{\sf Network number}}
		\rput(0.225\linewidth,0.17){\footnotesize ($\textrm{b}$) Scenarios~$3$, $4$, and~$5$}

		\rput[bl](0.07,5.06){\footnotesize \textbf{{\sf Communication steps}}}
                \rput[l](0.05,1.62){\scriptsize $\mathsf{10^0}$}
                \rput[l](0.05,2.41){\scriptsize $\mathsf{10^1}$}
                \rput[l](0.05,3.21){\scriptsize $\mathsf{10^2}$}
                \rput[l](0.05,4.00){\scriptsize $\mathsf{10^3}$}
                \rput[l](0.05,4.75){\scriptsize $\mathsf{10^4}$}

                \rput[t](0.604,1.42){\scriptsize $\mathsf{1}$}
                \rput[t](1.630,1.42){\scriptsize $\mathsf{2}$}
                \rput[t](2.673,1.42){\scriptsize $\mathsf{3}$}
                \rput[t](3.721,1.42){\scriptsize $\mathsf{4}$}
                \rput[t](4.764,1.42){\scriptsize $\mathsf{5}$}
                \rput[t](5.800,1.42){\scriptsize $\mathsf{6}$}
                \rput[t](6.84,1.42){\scriptsize $\mathsf{7}$}

                \rput[l](7,4.14){\tiny \textbf{\sf D-ADMM: 3}}
                \rput[l](7,4.40){\tiny \textbf{\sf D-Lasso: 3}}
		\rput[l](7,3.58){\tiny \textbf{\sf D-ADMM: 5}}
		\rput[l](7,3.83){\tiny \textbf{\sf D-Lasso: 5}}
		\rput[l](7,2.81){\tiny \textbf{\sf D-ADMM: 4}}
               	\rput[l](7,3.04){\tiny \textbf{\sf D-Lasso: 4}}
               	
            \end{pspicture}
	  }{
            \begin{pspicture}(0.46\linewidth,5.5cm)
                \rput(0.214\linewidth,3.135){\includegraphics[scale=0.38]{figures/RPBestRhoCen345.eps}}
                \rput(0.21\linewidth,0.79){\footnotesize \textbf{\sf Network number}}
		\rput(0.21\linewidth,0.17){\footnotesize ($\textrm{b}$) Scenarios~$3$, $4$, and~$5$}

		\rput[bl](0.08,5.06){\footnotesize \textbf{{\sf Communication steps}}}
                \rput[l](0.08,1.62){\scriptsize $\mathsf{10^0}$}
                \rput[l](0.08,2.41){\scriptsize $\mathsf{10^1}$}
                \rput[l](0.08,3.21){\scriptsize $\mathsf{10^2}$}
                \rput[l](0.08,4.00){\scriptsize $\mathsf{10^3}$}
                \rput[l](0.08,4.75){\scriptsize $\mathsf{10^4}$}

                \rput[t](0.652,1.34){\scriptsize $\mathsf{1}$}
                \rput[t](1.732,1.34){\scriptsize $\mathsf{2}$}
                \rput[t](2.803,1.34){\scriptsize $\mathsf{3}$}
                \rput[t](3.876,1.34){\scriptsize $\mathsf{4}$}
                \rput[t](4.960,1.34){\scriptsize $\mathsf{5}$}
                \rput[t](6.044,1.34){\scriptsize $\mathsf{6}$}
                \rput[t](7.111,1.34){\scriptsize $\mathsf{7}$}

                \rput[l](7.29,4.14){\scriptsize \textbf{\sf D-ADMM: 3}}
                \rput[l](7.29,4.40){\scriptsize \textbf{\sf D-Lasso: 3}}
		\rput[l](7.29,3.58){\scriptsize \textbf{\sf D-ADMM: 5}}
		\rput[l](7.29,3.83){\scriptsize \textbf{\sf D-Lasso: 5}}
		\rput[l](7.29,2.81){\scriptsize \textbf{\sf D-ADMM: 4}}
               	\rput[l](7.29,3.04){\scriptsize \textbf{\sf D-Lasso: 4}}
               	
            \end{pspicture}
        }
        }
       \caption{Type II experiments: number of communication steps to reach $10^{-3}\%$ of accuracy or~$3000$ communication steps.}			 \label{Fig:Type2ExperimentsRP}
	\end{figure*}

\mypar{Type II experiments} For the type II experiments we only
considered the two best algorithms: D-ADMM and D-Lasso.
Figure~\ref{Fig:Type2ExperimentsRP} shows for each network the number
of communication steps to reach an accuracy of~$10^{-3}\%$. We
allowed for maximally~$3000$ communication steps (these were only
achieved by D-Lasso in scenario~$3$ for networks~$3$, $4$, and~$5$, as
can be seen in
Figure~\ref{Fig:Type2ExperimentsRP}\subref{SubFig:TypeIIRPScen345}).
We observed that the best values of~$\rho$ for D-ADMM were
always~$10^{-2}$, $10^{-1}$, or~$1$. For example, D-ADMM had the best
performance for~$\rho=1$ for scenarios~$1$, $3$, and~$5$ when the
networks were either~$5$ or~$7$. For instance, for scenario~$1$,
network~$5$ D-ADMM took~$462$ communication steps (see Figure~\ref{Fig:Type2ExperimentsRP}\subref{SubFig:TypeIIRPScen12}),
the same number observed in the type I experiments, in right-hand plot
of
Figure~\ref{Fig:RPExperimentsScenarios}\subref{SubFig:RPExperimentsScen1}. Recall
that~$\rho$ was fixed at~$1$ for D-ADMM in the type~I experiments. This
also means that in the type~II experiments the number of
communications for D-ADMM decreased except for scenarios~$1$, $3$,
and~$5$ when the networks were either~$5$ or~$7$. The same phenomenon
happened for D-Lasso: the optimal~$\rho$ was~$1$ only in scenarios~$1$
and~$5$ for the $5$th network; and the optimal $\rho$'s were~$10^{-2}$, $10^{-1}$, or~$1$.

We conclude from Figure~\ref{Fig:Type2ExperimentsRP} that
D-ADMM requires less communication steps than D-Lasso, independently of the
scenario or network type. 
Excluding the cases D-Lasso reached the maximum number of iterations, we see that in average D-ADMM uses $51\%$ of D-Lasso's number of communications ($11\%$ of standard deviation). The largest difference occurred in scenario~$3$, network~$6$, where D-ADMM used~$35\%$ of the communications D-Lasso used; this number was~$78\%$ for scenario~$4$, network~$1$, the smallest difference that occurred.

\begin{figure}[h]
    \centering
    \isdraft{
			\begin{pspicture}(0.5\linewidth,5.0cm)
				\rput(0.225\linewidth,2.635){\includegraphics[scale=0.38]{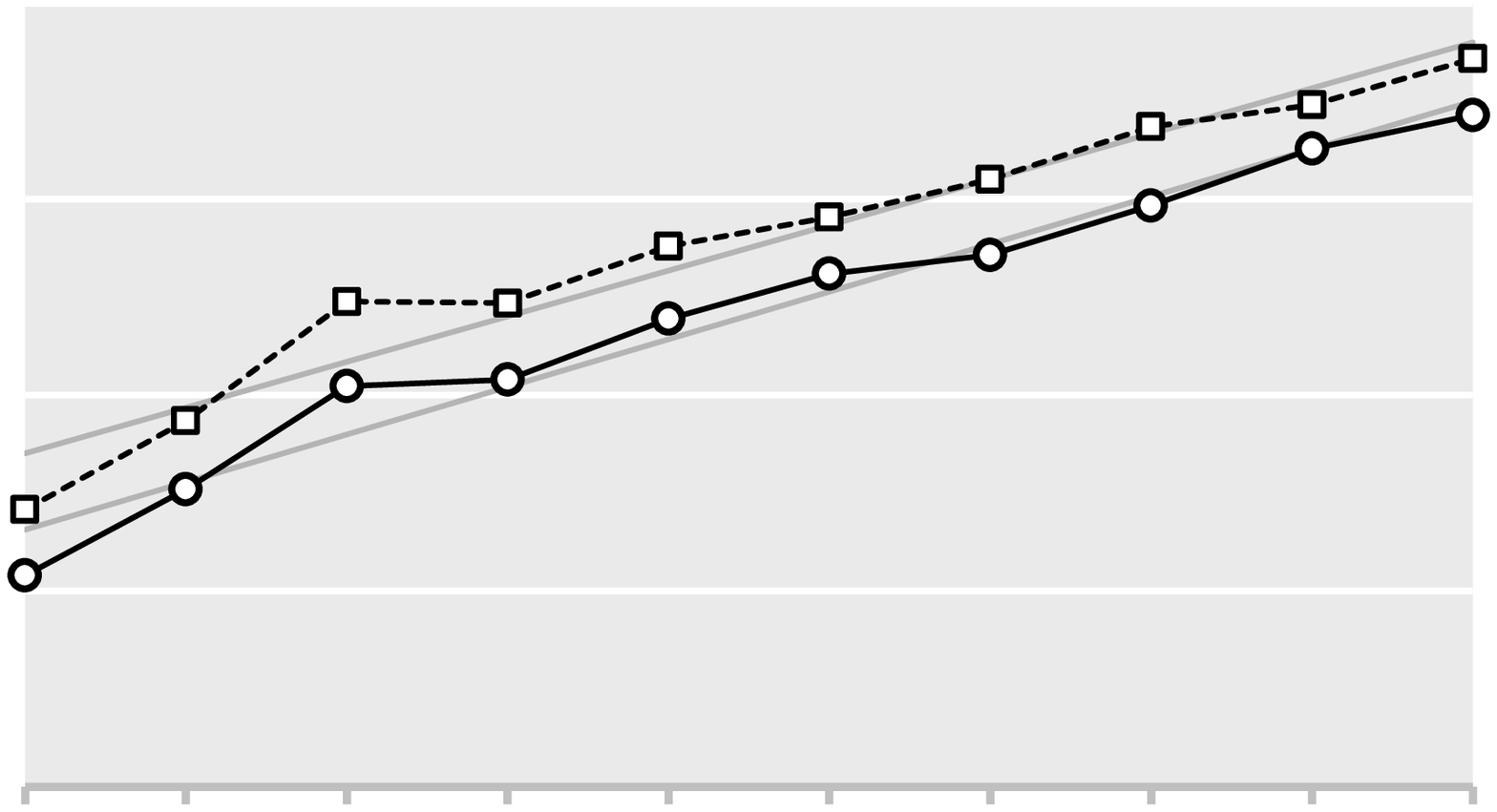}}
        \rput(0.225\linewidth,0.10){\footnotesize \textbf{\sf Number of nodes}}      

        \rput[bl](0.04,4.71){\footnotesize \textbf{{\sf Communication steps}}}
        \rput[l](0.02,0.960){\scriptsize $\mathsf{10^0}$}
        \rput[l](0.02,1.845){\scriptsize $\mathsf{10^1}$}
        \rput[l](0.02,2.705){\scriptsize $\mathsf{10^2}$}
        \rput[l](0.02,3.560){\scriptsize $\mathsf{10^3}$}
        \rput[l](0.02,4.400){\scriptsize $\mathsf{10^4}$}

        \rput[t](0.618,0.75){\scriptsize $\mathsf{2^1}$}
        \rput[t](1.329,0.75){\scriptsize $\mathsf{2^2}$}
        \rput[t](2.041,0.75){\scriptsize $\mathsf{2^3}$}
        \rput[t](2.758,0.75){\scriptsize $\mathsf{2^4}$}
        \rput[t](3.468,0.75){\scriptsize $\mathsf{2^5}$}
        \rput[t](4.181,0.75){\scriptsize $\mathsf{2^6}$}
        \rput[t](4.888,0.75){\scriptsize $\mathsf{2^7}$}
        \rput[t](5.605,0.75){\scriptsize $\mathsf{2^8}$}
        \rput[t](6.309,0.75){\scriptsize $\mathsf{2^9}$}
        \rput[t](7.037,0.75){\scriptsize $\mathsf{2^{10}}$}

				\rput[l](7.14,3.9){\scriptsize \textbf{\sf D-ADMM}}
        \rput[l](7.14,4.2){\scriptsize \textbf{\sf D-Lasso}}
        
        \rput(3.0,2.3){\psscalebox{1}{\rput(0,0){\scriptsize $11.7 \cdot P^{0.8}$}}}
        \rput(1.3,3.25){\psscalebox{1}{\rput(0,0){\scriptsize $29.5 \cdot P^{0.77}$}}}
        
        \psline[linewidth=0.5pt]{-}(1.2,3.05)(1.0,2.6)
				\psline[linewidth=0.5pt]{-}(2.5,2.4)(2.3,2.55)
     \end{pspicture}
		}{
			\begin{pspicture}(0.9\linewidth,5.0cm)
				\rput(0.425\linewidth,2.635){\includegraphics[scale=0.38]{figures/LargeNetworks.eps}}
        \rput(0.425\linewidth,0.10){\footnotesize \textbf{\sf Number of nodes}}

        \rput[bl](0.04,4.71){\footnotesize \textbf{{\sf Communication steps}}}
        \rput[l](0.02,0.960){\scriptsize $\mathsf{10^0}$}
        \rput[l](0.02,1.845){\scriptsize $\mathsf{10^1}$}
        \rput[l](0.02,2.705){\scriptsize $\mathsf{10^2}$}
        \rput[l](0.02,3.560){\scriptsize $\mathsf{10^3}$}
        \rput[l](0.02,4.400){\scriptsize $\mathsf{10^4}$}

        \rput[t](0.618,0.75){\scriptsize $\mathsf{2^1}$}
        \rput[t](1.329,0.75){\scriptsize $\mathsf{2^2}$}
        \rput[t](2.041,0.75){\scriptsize $\mathsf{2^3}$}
        \rput[t](2.758,0.75){\scriptsize $\mathsf{2^4}$}
        \rput[t](3.468,0.75){\scriptsize $\mathsf{2^5}$}
        \rput[t](4.181,0.75){\scriptsize $\mathsf{2^6}$}
        \rput[t](4.888,0.75){\scriptsize $\mathsf{2^7}$}
        \rput[t](5.605,0.75){\scriptsize $\mathsf{2^8}$}
        \rput[t](6.309,0.75){\scriptsize $\mathsf{2^9}$}
        \rput[t](7.037,0.75){\scriptsize $\mathsf{2^{10}}$}

				\rput[l](7.14,3.9){\scriptsize \textbf{\sf D-ADMM}}
        \rput[l](7.14,4.2){\scriptsize \textbf{\sf D-Lasso}}
        
        \rput(3.0,2.3){\psscalebox{1}{\rput(0,0){\scriptsize $11.7 \cdot P^{0.8}$}}}
        \rput(1.3,3.25){\psscalebox{1}{\rput(0,0){\scriptsize $29.5 \cdot P^{0.77}$}}}
        
        \psline[linewidth=0.5pt]{-}(1.2,3.05)(1.0,2.6)
				\psline[linewidth=0.5pt]{-}(2.5,2.4)(2.3,2.55)
     \end{pspicture}
		}
     \caption{Type II experiments for row partition: number of communication steps to reach $0.1\%$ of accuracy as a function of the network size. The straight lines represent a linear fit.}
     \label{Fig:RPTypeIIIncreaseNetworkSize}
  \end{figure}

	Figure~\ref{Fig:RPTypeIIIncreaseNetworkSize} shows another type~II experiment: we fixed the scenario and network type: Scenario~$3$, Watts-Strogatz with parameters~$(4,0.6)$; and observed how the number of communication steps varies as the size of the network increases. The number of nodes varied from~$2$ (each node stores~$512$ rows) to~$1024$ (each node stores~$1$ row) and was always a power of~$2$. D-ADMM and D-Lasso stopped after reaching~$0.1\%$ of accuracy. As shown by the gray straight lines in Figure~\ref{Fig:RPTypeIIIncreaseNetworkSize}, the communication steps in both algorithms increases approximately linearly in a log-log plot. The model we used to compute those lines was~$\log_{10}C = \alpha \log_2P + \beta$, where~$C$ is the number of communication steps, $P$ the number of nodes, and~$\alpha$ and~$\beta$ are the parameters to be found for each line. The minimum least squares error yielded~$(\alpha,\beta) = (0.243,1.07)$ for D-ADMM and~$(\alpha,\beta) = (0.233,1.47)$ for D-Lasso. Therefore, $C \simeq 11.7 \cdot P^{0.8}$ for D-ADMM and $C \simeq 29.5 \cdot P^{0.77}$ for D-Lasso, showing a less-than-linear increase of the communication steps with the number of nodes, for both algorithms. Also, the difference between the lines' offsets reveals that D-Lasso took in average~$2.5$ times more communications than D-ADMM. The average number of colors was~$4.6$, which means that in a collision-free network D-ADMM would be~$1.8$ times slower than D-Lasso. Again, the optimal $\rho$'s were~$10^{-2}$, $10^{-1}$, or~$1$, but we noticed a curious pattern on both algorithms: the optimal value for~$\rho$ decreased as the size of the network increased.
    
  \mypar{Results for the column partition} For the column partition we only executed type~II experiments. While the scenarios were the same as before (Table~\ref{Tab:Scenarios}), we changed the networks: they now have~$10$ nodes (for scenarios~$1$, $2$, and~$4$) or~$8$ nodes (for scenarios~$3$ and~$5$). All nodes thus store the same number of columns, i.e., the number of columns~$n$ is divisible by the number of nodes~$P$. The model for generating these networks is the same as in Table~\ref{Tab:Networks}. 
    In all experiments we set the regularization parameter to~$\delta = 10^{-3}$, a value that always allowed the recovery of the solution of BP, as we will see. 
    
\begin{figure*}
	\centering
	\subfigure{\label{SubFig:CPTypeIIRPScen12}
	\isdraft{
	    \begin{pspicture}(0.46\linewidth,5.3cm)
                \rput(0.225\linewidth,3.135){\includegraphics[scale=0.37]{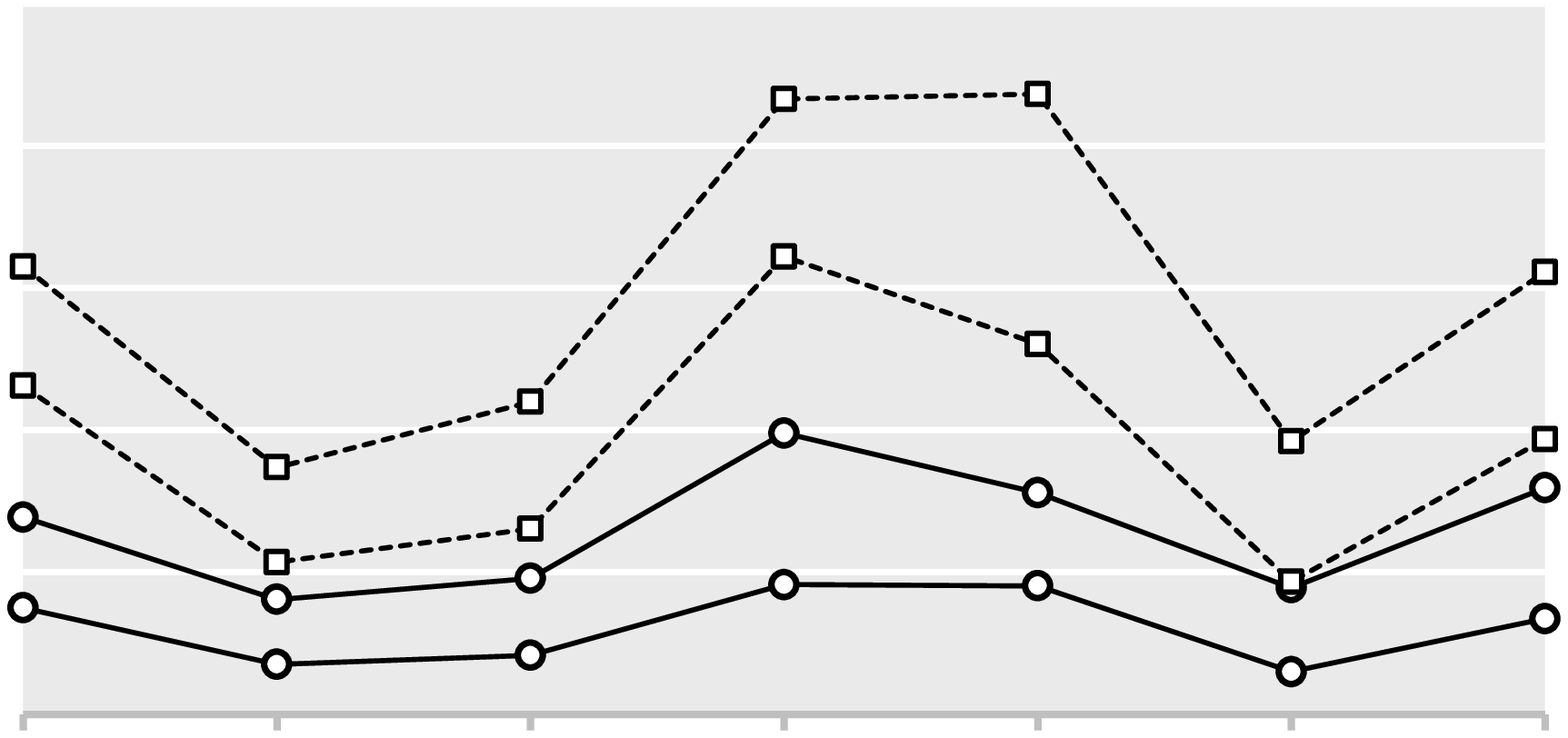}}
                \rput(0.225\linewidth,0.79){\footnotesize \textbf{\sf Network number}}
                \rput(0.225\linewidth,0.17){\footnotesize ($\textrm{a}$) Scenarios~$1$ and~$2$}

                \rput[bl](-0.08,5.02){\footnotesize \textbf{{\sf Communication steps}}}
                \rput[r](0.48,1.67){\scriptsize $0$}
                \rput[r](0.48,2.279){\scriptsize $400$}
                \rput[r](0.48,2.884){\scriptsize $800$}
                \rput[r](0.48,3.505){\scriptsize $1200$}
                \rput[r](0.48,4.09){\scriptsize $1600$}
                \rput[r](0.48,4.66){\scriptsize $2000$}

                \rput[t](0.604,1.42){\scriptsize $\mathsf{1}$}
                \rput[t](1.630,1.42){\scriptsize $\mathsf{2}$}
                \rput[t](2.673,1.42){\scriptsize $\mathsf{3}$}
                \rput[t](3.721,1.42){\scriptsize $\mathsf{4}$}
                \rput[t](4.764,1.42){\scriptsize $\mathsf{5}$}
                \rput[t](5.800,1.42){\scriptsize $\mathsf{6}$}
                \rput[t](6.84,1.42){\scriptsize $\mathsf{7}$}

                \rput[l](7,2.613){\tiny \textbf{\sf D-ADMM: 1}}
                \rput[l](7,3.55){\tiny \textbf{\sf D-Lasso: 1}}
                \rput[l](7,2.07){\tiny \textbf{\sf D-ADMM: 2}}
                \rput[l](7,2.84){\tiny \textbf{\sf D-Lasso: 2}}
            \end{pspicture}
	}{
            \begin{pspicture}(0.46\linewidth,5.4cm)
                \rput(0.214\linewidth,3.135){\includegraphics[scale=0.38]{figures/CPBestRhoCen12.eps}}
                \rput(0.21\linewidth,0.79){\footnotesize \textbf{\sf Network number}}
                \rput(0.21\linewidth,0.17){\footnotesize ($\textrm{a}$) Scenarios~$1$ and~$2$}

                \rput[bl](-0.08,5.02){\footnotesize \textbf{{\sf Communication steps}}}
                \rput[r](0.48,1.67){\scriptsize $0$}
                \rput[r](0.48,2.279){\scriptsize $400$}
                \rput[r](0.48,2.884){\scriptsize $800$}
                \rput[r](0.48,3.505){\scriptsize $1200$}
                \rput[r](0.48,4.09){\scriptsize $1600$}
                \rput[r](0.48,4.66){\scriptsize $2000$}

                \rput[t](0.652,1.34){\scriptsize $\mathsf{1}$}
                \rput[t](1.732,1.34){\scriptsize $\mathsf{2}$}
                \rput[t](2.803,1.34){\scriptsize $\mathsf{3}$}
                \rput[t](3.876,1.34){\scriptsize $\mathsf{4}$}
                \rput[t](4.960,1.34){\scriptsize $\mathsf{5}$}
                \rput[t](6.044,1.34){\scriptsize $\mathsf{6}$}
                \rput[t](7.111,1.34){\scriptsize $\mathsf{7}$}

                \rput[l](7.29,2.613){\scriptsize \textbf{\sf D-ADMM: 1}}
                \rput[l](7.29,3.55){\scriptsize \textbf{\sf D-Lasso: 1}}
                \rput[l](7.29,2.07){\scriptsize \textbf{\sf D-ADMM: 2}}
                \rput[l](7.29,2.84){\scriptsize \textbf{\sf D-Lasso: 2}}
            \end{pspicture}
        }
        }
        \isdraft{\hfill}{\hspace{0.12cm}}
        \subfigure{\label{SubFig:CPTypeIIRPScen345}
	  \isdraft{
	    \begin{pspicture}(0.46\linewidth,5.3cm)
                \rput(0.225\linewidth,3.135){\includegraphics[scale=0.37]{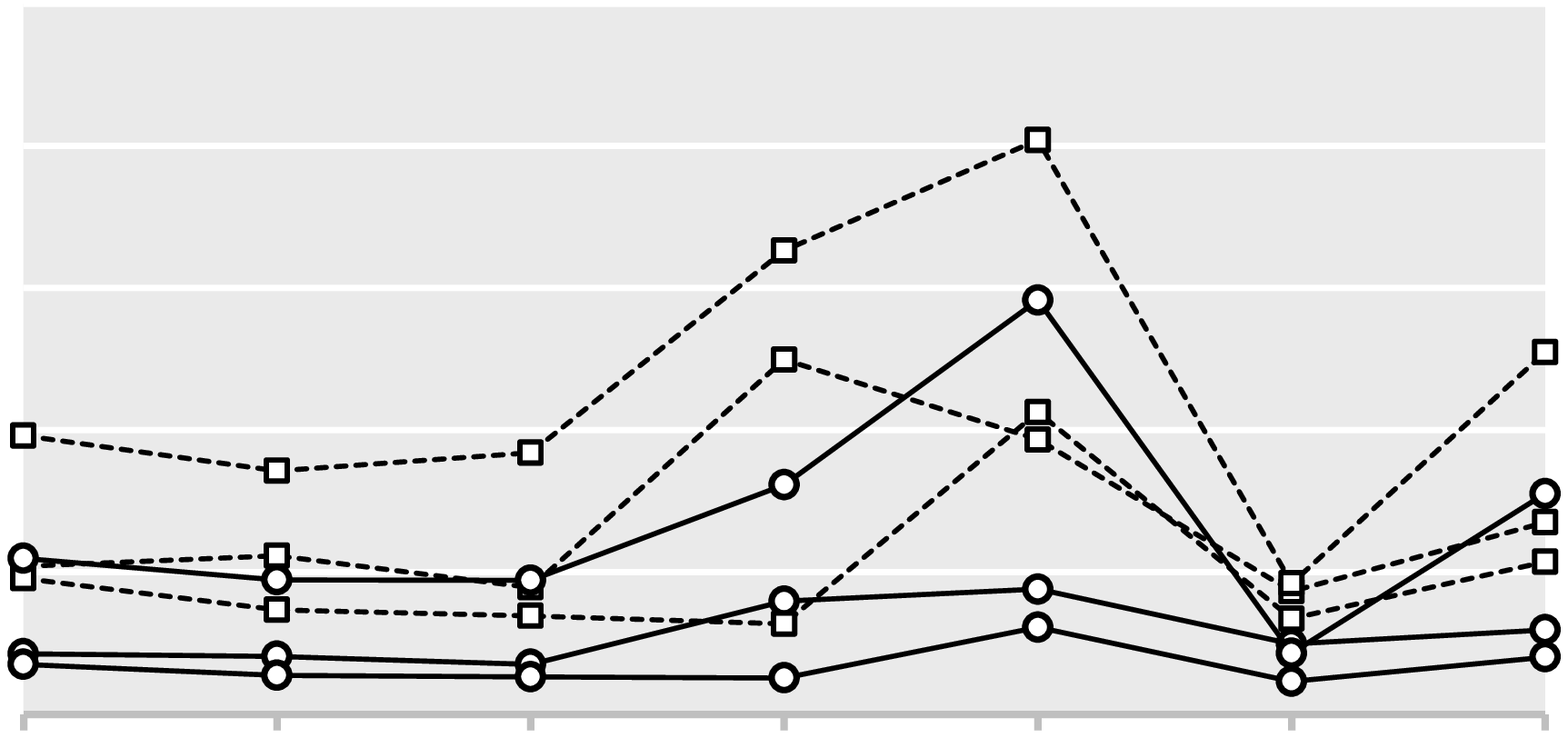}}
                \rput(0.225\linewidth,0.79){\footnotesize \textbf{\sf Network number}}
		\rput(0.225\linewidth,0.17){\footnotesize ($\textrm{b}$) Scenarios~$3$, $4$, and~$5$}

		\rput[bl](-0.08,5.02){\footnotesize \textbf{{\sf Communication steps}}}
                \rput[r](0.48,1.67){\scriptsize $0$}
                \rput[r](0.48,2.279){\scriptsize $400$}
                \rput[r](0.48,2.884){\scriptsize $800$}
                \rput[r](0.48,3.505){\scriptsize $1200$}
                \rput[r](0.48,4.09){\scriptsize $1600$}
                \rput[r](0.48,4.66){\scriptsize $2000$}

                \rput[t](0.604,1.42){\scriptsize $\mathsf{1}$}
                \rput[t](1.630,1.42){\scriptsize $\mathsf{2}$}
                \rput[t](2.673,1.42){\scriptsize $\mathsf{3}$}
                \rput[t](3.721,1.42){\scriptsize $\mathsf{4}$}
                \rput[t](4.764,1.42){\scriptsize $\mathsf{5}$}
                \rput[t](5.800,1.42){\scriptsize $\mathsf{6}$}
                \rput[t](6.84,1.42){\scriptsize $\mathsf{7}$}

                \rput[l](7,2.05){\tiny \textbf{\sf D-ADMM: 3}}
                \rput[l](7,2.5){\tiny \textbf{\sf D-Lasso: 3}}
                \rput[l](7,1.85){\tiny \textbf{\sf D-ADMM: 4}}
               	\rput[l](7,2.28){\tiny \textbf{\sf D-Lasso: 4}}
		\rput[l](7,2.69){\tiny \textbf{\sf D-ADMM: 5}}
		\rput[l](7,3.21){\tiny \textbf{\sf D-Lasso: 5}}
		
            \end{pspicture}
	  }{
            \begin{pspicture}(0.46\linewidth,5.4cm)
                \rput(0.214\linewidth,3.135){\includegraphics[scale=0.38]{figures/CPBestRhoCen345.eps}}
                \rput(0.21\linewidth,0.79){\footnotesize \textbf{\sf Network number}}
		\rput(0.21\linewidth,0.17){\footnotesize ($\textrm{b}$) Scenarios~$3$, $4$, and~$5$}

		\rput[bl](-0.08,5.02){\footnotesize \textbf{{\sf Communication steps}}}
                \rput[r](0.48,1.67){\scriptsize $0$}
                \rput[r](0.48,2.279){\scriptsize $400$}
                \rput[r](0.48,2.884){\scriptsize $800$}
                \rput[r](0.48,3.505){\scriptsize $1200$}
                \rput[r](0.48,4.09){\scriptsize $1600$}
                \rput[r](0.48,4.66){\scriptsize $2000$}

                \rput[t](0.652,1.34){\scriptsize $\mathsf{1}$}
                \rput[t](1.732,1.34){\scriptsize $\mathsf{2}$}
                \rput[t](2.803,1.34){\scriptsize $\mathsf{3}$}
                \rput[t](3.876,1.34){\scriptsize $\mathsf{4}$}
                \rput[t](4.960,1.34){\scriptsize $\mathsf{5}$}
                \rput[t](6.044,1.34){\scriptsize $\mathsf{6}$}
                \rput[t](7.111,1.34){\scriptsize $\mathsf{7}$}

                \rput[l](7.29,2.05){\scriptsize \textbf{\sf D-ADMM: 3}}
                \rput[l](7.29,2.5){\scriptsize \textbf{\sf D-Lasso: 3}}
                \rput[l](7.29,1.85){\scriptsize \textbf{\sf D-ADMM: 4}}
               	\rput[l](7.29,2.28){\scriptsize \textbf{\sf D-Lasso: 4}}
		\rput[l](7.29,2.69){\scriptsize \textbf{\sf D-ADMM: 5}}
		\rput[l](7.29,3.21){\scriptsize \textbf{\sf D-Lasso: 5}}
		
            \end{pspicture}
        }
        }
       \caption{Type II experiments for the column partition: number of communication steps to reach $10^{-3}\%$ of accuracy.}
       \label{Fig:CPType2Experiments}
	\end{figure*}

    Figure~\ref{Fig:CPType2Experiments} shows the plots with the results of the type II experiments.  As before, D-ADMM always required less communication steps to achieve a $10^{-3}\%$ of accuracy. In particular, D-ADMM used in average~$42\%$ of the communications D-Lasso used; the standard deviation was $10\%$. The largest difference in the number of communications occurred in scenario~$2$, network~$4$, where D-ADMM only used $28\%$ of the communications that D-Lasso used. The smallest difference was~$72\%$ and it occurred in scenario~$5$, network~$5$. We mention that, in contrast with the row partition, there were cases in which the optimal value for~$\rho$ was~$10^{-3}$ and~$10$ (cf. Table~\ref{Tab:TypeExperiments}), the ``boundary'' values of the set of $\rho$'s we used. Therefore, we might improve the results if we try a wider range of~$\rho$'s.
        
\section{Final Remarks and Conclusions}
\label{Sec:Conclusions}

    We proposed an algorithm for solving BP in two distributed frameworks. In one framework, the BP matrix is partitioned by rows, with its rows distributed over a network with an arbitrary number of nodes; in the other framework, it is the columns of the matrix that are distributed. The only requirement on the topology of the network through which the nodes communicate is connectivity (and we also assume that this topology does not change along the algorithm). Therefore, our algorithms can be applied to several scenarios, ranging from sensor networks, where the communication network is usually sparse, to super-computing platforms, characterized by dense networks.

    We simulated our algorithms for several types of data and networks and conclude that they always require less communications than competing algorithms. This is paramount in energy-constrained environments such as sensor networks.

\bibliographystyle{IEEEbib}

{ \isdraft{\singlespace}{}
\bibliography{paper}
}

\appendices

\section{Alternating Direction Method of Multipliers}
\label{App:ADMM}

    Let~$f$ and $g$ be two real-valued convex functions and~$X$ and~$Y$ two polyhedral sets. Let also~$A$ and~$B$ be two full column-rank matrices, and consider the problem
    \begin{equation}\label{Eq:BackOptimPromADMM}
        \begin{array}{ll}
          \underset{x \in X, y \in Y}{\textrm{minimize}} & f(x) + g(y) \\
          \textrm{subject to} & Ax + By = 0\,,
        \end{array}
    \end{equation}
    with variables~$x$ and~$y$. The \emph{alternating direction method of multipliers} (ADMM)~\cite{BoydADMM,Bertsekas:Parallel} solves~\eqref{Eq:BackOptimPromADMM} by applying the method of multipliers~\cite[p.408]{Bertsekas:Nonlinear} concatenated with one single loop of the nonlinear Gauss-Seidel~\cite[p.272]{Bertsekas:Nonlinear}: 
    \begin{eqnarray}
        x^{(k+1)} &=& \arg\min_{x \in X} f(x) + \phi_{\lambda^{(k)}}(Ax + By^{(k)})
        \label{Alg:ADMM1}
        \\
        y^{(k+1)} &=& \arg\min_{y \in Y} g(y) + \phi_{\lambda^{(k)}}(Ax^{(k+1)} + By)
        \label{Alg:ADMM2}
        \\
        \lambda^{(k+1)} &=& \lambda^{(k)} + \rho (Ax^{(k+1)} + By^{(k+1)})\,,
        \label{Alg:ADMM3}
    \end{eqnarray}
    where~$\phi_\lambda(z) = \lambda^\top z + \frac{\rho}{2}\|z\|^2$ and~$\rho > 0$ is a tuning parameter. In words, the augmented Lagrangian
    $$
        L(x,y;\lambda) = f(x) + g(y) + \lambda^\top (Ax + By) + \frac{\rho}{2}\|Ax + By\|^2\,,
    $$
    is first minimized with respect to~$x$ and then, keeping the value of~$x$ fixed at the just computed value~$x^{(k+1)}$, the augmented Lagrangian is minimized with respect to~$y$. Thus, \eqref{Alg:ADMM1} and~\eqref{Alg:ADMM2} cannot be carried out simultaneously. After these minimization steps, the dual variable~$\lambda$ is updated in a gradient-based way via~\eqref{Alg:ADMM3}. The following theorem guarantees its convergence.
    \begin{Theorem}[\cite{BoydADMM,Bertsekas:Parallel,Mota11-ADMMProof}]
    \label{Thm:ADMMConvergence}
				Let~$f:\mathbb{R}^{n_1} \xrightarrow{} \mathbb{R}$ and~$g:\mathbb{R}^{n_2} \xrightarrow{} \mathbb{R}$ be convex over~$\mathbb{R}^{n_1}$ and~$\mathbb{R}^{n_2}$, respectively. Let~$X \subset \mathbb{R}^{n_1}$ and~$X \subset \mathbb{R}^{n_2}$ be polyhedral sets and let~$A$ and~$B$ two full column-rank matrices. Assume~\eqref{Eq:BackOptimPromADMM} is solvable. Then,
        \begin{enumerate}
          \item $\{(x^{(k)}, y^{(k)})\}$ converges to a solution of~\eqref{Eq:BackOptimPromADMM};
          \item $\{\lambda^{(k)}\}$ converges to a solution of the dual problem
            $$
							\isdraft{
								\underset{\lambda}{\text{maximize}}\,\,\, F(\lambda) + G(\lambda)\,,
							}{
                \begin{array}{cl}
                  \textrm{\emph{maximize}} & F(\lambda) + G(\lambda) \\
                  \lambda &
                \end{array}\,,
							}
            $$
            where $F(\lambda) = \inf_{x \in X} f(x) + \lambda^\top Ax$ and $G(\lambda) = \inf_{y \in Y} g(y) + \lambda^\top By$\,.          
        \end{enumerate}
    \end{Theorem}
    Furthermore, \cite{He11-OnTheConvergenceRateADM} recently proved that ADMM converges with rate~$O(1/k)$. This rate holds even if the quadratic term of~$\phi_\lambda(z)$ in~\eqref{Alg:ADMM1} is linearized, which can many times simplify the solution of that optimization problem. For more properties of ADMM and its relation to other algorithms see \cite{ProximalSplittingInSP,Chambolle11-ADMM}.
    
    We now present a generalization of ADMM, which we call ``generalized ADMM.'' The generalized ADMM solves:
    \begin{equation}\label{Eq:AppGeneralizedADMM}
        \begin{array}{ll}
          \textrm{minimize} & \sum_{i=1}^I f_i(x_i)\\
          \textrm{subject to} & x_i \in X_i \,,\quad i=1,\ldots,I \\
															& \sum_{i=1}^I A_i x_i = 0\,,
        \end{array}
    \end{equation}
    where~$(x_1,\ldots,x_I)$ is the variable, $I>2$, the functions~$f_i$ are convex, $A_i$ are full column-rank matrices, and~$X_i$ are polyhedral sets. The generalized ADMM solves~\eqref{Eq:AppGeneralizedADMM} with:
    \begin{eqnarray*}
        x_1^{(k+1)} &=& \arg\min_{x_1 \in X_1} f_1(x_1) + \phi_{\lambda^{(k)}}(A_1x_1 + \sum_{j> 1}A_jx_j^{(k)})
        \\
        &\vdots&
        \\
        x_i^{(k+1)} &=& \arg\min_{x_i \in X_i} f_i(x_i) + \phi_{\lambda^{(k)}}(A_ix_i + \sum_{j < i}A_jx_j^{(k+1)}
        \isdraft{}{\\ && \phantom{ssssssssssssssssssssssssssss}}+ \sum_{j > i}A_jx_j^{(k)})
        \\
        &\vdots&
        \\
        x_I^{(k+1)} &=& \arg\min_{x_I \in X_I} f_I(x_I) + \phi_{\lambda^{(k)}}(A_Ix_I + \sum_{j < I}A_jx_j^{(k+1)})
        \\
        \lambda^{(k+1)} &=& \lambda^{(k)} + \rho \sum_{i=1}^I A_ix_i^{(k+1)}\,.
    \end{eqnarray*}
    This algorithm is then the natural generalization of~\eqref{Alg:ADMM1}-\eqref{Alg:ADMM3}. It is not known yet if Theorem~\ref{Thm:ADMMConvergence} also applies to the generalized ADMM. The latest efforts for doing that can be found in~\cite{He10ADMGaussianBack,He10SplittingMethodSeparateConvexProgramming,He10AlternatingdirectionsBasedContraction,He11-LinearizedADMGaussianBackSubstitution}.
    In spite of this fact, we apply the generalized ADMM in this paper and the resulting algorithm never failed to converge in our simulations.

\section{Problem for Each Node: Row Partition}
\label{App:OptimForEachNode}

    In the distributed algorithm we propose, each node has to solve, in each iteration, the problem
    \begin{equation}\label{Eq:CommonProblem}
        \begin{array}{ll}
          \textrm{minimize} & \|x\|_1 + v^\top x + c \|x\|^2 \\
          \textrm{subject to} & A x = b\,,
        \end{array}
    \end{equation}
    where~$x \in \mathbb{R}^n$ is the variable, and~$v \in \mathbb{R}^n$, $c >0$, $A \in \mathbb{R}^{m \times n}$, and $b \in \mathbb{R}^{m}$ are given. We propose to solve~\eqref{Eq:CommonProblem} by solving its dual problem:
    \begin{equation}\label{Eq:CommonProbDual}
			\isdraft{
						\underset{\lambda}{\textrm{maximize}} \,\,\,
						\lambda^\top b + \sum_{i=1}^n \inf_{x_i}
            \left(
                |x_i| + u_i(\lambda) x_i + c x_i^2
            \right)\,,
					}{
        \begin{array}{cl}
          \textrm{maximize} & \lambda^\top b + \sum_{i=1}^n \inf_{x_i}
            \left(
                |x_i| + u_i(\lambda) x_i + c x_i^2
            \right) \\
          \lambda &
        \end{array},
				}	
    \end{equation}
    where the dual variable is~$\lambda \in \mathbb{R}^{m}$ and~$u(\lambda) = v - A^\top \lambda$. To compute the objective of this dual problem for a fixed~$\lambda$, we need to find the minimum~$x_i(\lambda)$ of the function $|x_i| + u_i(\lambda) x_i + c x_i^2$ for~$i=1,\ldots,n$. Each one of these functions is strictly convex due to~$c > 0$, and hence it has a unique minimizer~$x_i(\lambda)$. It follows from Danskin's theorem~\cite[prop. B25]{Bertsekas:Nonlinear} that the objective of~\eqref{Eq:CommonProbDual} is differentiable with gradient~$b - A x(\lambda)$, where~$x(\lambda) = (x_1(\lambda),\ldots,x_n(\lambda))$. By the optimal conditions for convex problems~\cite[prop.B24]{Bertsekas:Nonlinear},
    $$
        x_i(\lambda) =
        \left\{
            \begin{array}{ll}
              0 &, -1\leq u_i(\lambda) \leq 1 \\
              -(u_i(\lambda)+1)/(2c) &, u_i(\lambda) < -1 \\
              -(u_i(\lambda) - 1)/(2c) &, u_i(\lambda) >1
            \end{array}
        \right.\,.
    $$
    The unicity of the minimizers~$x_i(\lambda)$ also implies that, once a solution~$\lambda^\star$ of~\eqref{Eq:CommonProbDual} is known, the solution of~\eqref{Eq:CommonProblem} is given by~$x(\lambda^\star)$. To solve~\eqref{Eq:CommonProbDual}, we propose using the method in~\cite{Barzilai-Borwein}, a very efficient Barzilai-Borwein (BB) algorithm. Per iteration, BB consumes~$O(n)$ flops plus the flops necessary to compute the gradient. Furthermore, BB is known to converge $R$-superlinearly for generic unconstrained optimization problems~\cite[Th.4]{BBNarushima}.
    
    As a final note, the number of iterations to solve~\eqref{Eq:CommonProblem} can be drastically reduced by using warm-starts. This means that, at iteration~$k+1$, node~$p$ will solve~\eqref{Eq:CommonProblem} by starting the BB algorithm with the solution found in iteration~$k$. The solutions of these two consecutive problems are expected to be close, since only~$v$ and~$c$ changed, possibly just by a small quantity.

\end{document}